\documentclass[11pt]{amsart}
\usepackage{amsmath,amsthm, amscd, amssymb, amsfonts, mathrsfs,pb-diagram,pb-xy}
\usepackage[all]{xy}
\usepackage{enumitem}
\usepackage{mathrsfs}
\usepackage{srcltx}
\usepackage{hyperref}

\allowdisplaybreaks

\usepackage{multicol}
\usepackage[dvips, dvipsnames, usenames]{color}

\newcommand{\superda}[1]{{\mathbf D}(2,1;#1)}
\newcommand{\superf}{{\mathbf F}(4)}
\newcommand{\superg}{{\mathbf G}(3)}

\newcommand{\superqa}[3]{{\bf A}_{#1}(#2|#3)}
\newcommand{\superqb}[3]{{\bf B}_{#1}(#2|#3)}

\newcommand{\superqd}[3]{{\bf D}_{#1}(#2|#3)}

\newcommand{\ad}{\operatorname{ad}}
\newcommand{\adc}{\operatorname{ad}_c}
\newcommand{\co}{\operatorname{co}}

\newcommand{\gr}{\operatorname{gr}}

\newcommand{\id}{\operatorname{id}} 
\newcommand{\ord}{\operatorname{ord}}
\newcommand{\ztu}{\overline{\zeta}}
\newcommand{\GK}{\operatorname{GKdim}}

\newcommand{\ku}{ \Bbbk}
\newcommand{\I}{\mathbb I}
\newcommand{\G}{\mathbb G}
\newcommand{\N}{\mathbb N}
\newcommand{\Z}{\mathbb Z}
\newcommand{\J}{\mathbb J}

\newcommand{\toba}{\mathscr{B}}
\newcommand{\wtoba}{\widetilde{\mathscr{B}}}
\newcommand{\htoba}{\widehat{\toba}}
\newcommand{\bq}{\mathfrak{q}}
\newcommand{\bp}{\mathfrak{p}}

\newcommand{\ydh}{{}^{H}_{H}\mathcal{YD}}
\newcommand{\ydkg}{{}^{\ku\Gamma}_{\ku\Gamma}\mathcal{YD}}
\newcommand{\ydkzt}{{}^{\ku\Z^{\theta}}_{\ku\Z^{\theta}}\mathcal{YD}}

\newcommand{\Ic}{\mathcal{I}}

\newcommand{\Pc}{\mathcal{P}}

\newcommand{\Zc}{\mathcal{Z}}
\newcommand{\Jc}{\mathcal{J}}

\newcommand{\Xc}{\mathcal{X}}
\newcommand{\Wc}{\mathcal{W}}

\newcommand{\fO}{\mathfrak O}

\newcommand{\pre}{\mathfrak{Pre}(V)}
\newcommand{\prefd}{\mathfrak{Pre}_{\operatorname{GKd}}(V)}

\newcommand{\qti}{\widetilde{q}}
\newcommand{\wbq}{\widetilde{\bq}}

\newcommand{\pf}{\begin{proof}}
\newcommand{\epf}{\end{proof}}

\newcommand{\GL}{\operatorname{GL}}

\numberwithin{equation}{section}
\theoremstyle{plain}

\newtheorem{theorem}{Theorem}[section]
\newtheorem{lemma}[theorem]{Lemma}
\newtheorem{definition-theorem}[theorem]{Definition-Theorem}
\newtheorem{prop}[theorem]{Proposition}

\theoremstyle{definition}

\newtheorem{remark}[theorem]{Remark}

\newtheorem{conjecture}[theorem]{Conjecture}

\newtheorem{step}{Step}
\newtheorem{stepo}{Step}

\newcommand{\ot}{\otimes}

\begin{document}

\noindent
\title[Pre-Nichols algebras of super and standard type]
{Finite GK-dimensional pre-Nichols algebras of super and standard type}

\author[Angiono]{Iv\'an Angiono}
\address{\noindent FaMAF-CIEM (CONICET) \\ 
	Universidad Nacional de C\'ordoba \\
	Medina Allende s/n, Ciudad Universitaria \\
	5000 C\'ordoba \\
	Rep\'ublica Argentina
}
\email{(angiono|gsanmarco)@famaf.unc.edu.ar}
\author[Campagnolo]{Emiliano Campagnolo}
\email{emiliano.campagnolo@mi.unc.edu.ar}
\author[Sanmarco]{Guillermo Sanmarco}

\keywords{Hopf algebras, Nichols algebras, Gelfand-Kirillov dimension.
\\
MSC2020: 16T05, 16T20, 17B37, 17B62.}

\thanks{The work of the three authors was partially supported by CONICET and Secyt (UNC)}

\begin{abstract}
We prove that finite GK-dimensional pre-Nichols algebras of super and standard type are quotients of the corresponding distinguished pre-Nichols algebras, except when the braiding matrix is of type super A and the dimension of the braided vector space is three. For these two exceptions we explicitly construct substitutes as braided central extensions of the corresponding pre-Nichols algebras by a polynomial ring in one variable. Via bosonization this gives new examples of finite GK-dimensional Hopf algebras.
\end{abstract} 

\maketitle

\section{Introduction}

Fix an algebraically closed field $\Bbbk$ of characteristic zero. We abbreviate $\GK A$ for the Gelfand-Kirillov dimension of an associative algebra $A$.

The study of Hopf algebras with finite $\GK$ can be approached from two different perspectives. A first possibility consists on looking for particular features of the underlying algebra structure. For instance noetherian Hopf algebras were intensively studied \cite{BG,G-survey} and many authors gave classification results under this or similar hypothesis, see for example
\cite{BZ,Liu,GZ,WZZ,G-survey}.

A second approach comes from the underlying coalgebra structure. Following the Lifting Method \cite{AS}, the first invariant to be considered is the coradical (the sum of all simple subcoalgebras). An intensively studied family is that of pointed Hopf algebras --- those whose coradical coincide with the group algebra of the so called group-like elements. If $H$ is a pointed Hopf algebra with group of group-like elements $\Gamma$ then the coradical filtration gives raise to a graded Hopf algebra $\gr H$ endowed with maps $\ku\Gamma \hookrightarrow \gr H \twoheadrightarrow \ku\Gamma$ composing to the identity of $\ku\Gamma$. Thus $\gr H$ decomposes as the \emph{bosonization} $\gr H\simeq R\#\Bbbk\Gamma$, where $R$ (the so called \emph{diagram} of $H$) is a coradically graded Hopf algebra in $\ydkg$ such that $R^0=\Bbbk1$. The degree $1$ component $V:=R^1\in\ydkg$ is the \emph{infinitesimal braiding} of $H$, the subalgebra $\toba(V)$ generated by $V$ is the Nichols algebra of $V$ and $R$ is a post-Nichols algebra of $V$, see \S \ref{subsec:Nichols}. 

Given a group $\Gamma$ with finite growth, in order to classify finite $\GK$ pointed Hopf algebras having $\Gamma$ as group of group-like elements the Lifting Method proposes: 
\begin{enumerate}[leftmargin=*,label=\rm{(\Roman*)}]
\item\label{item:lifting-I} classify all Yetter-Drinfeld modules $V$ over  $\Gamma$ such that $\GK\toba(V)<\infty$, 
\item\label{item:lifting-II} for each $V$ as in \ref{item:lifting-I}, determine all finite $\GK$ post-Nichols algebras $R$ of $V$;
\item\label{item:lifting-III} for all such $R$, obtain all possible $H$ such that $\gr H\simeq R\#\Bbbk\Gamma$.	
\end{enumerate}

Assume from now on that $\Gamma$ is abelian and finitely generated. Each object $V \in \ydkg$ either splits as direct sum of one-dimensional submodules or it contains an indecomposable submodule of dimension $>1$. In the former scenario we say that $V$ is of \emph{diagonal type}; the braiding is described by a matrix $\bq\in(\Bbbk^\times)^{\theta\times\theta}$, see \S \ref{subsec:diagonal}, and the Nichols algebra is denoted by $\toba_{\bq}$. Regarding \ref {item:lifting-I}, based on the evidence from \cite{AAH-diag}, the same authors proposed:

\begin{conjecture}{\cite[1.5]{AAH-memoirs}}\label{conj:AAH}
	If $\GK\toba_\bq<\infty$, then the root system $\varDelta^\bq$ is finite.
\end{conjecture}
\emph{Throughout this paper we will assume the validity of this statement}:
it means that $\bq$ belongs to the classification given in \cite{H-classif}. About non-diagonal braided vector spaces, a wide classification result was performed in \cite{AAH-memoirs}.

\medbreak
For \ref{item:lifting-II} the situation is quite different from the finite-dimensional context, where the \emph{generation in degree-one problem} has a positive answer \cite{An-crelle}: it means that the unique finite dimensional post-Nichols algebra of each $\bq$ is $\toba_\bq$ itself. This is not longer true when we move to finite $\GK$. For example, if $\dim\toba_\bq<\infty$, the so-called \emph{distinguished} pre-Nichols algebra $\wtoba_\bq$ is a finite $\GK$ pre-Nichols algebra different from $\toba_\bq$, see \S \ref{subsec:distinguished}.

\medbreak

As it was discussed in \cite[\S 2.6]{ASa}, the question in \ref{item:lifting-II} translates to
\begin{enumerate}[leftmargin=*]
\item[(II')]\label{item:lifting-II'} For each $\bq$ as in \ref{item:lifting-I}, determine all pre-Nichols algebras $\toba$ of $V$ with $\GK \toba<\infty$.
\end{enumerate}

Assume that $\bq$ is as in \ref{item:lifting-I}.
Pre-Nichols algebras of a fixed $\bq$ form a poset with minimal element $T(V)$ and maximal one $\toba_\bq$. Those with finite $\GK$ determine a sub-poset and a first question is if there exist minimal elements (called \emph{eminent} pre-Nichols algebras). If $\dim\toba_\bq<\infty$, then it is natural to ask whether the distinguished pre-Nichols algebra is eminent, as it was formulated in \cite{An-dist}.

Following \cite{AA17}, the list given in \cite{H-classif} containing all matrices $\bq$ with connected Dynkin diagram and finite generalized root system can be organized as follows:

\begin{itemize} [leftmargin=*]\renewcommand{\labelitemi}{$\diamond$}
	\smallbreak\item  Standard type: the Cartan matrix is constant in the Weyl equivalence class. This family includes Cartan type \cite{AS}, but also braidings with Cartan matrices $B_\theta$ and $G_2$.
	
	\smallbreak\item Super type: the Weyl groupoid \emph{coincides} with the one of a finite-dimensional contragredient Lie superalgebras in characteristic 0.
	
	\smallbreak\item (Super) modular type: related to finite-dimensional contragredient Lie (super)algebras in positive characteristic.
	
	\smallbreak\item \textit{UFO} type: a list of twelve Weyl equivalence classes of braiding matrices.
\end{itemize}

Eminent pre-Nichols algebras of Cartan type were considered in \cite{ASa}. 
Here we give another step towards the classification of pointed Hopf algebras  with finite $\GK$ and abelian group of group-likes: we answer the question above for some exceptional cases of Cartan type $G_2$, braidings of standard type (but not Cartan) and the whole family of super type. Super modular and UFO type will be considered in a forthcoming paper. 

As in \cite{ASa} we find that the distinguished pre-Nichols algebra is eminent for almost all braidings considered here and determine the eminent pre-Nichols algebras in these exceptional cases. Coming back to \ref{item:lifting-II}, this means that in almost all the cases, post-Nichols algebras with finite $\GK$ are subalgebras of the Lusztig algebras introduced in \cite{AAR}, containing the corresponding Nichols algebras.

More explicitly, the main result of this paper can be summarized as follows.

\begin{theorem}\label{thm:main}
\begin{enumerate}[leftmargin=*]
\item The distinguished pre-Nichols algebra $\wtoba_\bq$ is eminent when $\bq$ is:
\smallbreak

\begin{itemize}
\item of Cartan type $G_2$ with parameter $q\in\G_4'\cup\G_6'$.
\smallbreak

\item of super type, except $\superqa{3}{q}{\J}$, $q\in\G_{\infty}$, $\J=\{2\}$ or $\J=\{1,2,3\}$.
\smallbreak

\item of standard type $B_{\theta}$ and $G_2$.
\smallbreak

\end{itemize}
\item Suppose that $\bq$ is of type $\superqa{3}{q}{\{2\}}$, $q\in\G_{\infty}$. Then
\begin{align*}
\htoba_{\bq}=T(V)/\langle x_2^2,x_{13},x_{112},x_{332} \rangle
\end{align*}
is an eminent pre-Nichols algebra of $\bq$, and $\GK \htoba_{\bq}= 3$.

\item Suppose that $\bq$ is of type $\superqa{3}{q}{\{1,2,3\}}$, $q\in\G_{\infty}$. Then
\begin{align*}
\htoba_{\bq}=T(V)/\langle x_1^2, x_2^2, x_3^2, x_{213}, [x_{123},x_2]_c \rangle
\end{align*}
is an eminent pre-Nichols algebra of $\bq$, and $\GK \htoba_{\bq}= 3$.
\end{enumerate}
	
\end{theorem}

The structure of the paper is the following. In \S \ref{sec:preliminaries} we recall the definitions and basic results on Nichols algebras, pre and post-Nichols algebras, with emphasis on the diagonal case; we also give a relation between Hilbert series and extensions of braided Hopf algebras.

Distinguished pre-Nichols algebras are defined in terms of the specific presentation by generators and relations of Nichols algebras  given in \cite[Theorem 3.1]{An-crelle}. In \S \ref{sec:def-relations} we determine sufficient conditions for some of these relations to hold in any finite $\GK$ pre-Nichols algebra. This is treated within a general setting that also includes defining relations for braidings of (super) modular or UFO type.

The proof of Theorem \ref{thm:main} for Cartan type $G_2$ is contained in \S \ref{section:Cartan-typeG2}, where first we determine a minimal presentation when the parameter is a root of unity of order $4,6$. The proof for super type braidings is given in \S \ref{section:super-type} with a case-by-case analysis. For the eminent pre-Nichols algebras which are not the distinguished ones we give a presentation by generators and relations and a PBW basis; moreover we show that 
$\htoba_{\bq}$ contains a central Hopf subalgebra $\widehat\Zc_{\bq}$ that fits in an extension of braided Hopf algebras $\Bbbk \to \widehat\Zc_{\bq} \hookrightarrow \htoba_{\bq} \twoheadrightarrow \wtoba_{\bq}\to\Bbbk$. Finally, \S \ref{section:standard-type} is devoted to standard type, where we also give a minimal presentation for $G_2$ type.

\subsection*{Notation}
For a positive integer $\theta$ we set $\I_{\theta}=\{1,\dots ,\theta\}$; 
if $\theta$ is understood we shall write simply $\I$. 
The canonical basis of $\Z^{\I}$ is denoted $\{\alpha_i: i\in\I\}$.

Let $N\in\N$. The subgroup of $\ku^{\times}$ of $N$-th roots of unity is denoted by $\G_N$, and $\G_ {N}'$ denotes the subset of those of order $N$. We also set $\G_{\infty}=\bigcup_{N \in \N} {\G_{N}}$.

We always consider Hopf algebras over $\Bbbk$ and assume they have bijective antipode. 
The subspace of primitive elements of a (braided) Hopf algebra $H$ is $\Pc(H)$. The category of (left) Yetter-Drinfeld modules over $H$ is denoted by $\ydh$.
The subalgebra generated by a subset $X\subseteq H$ is denoted by $\Bbbk\langle X\rangle$. 

We refer to \cite{Rad-libro} for any unexplained notion on Hopf algebras and to \cite{KL} for the definition and basic properties of Gelfand-Kirillov dimension of associative algebras.

\section{Preliminaries}\label{sec:preliminaries}

\subsection{Nichols algebras}\label{subsec:Nichols}
A braided vector space is a pair $(V,c)$, where $V$ is a vector space and the \emph{braiding} $c\in\GL(V\ot V)$  satisfies the Yang-Baxter equation:
\begin{align}\label{eq-YB}
 (c \ot \id)(\id \ot c)(c \ot \id)=(\id \ot c)(c \ot \id)(\id \ot c).
\end{align}

The tensor algebra $T(V)$ is a braided Hopf algebra in the sense of \cite{Tk1};
a \emph{pre-Nichols algebra} of $V$ is a quotient of $T(V)$ 
by a braided graded Hopf ideal contained in $\oplus_{n\geqslant 2}V^{\otimes n}$. 

Dually, the tensor coalgebra $T^c(V)$ 
is a braided Hopf algebra; 
a \emph{post-Nichols algebra} of $V$ is a graded subalgebra of $T^c(V)$ containing $V$.

By the universal properties of $T(V)$ and $T^c(V)$, there exists a (unique) graded Hopf algebra map $\Omega:T(V)\to T^c(V)$ such that $\Omega_{|V}=\id_V$.
The \emph{Nichols algebra} of $V$ is $\Omega(T(V))$; thus, there exists a braided graded Hopf ideal $\Jc(V)$ contained in $\oplus_{n\geqslant 2}V^{\otimes n}$ such that 
$\toba(V) = T(V) /\Jc(V)$. Moreover, $\Jc(V)$ is the maximal Hopf ideal among those.

\medbreak

Let $H$ be a Hopf algebra and $V \in \ydh$. As $\ydh$ is a braided tensor category, the braiding $c:V \ot V\rightarrow V\ot V $ of $\ydh$ makes $(V,c)$ a braided vector space. Hence we may consider $\toba(V)$, which is a Hopf algebra in $\ydh$.

\subsection{Braidings of diagonal type}\label{subsec:diagonal}

Let $\bq=(q_{ij})_{1 \leq i,j \leq \theta}\in (\Bbbk^{\times})^{\I\times\I}$. Consider the pair $(V,c^\bq)$, where $V$ has basis $(x_i)_{i \in\I}$ and $c^{\bq}\in \GL(V\ot V)$ is determined by:
\begin{align}\label{eq:diag-type-def}
c^{\bq}(x_i\ot x_j)&=q_{ij}x_j\ot x_i, & i,j&\in\I.
\end{align}
Then $c^{\bq}$ satisfies (\ref{eq-YB}). 
The braided vector space $(V,c^{\bq})$ and the associated Nichols algebra $\toba_{\bq}:=\toba(V)$ are said \emph{of diagonal type}. 
The generalized Dynkin diagram of $\bq$ is a graph with labelled vertices and edges, defined as follows:
\begin{itemize}[leftmargin=*]
\item the set of vertices is $\I$; the vertex $i$ is labelled with $q_{ii}$.
\item two vertices $i\ne j\in \I$ are connected by an edge if and only if $\qti_{ij}:=q_{ij}q_{ji}\neq 1$. If so, the edge is labelled with $\qti_{ij}$.
\end{itemize}

\smallbreak
We also denote by $\bq$ the $\Z$-bilinear form $\Z^{\I}\times\Z^{\I}\to \Bbbk^\times$
associated to the matrix $\bq$, that is $\bq(\alpha_j,\alpha_k) :=q_{jk}$, $j,k \in\I$.
If $\alpha,\beta  \in \Z^{\I}$ and $i\in \I$, then we set
\begin{align}\label{eq:notation-qab}
	q_{\alpha\beta} &= \bq(\alpha,\beta), & q_{\alpha} &= \bq(\alpha,\alpha),&
	N_{\alpha} &= \ord q_{\alpha}, & N_{i} &= \ord q_{\alpha_i} = N_{\alpha_i}.
\end{align}

\medbreak
The braided vector space $(V,c^{\bq})$ admits a realization over $\ydkzt$, where 
\begin{itemize}
	\item the coaction is given by $\delta(x_i)=\alpha_i\otimes x_i$, $i\in\I$;
	\item the action is defined by $\beta\cdot x_i=\bq(\beta,\alpha_i)x_i$, $i\in\I$, $\beta\in\Z^{\theta}$.	
\end{itemize}
In this context, $T(V)$ is a Hopf algebra in $\ydkzt$, $\Jc(V)$ is a Yetter-Drinfeld submodule of $T(V)$ and $\toba_{\bq}$ is a Hopf algebra in $\ydkzt$.
In particular, $T(V)$ is a $\Z^{\I}$-graded braided Hopf algebra, with grading $\deg x_i = \alpha_i$, $i\in \I$, and $\toba_{\bq}$ inherits the graduation. 

\medspace

Following \cite{Kh}, any $\Z^{\I}$-graded pre-Nichols algebra $\toba$ (and particularly $\toba_{\bq}$) has a PBW basis with homogeneous PBW generators. This means that there exists a  subset $\emptyset \neq S \subset \toba$ of homogeneous elements (the \emph{PBW-generators}) provided with a total order $<$, and a
function $h: S \to \N \cup \{ \infty \}$ (the \emph{height}) such that the following set is a $\Bbbk$-basis of $\toba$:
\begin{align*}
B(S,<,h):= \big\{s_1^{e_1}\dots s_t^{e_t}: \, t \in \N_0,\, s_i \in S, \, s_1>\dots >s_t, \, 0<e_i<h(s_i) \big\}.
\end{align*}

\subsection{Adjoint action and braided bracket}
Any braided Hopf algebra $R$ admits a \emph{left adjoint representation} $\adc: R\to\operatorname{End} R$,
\begin{align*}
	(\adc x)y&=m(m \ot S)(\id\ot c)(\Delta \ot \id)(x \ot y), && x,y\in R.
\end{align*}
Also, the \emph{braided bracket} $[\cdot, \cdot]_c:R \ot R \to R$ is the map given by 
\begin{align*}
	[x,y]_c &=m(\id-c) (x\ot y), && x,y\in R.
\end{align*}
Notice that $(\adc x)y=[x,y]_c$ if $x\in\Pc(R)$.

\medspace

We are interested in $\Z^{\I}$-graded pre-Nichols $R$ algebras of $V$. In this case the braided commutator satisfies
\begin{align}
	\label{eq:braided-commutator-right-mult}
	[u,vw]_c &= [u,v]_c w + \bq_{\alpha \beta} v [u,w]_c,
	\\
	\label{eq:braided-commutator-left-mult} 
	[uv,w]_c &= \bq_{ \beta \gamma}[u,w]_c v + u[v,w]_c,
	\\
	\label{eq:braided-commutator-iteration}
	\big[ [u,v]_c, w \big]_c &= \big[ u, [v, w]_c \big]_c 
	- \bq_{\alpha \beta} v [u,w]_c 
	+ \bq_{ \beta \gamma}[u,w]_c v,
\end{align}
for all homogeneous elements $u\in R^{\alpha}$,  $v \in R^{\beta}$, $w \in R^{\gamma}$. 
Given $i_1,\cdots, i_k\in\I$, $j\ge 2$, we set 
\begin{align*}
x_{i_1\cdots i_k}:=(\adc x_{i_1})x_{i_2\cdots i_k}=x_{i_1}x_{i_2\cdots i_k}-q_{i_1i_2}\cdots q_{i_1i_k}x_{i_2\cdots i_k}x_{i_1}.
\end{align*}

\subsection{Weyl groupoids and Cartan roots}
Next we briefly recall the notions of Weyl groupoid and generalized root systems \cite{H-Weyl gpd,HY-groupoid}.
We assume here that $\GK\toba_{\bq}<\infty$.

Let  $C^{\bq}=(c_{ij}^{\bq})_{i,j\in \I}\in\Z^{\I\times\I}$ be the generalized Cartan matrix
defined by $c_{ii}^{\bq} := 2$ and
\begin{align}\label{eq:defcij}
	c_{ij}^{\bq}&:= -\min \left\{ n \in \N_0: (n+1)_{q_{ii}}
	(1-q_{ii}^n q_{ij}q_{ji} )=0 \right\},  & i & \neq j.
\end{align}
When $\bq$ is fixed we simply write $(c_{ij})$.
Let $i\in \I$. 
The reflection $s_i^{\bq}\in GL(\Z^\I)$ is given by 
\begin{align}\label{eq:siq-definition}
	s_i^{\bq}(\alpha_j)&:=\alpha_j-c_{ij}^{\bq}\alpha_i, & &j\in \I.
\end{align}
Next we introduce the matrix  $\rho_i(\bq)$, given by
\begin{align}\label{eq:rhoiq-definition}
	(\rho_i (\bq))_{jk}&:= \bq(s_i^{\bq}(\alpha_j),s_i^{\bq}(\alpha_k)) = q_{jk} q_{ik}^{-c_{ij}^{\bq}} q_{ji}^{-c_{ik}^{\bq}} q_{ii}^{c_{ij}^{\bq}c_{ik}^{\bq}}, & &j, k \in \I,
\end{align}
and $\rho_i(V)$ is the braided vector space of diagonal type with matrix  $\rho_i (\bq)$. Finally we set
\begin{align*}
	\Xc := \{\rho_{j_1} \dots \rho_{j_n}(\bq): j_1, \dots, j_n \in \I, n \in \N \}.
\end{align*}
This set is called the Weyl-equivalence class of $\bq$.

\medbreak
For each $\bp\in\Xc$ the set $\varDelta_+^{\bp}$ of \emph{positive roots} consists of the $\Z^{\I}$-degrees of the generators of a PBW-basis of $\toba_\bp$, 
counted with multiplicities.
Let $\varDelta^{\bp} := \varDelta_+^{\bp} \cup -\varDelta_+^{\bp}$.
The  generalized root system of $\bq$ is the fibration $\Delta=(\Delta^{\bp})_{\bp\in\Xc}$ over $\Xc$.

\medbreak
The reflections $s_i^{\bp}$, $i\in \I$ and $\bp \in \Xc$, generate a subgroupoid of $\Xc \times GL(\Z^{\theta}) \times \Xc$ called
the Weyl groupoid $\Wc$ of $\bq$. The groupoid
$\Wc$ acts on $(C^{\bp})_{\bp \in \Xc}$ and on the generalized root system $(\varDelta^{\bp})_{\bp \in \Xc}$, generalizing the classical Weyl group.

\medspace
Next we assume that $\varDelta_+^{\bq}$ is finite.
Let $\omega_0^{\bq} \in \Wc$ be an element of maximal length and $\omega_0^{\bq}=\sigma_{i_1}^{\bq} \sigma_{i_2}\cdots \sigma_{i_\ell}$ 
be a reduced expression. Then
\begin{align} \label{eq:betak}
	\beta_k &:= s_{i_1}^{\bq}\cdots s_{i_{k-1}}(\alpha_{i_k}), & 
	k\in\I_{\ell}
\end{align}
are pairwise different; moreover $\Delta_+^{\bq}=\{\beta_k : k\in\I_{\ell}\}$ \cite{CH-rank 3}, so $\vert \Delta_+^{\bq} \vert = \ell$. 
In this setting, there exist PBW generators $x_{\beta}$ of degree $\beta$ such that the following set is a basis of $\toba_{\bq}$:
\begin{align*}
\{ x_{\beta_1}^{n_1} \cdots x_{\beta_{\ell}}^{n_{\ell}}: \, 0\le n_k<N_{\beta_k} \}.
\end{align*}

\subsection{Distinguished and eminent pre-Nichols algebras}\label{subsec:distinguished}

We start with the definition of Cartan roots \cite{An-crelle}. An element $i\in\I$ is a \emph{Cartan vertex} of $\bq$ if
\begin{align}\label{eq:cartan-vertex}
	q_{ij}q_{ji} &= q_{ii}^{c_{ij}^{\bq}}, & \text{for all } j \neq i.
\end{align}
The set of \emph{Cartan roots} of $\bq$ is the orbit of Cartan vertices under the action of $\Wc$; explicitly,
\begin{align*}
	\fO^{\bq} &= \{s_{i_1}^{\bq} s_{i_2} \dots s_{i_k}(\alpha_i) \in \Delta^{\bq}:
	i\in \I  \text{ is a Cartan vertex of } \rho_{i_k} \dots \rho_{i_2}\rho_{i_1}(\bq) \}.
\end{align*}

\medbreak

Now we assume that $\dim \toba_\bq<\infty$ and recall the definition of the distinguished pre-Nichols algebra introduced in \cite{An-crelle,An-dist}.
The presentation of $\toba_\bq$ given in \cite[Theorem 3.1]{An-crelle} includes a long list of relations in two, three and four letters $x_i$, and powers of root vectors $x_\beta^{N_\beta}$ for $\beta\in\fO^{\bq}_+$. Let $\Ic(V)$ denote the ideal of $T(V)$  generated by all the defining relations of $\toba_{\bq}$ except $x_\beta^{N_\beta}$, $\beta\in\fO^{\bq}_+$,
and adding the quantum Serre relations $(\adc x_i)^{1-c_{ij}}(x_j)$ when $ i\neq j$ are such that $q_{ii}^{c_{ij}}=q_{ij}q_{ji}=q_{ii}$. This happens to be a Hopf ideal.
The \emph{distinguished pre-Nichols algebra} is the braided Hopf algebra $\wtoba_{\bq}:=T(V)/\Ic(V)$.

When $\bq$ is of Cartan type, $\wtoba_{\bq}$ is the positive part of a multiparametric version of the de Concini-Procesi quantum group $U_q(\mathfrak g)$ (under restrictions on the order of $q$).

\medbreak

We recall some properties of $\wtoba_{\bq}$:
\begin{enumerate}[leftmargin=*,label=\rm{(\Alph*)}]
\smallbreak

\item Set $\fO^{\bq}_+ = \fO^{\bq} \cap \N_{0}^{\theta}$. 
Define $\widetilde N_\beta :=  N_{\beta}$ if $\beta\notin\fO^{\bq}_+$, and 
$\widetilde N_\beta := \infty$ if $\beta\in\fO^{\bq}_+$ (see \eqref{eq:notation-qab}).  
Let $x_{\beta}$ denote the canonical preimage in $\wtoba_\bq$ of the PBW generator of degree $\beta$ for $\toba_{\bq}$.
The set $\{ x_{\beta_1}^{n_1} \cdots x_{\beta_{\ell}}^{n_{\ell}}: \, 0\le n_k<\widetilde N_{\beta_k} \}$ is a basis of $\wtoba_{\bq}$.

\smallbreak

\item Let $\Zc_{\bq}$ be the subalgebra generated by $x_\beta^{N_\beta}$, $\beta\in\fO^{\bq}_+$. Then $\Zc_{\bq}$ is a $q$-polynomial ring in variables $x_\beta^{N_\beta}$, $\beta\in\fO^{\bq}_+$, which is also a Hopf subalgebra.
\smallbreak

\item We have that $\GK \wtoba_{\bq}= |\fO^{\bq}_+|$.
\end{enumerate}

\medbreak

Let $\bq$ be such that $\GK \toba_{\bq}< \infty$. The set of all pre-Nichols algebras of $\bq$ is a poset $\pre$ with $T(V)$ minimal and $\toba_{\bq}$ maximal. Let $\prefd$ be the subposet of $\pre$ of all finite GK-dimensional pre-Nichols algebras. 
We say that a pre-Nichols algebra $\htoba$ is \emph{eminent} if it is a minimum in $\prefd$; that is, for any $\toba\in\prefd$,there is an epimorphism of braided Hopf algebras which is $\id_V$ when we restrict to degree one.

\begin{remark}\label{rem:subalg-pre-Nichols}
Let $\toba$ be a pre-Nichols algebra of $\bq$. Fix $\J\subseteq \I$ and $\bq'=(q_{j,k})_{j,k\in\J}$, $V'$ the braided vector subspace of $V$ with basis $x_j$, $j\in\J$. The subalgebra $\toba'$ of $\toba$ generated by $x_j$, $j\in\J$, is a pre-Nichols algebra of $\bq'$, and $\GK\toba'\le \GK\toba$.

In particular, if $\toba\in\prefd$, then $\toba'\in \mathfrak{Pre}_{\operatorname{GKd}}(V')$.
\end{remark}

\begin{remark}\label{rem:dist-pre-Nichols}
Since $\wtoba_{\bq}\in\prefd$, in order to prove that $\wtoba_{\bq}$ is eminent it suffices to show that each defining relation of $\wtoba_{\bq}$ (of the presentation fixed above) holds in any $\toba\in\prefd$.
\end{remark}

\subsection{Extensions of graded braided Hopf algebras}\label{subsec:extensions}

Let $H$ be a Hopf algebra with bijective antipode. Recall that a sequence of morphisms of Hopf algebras in $\ydh$
\begin{align}\label{eq:def-exseq}
\ku \rightarrow A \overset{\iota}{\to} C \overset{\pi}{\to} B \rightarrow \ku
\end{align}
is an \emph{extension of braided Hopf algebras} (cf. \cite[\S 2.5]{AN}) if
\begin{enumerate}[leftmargin=*,label=\rm{(\alph*)}]
\begin{multicols}{2}
\item $\iota$ is injective,
\item $\pi$ is surjective,
\item $\ker \pi = C\iota(A^+)$ and
\item $A=C^{\,\co\pi}$, or equivalently $A=\,^{\co\pi}C$.
\end{multicols}
\end{enumerate}
For simplicity, we shall write $A \overset{\iota}{\hookrightarrow} C \overset{\pi}{\twoheadrightarrow} B$ instead of \eqref{eq:def-exseq}.

Assume further that $C$ is connected (i.e. the coradical of $C$ is $\Bbbk$), so for any extension $A \overset{\iota}{\hookrightarrow} C \overset{\pi}{\twoheadrightarrow} B$, also $B$ and $A$ are connected. In this context we recall the following result: 

\begin{prop}\label{prop:equiv-coinv}\cite[3.6]{A+} Let $C\in \ydh$ be a connected Hopf algebra. The map 
\begin{align*}
&\{\text{right coideal subalgebras of }C\}\to\{\text{quotient left }C-\text{module coalgebras}\}, &
A&\mapsto C/CA^{+},
\end{align*}
is bijective, with inverse $B \longmapsto C^{\operatorname{co}B}$.

If $A$ is a right coideal subalgebras of $C$ and $B=C/CA^{+}$, 
then there exists a left $B$-colinear and right $A$-linear isomorphism
$B\ot A \overset{\simeq}{\longrightarrow} C$.\qed
\end{prop}

Hence in order to get extensions of a given connected $C$, it is enough to consider either
\begin{itemize}[leftmargin=*]
\item a surjective Hopf algebra morphism $C \overset{\pi}{\twoheadrightarrow} B$ and set $A=C^{\,\co\pi}$, or 
\item a normal Hopf subalgebra $A$ of $C$, $\iota$ the inclusion and set $B=CA^+$.
\end{itemize}
In either case there exists a left $B$-colinear and right $A$-linear isomorphism
$ B\ot A \simeq C$.

An example of extensions of connected braided Hopf algebras is $\Bbbk \to \Zc_{\bq} \hookrightarrow \wtoba_{\bq} \twoheadrightarrow \toba_{\bq}\to\Bbbk$,where $\wtoba_{\bq}$ and $\Zc_{\bq}$ are defined as in the previous subsection.

\medbreak
It is more suitable to our purposes to consider $\N_0^{\theta}$-graded versions of these objects. 
For each $\alpha=(a_1,\cdots, a_{\theta})\in\Z^{\theta}$, set $t^{\alpha}=t_1^{a_1}\cdots t_{\theta}^{a_{\theta}}$.
Given a $\N_0^{\theta}$-graded object $U$ with finite-dimensional homogeneous components, the Hilbert (or Poincar\'e) series is
\begin{align*}
\mathcal{H}_{U} = \sum_{\alpha\in\N_0^{\theta}} \dim U_{\alpha} \, t^{\alpha}
\in \N_0[[t_1,\dots,t_{\theta}]].
\end{align*}

For example, from the definition of $\varDelta^{\bq}$ and the PBW bases of $\toba_{\bq}$, $\wtoba_{\bq}$, we have that
\begin{align*}
\mathcal{H}_{\toba_{\bq}} &=\prod_{\alpha\in\varDelta^{\bq}_+} \frac{1-t^{N_{\alpha}\alpha}}{1-t^{\alpha}},
&
\mathcal{H}_{\wtoba_{\bq}} &=\prod_{\alpha\in\fO^{\bq}_+} \frac{1}{1-t^{\alpha}}
\prod_{\alpha\in\varDelta^{\bq}_+-\fO^{\bq}_+} \frac{1-t^{N_{\alpha}\alpha}}{1-t^{\alpha}}.
\end{align*}

If $U'$ is another $\N_0^{\theta}$-graded object, then we say that $\mathcal{H}_U\le \mathcal{H}_{U'}$ if 
\begin{align*}
	\dim U_{\alpha} &\le \dim U'_{\alpha} & \text{for all }& \alpha\in\N_0^{\theta}.
\end{align*}

Next we relate extensions of (braided) connected Hopf algebras and Hilbert series.

\begin{lemma}\label{lem:extension-braided-graded}
Assume there is a degree-preserving extension $A \overset{\iota}{\hookrightarrow} C \overset{\pi}{\twoheadrightarrow} B$  of $\N_0^{\theta}$-graded connected Hopf algebras in $\ydh$ with finite-dimensional homogeneous components. Then $\mathcal{H}_C= \mathcal{H}_A\mathcal{H}_B$.
\end{lemma}
\pf
Consider the Hopf algebra $\widehat{H}=\Bbbk \Z^{\theta}\ot H$. Then every $\N_0^{\theta}$-graded object $U\in\ydh$ is canonically a Yetter-Drinfeld module over $\widehat{H}$, where
\begin{itemize}[leftmargin=*]
\item the action of $H$ on $U$ is extended to an action of $\widehat{H}$, where $\Z^{\theta}$ acts trivially;
\item for each $\beta\in\N_0^{\theta}$ and $u\in U_{\beta}$, the coaction on $u$ is given by
$$ \delta(u)=\beta\ot u_{-1}\ot u_0\in \N_0^{\theta} \ot H\ot U \subset \widehat{H}\ot U. $$
\end{itemize}

Since $\iota$ and $\pi$ preserve the degrees, $A \overset{\iota}{\hookrightarrow} C \overset{\pi}{\twoheadrightarrow} B$ is an extension in the category of $\widehat{H}$-Yetter-Drinfeld modules and \cite[Proposition 3.6 (d)]{A+} applies: there is an $\widehat{H}$-colinear, isomorphism  $B\ot A \overset{\simeq}{\longrightarrow} C$, which implies that 
$\mathcal{H}_C= \mathcal{H}_A\mathcal{H}_B$.
\epf

\subsection{Tool box for pre Nichols algebras of diagonal type}\label{subsec:useful-results}
We collect here some technical results which help to conclude that a given pre-Nichols algebra has infinite $\GK$.

\begin{remark}\label{lem:q-q2-(-q)}
If $W$ is of diagonal type and has Dynkin diagram $\xymatrix@C=30pt{\overset{q}{\circ} \ar@{-}[r]^{q^2}\ &\overset{-q}{\circ}}$ where $q\in \ku^\times$ is a root of unity of order larger than $2$, then $\GK \toba(W)=\infty$.

Indeed, one can see that this diagram is of Cartan type with non-finite associated Cartan matrix, so the rank-two result of \cite{AAH-diag} applies.
\end{remark}

\begin{lemma}\label{lem:1connected} \cite[Proposition 4.16]{AAH-memoirs}
Let $W$ be a braided vector space of diagonal type with diagram $\xymatrix@C=30pt{\overset{1}{\circ} \ar@{-}[r]^{p}\ &\overset{q}{\circ}}$, $p\ne 1$. Then $\GK \toba(W)=\infty$.
\end{lemma}

Next we state a well-known result: for a proof we refer to \cite[Lemma 2.8]{ASa}.

\begin{lemma}\label{lem:subspace-primitives}
Let $R$ be a graded braided Hopf algebra and $W$ a braided subspace of $\Pc(R)$. Then $\GK \toba(W) \le \GK R$.\qed
\end{lemma}

As we assume Conjecture \ref{conj:AAH}, it is useful to describe the shape of some Dynkin diagrams with finite root systems.

\begin{lemma}\label{lem:rank3}\cite[Lemma 9 (ii)]{H-classif} 
Assume $\bq=(q_{ij})_{i,j\in\I_3}$ has connected Dynkin diagram and finite root system. Then $\qti_{12}\qti_{13}\qti_{23}=1$ and $(\qti_{12}+1)(\qti_{13}+1)(\qti_{23}+1)=0$. 

Moreover, if $q_{11}=-1$, then $q_{22}\qti_{12}=q_{33}\qti_{13}=1$.
\end{lemma}

\begin{lemma}\label{lem:cycles} \cite[Lemma 23]{H-classif}
Let $\bq$ be such that $\varDelta^{\bq}$ is finite. Then the Dynkin diagram of $\bq$ does not contain cycles of length larger than $3$.
\end{lemma}

Finally we summarize the main result of \cite{ASa} on pre-Nichols algebras of Cartan type:

\begin{theorem}\label{thm:AndSanmarco}
Let $\bq$ be of Cartan type with connected Cartan matrix $X$ and parameter $q\in\G_N'$, $N\ge2$. If $(X,N)$ is different from $(A_{\theta},2)$, $(D_{\theta},2)$, $(A_2,3)$, $(G_2,4)$, $(G_2,6)$, then the distinguished pre-Nichols algebra $\toba_{\bq}$ is eminent. 
\qed
\end{theorem}

\section{Defining relations and finite GK dimensional pre-Nichols algebras}\label{sec:def-relations}

Fix $\theta\ge2$ and $\bq=(q_{ij})_{i,j\in\I_\theta}$  with connected Dynkin diagram such that $\GK \toba_\bq <\infty$. By Lemma \ref{lem:1connected}, $q_{ii}\neq 1$ for all $i\in\I_{\theta}$. Let $V=(V,c^\bq)$ be the associated braided vector space with basis $(x_i)_{i\in\I_\theta}$, and $\toba$ a pre-Nichols algebra of $\bq$ such that $\GK \toba<\infty$. 

We determine sufficient conditions under which some defining relations from the presentation of Nichols algebras in \cite[Theorem 3.1]{An-crelle} are annihilated in $\toba$. 

For a given relation $x_u$ the strategy to show that $x_u=0$ in $\toba$ is the following: 
\begin{enumerate}[leftmargin=*,label=\rm{(\alph*)}]
\item We suppose that $x_u\ne 0$, and either we check that $x_u\in\Pc(\toba)$ or we assume this fact. 
\item By Lemma \ref{lem:subspace-primitives}, the braiding matrix $\bq'$ of $V\oplus \Bbbk x_u \subset \Pc(\toba)$ satisfies $\GK \toba_{\bq'}<\infty$.
\item We compute a subdiagram of $\bq'$  and prove that $\GK \toba_{\bq'}=\infty$ using results in \S \ref{subsec:useful-results} (sometimes we invoke Conjecture \ref{conj:AAH} and the classification in \cite{H-classif}), a contradiction.
\end{enumerate}

\medbreak
We start with the relations $x_i^{N_i}$ when $i$ is not of Cartan type.
After that we take care of quantum Serre relations $(\adc x_i)^{1-c_{ij}}x_j$, which are primitive in $T(V)$ by \cite[(4.45)]{H-Lusztig-isos}, see also \cite[Lemma A.1]{AS}; first we give sufficient conditions for them to vanish when $c_{ij}=0$, and later we do the same when $c_{ij}<0$. 
Finally we consider other relations in \cite[Theorem 3.1]{An-crelle} which appear for braidings of super or standard type.

\subsection{Powers of non Cartan vertices}

\begin{lemma}\label{lem:xi-no-Cartan} 
Let $i\in\I_\theta$ be a non-Cartan vertex such that $q_{ii}\in\G_{N}'$. Then $x_i^N=0$. 
\end{lemma}

\pf
As $i$ is not Cartan, there exists $j\neq i$ such that
$\qti_{ij}\notin \{q_{ii}^{-n}: n\in\N_0\}=\G_N$. 

Note that $x_{i}^N$ is a primitive element of $\toba$ since $q_{ii}\in\G_N$. Suppose that $x_i^N\neq 0$. Then $\Pc(\toba)$ contains $\ku x_{i}^N\oplus \ku x_j$, which has Dynkin diagram 
$\xymatrix{ \overset{1}{\circ} \ar@{-}[r]^{\qti_{ij}^N} & \overset{q_{jj}}{\circ}}$. As $\qti_{ij}^N\neq 1$, we get $\GK\toba_{\bq'}=\infty$ from Lemma \ref{lem:subspace-primitives}, a contradiction. Hence $x_i^N=0$.
\epf

\begin{remark}
Let $i,j\in \I_{\theta}$ be such that $q_{ii}\in \G_{N}'$, $\qti_{ij}\notin \G_N$. Then $(\adc x_i)^Nx_j=0$.

This follows by Lemma \ref{lem:xi-no-Cartan} and \cite[(4.45)]{H-Lusztig-isos}.
\end{remark}

\subsection{Quantum Serre relations, $c_{ij}=0$}

\begin{lemma}\label{lem:ij} 
Let $i,j \in \I_{\theta}$ be such that $\qti_{ij}=1$ and either $q_{ii}^2 \neq 1$ or $q_{jj}^2 \neq 1$. Then $x_{ij} = 0$.
\end{lemma}
\pf
Recall that $q_{ii}, q_{jj} \ne 1$. Hence the statement follows by \cite[Proposition 3.2]{ASa}.
\epf

\begin{lemma}\label{lem:ij3} 
Let $i,j \in \I_{\theta}$ be such that $q_{ii}q_{jj} = 1$, $\qti_{ij}=1$ and there exists $\ell \in \I_{\theta}-\{i,j\}$ such that $\qti_{i\ell}\qti_{j\ell} \neq 1$. Then $x_{ij} = 0$.
\end{lemma}
\pf
Suppose that $x_u:=x_{ij} \ne 0$. The subdiagram of $\bq'$ with vertices $u$ and $l$ is $\xymatrix@C=30pt{\overset{q_{\ell \ell}}{\circ}  \ar@{-}[r]^{\qti_{i\ell}\qti_{j\ell}}& \overset{1}{\circ}}$, so $\GK\toba_{\bq'}=\infty$ by Lemma \ref{lem:1connected}, a contradiction.
\epf

\subsection{Quantum Serre relations, $c_{ij}<0$, $q_{ii}^{1-c_{ij}}\neq 1$}
Here the condition $q_{ii}^{1-c_{ij}}\neq 1$ assures that the relation $(\adc x_i)^{1-c_{ij}}x_j=0$ appears in the essentially minimal presentation \cite[Theorem 3.1]{An-crelle}. Following the strategy in \cite[Proposition 4.1]{An-crelle}, we prove:

\begin{lemma}\label{lem:qs}
Let $i,j\in\I_\theta$ be such that $c_{ij}<0$, $q_{ii}^{1-c_{ij}}\neq 1$ and one of the following hold:
\begin{enumerate}[leftmargin=*,label=\rm{(\alph*)}]
\begin{multicols}{2}
\item \label{item:lem-qs-1}$q_{ii}^{2-c_{ij}}\ne 1$, or
\item \label{item:lem-qs-2}$q_{ii}^{c_{ij}(1-c_{ij})}q_{jj}^2 \ne 1$.
\end{multicols}
\end{enumerate}
Then $(\adc x_i)^{1-c_{ij}}x_j=0$ in $\toba$.
\end{lemma}

\pf
We set $q=q_{ii}$, $m=-c_{ij}$. Suppose that $x_u:=(\adc x_i)^{m+1}x_j \neq 0$. By definition of $m$, $q^{m}\qti_{ij}=1$, hence the Dynkin diagram of $\ku x_j \oplus \ku x_i \oplus \ku x_u$ is
\begin{align*}
\xymatrix@C=50pt@R=15pt{&\overset{q^{m+1}q_{jj}}{\circ}&\\
\overset{q_{jj}\,}{\circ} \ar@{-}[rr]^{q^{-m}} \ar@{-}[ur]^{q^{-m(m+1)}q_{jj}^2}&  & \overset{q}{\circ} \ar@{-}[ul]_{q^{m+2}}.}
\end{align*}
Next we use the hypothesis to split the proof in three cases.

\begin{enumerate}[leftmargin=*]
\item $q_{ii}^{m+2}\ne 1$, $q_{ii}^{-m(m+1)}q_{jj}^2 \ne 1$.
By Lemma \ref{lem:rank3}, 
$$ 1=\qti_{ij}\qti_{ju}\qti_{iu}=q^{-m}q^{m+2}q^{-m(m+1)}q_{jj}^2=q^{2-m(m+1)}q_{jj}^2=q^{2-m(m+1)}q_{jj}^2.$$
Moreover, at least one of the vertices has label $-1$. We study each case: notice that $q\ne -1$ since $0<m<\ord(q_{ii})-1$ by hypothesis.
\begin{itemize}[leftmargin=*]
\item $q_{jj}=-1$. By Lemma \ref{lem:rank3}, $m=1$ and $1=(q^{m+1}q_{jj})(q^{-m(m+1)}q_{jj}^2)=(-q^2)(q^{-2})=-1$, a contradiction.

\item $q^{m+1}q_{jj}=-1$. By Lemma \ref{lem:rank3} we also have that $1=qq^{m+2}=q^{m+3}$ and 
$$ 1=(q_{jj})(q^{-m(m+1)}q_{jj}^2)=q^{-m(m+3)+2m}q_{jj}^3=q^{2m}q_{jj}^3 .$$
Hence $-1=(-1)^3=(q^{m+1}q_{jj})^3=(q^{3m+3}q^{-2m})=q^{m+3}=1$, a contradiction.
\end{itemize}

\item $q^{m+2}=1$. By hypothesis $\qti_{ju}=q^{-m(m+1)}q_{jj}^2\neq 1$, and we have the diagram
\begin{align*}
\xymatrix@C=30pt{\overset{q}{\underset{i}{\circ}} \ar@{-}[rr]^{q^2}& &\overset{q_{jj}}{\underset{j}{\circ}}   \ar@{-}[rr]^{q^{-2}q_{jj}^2}&&\overset{q^{-1}q_{jj}}{\underset{u}{\circ}}.
}
\end{align*}
We see that this diagram does not belong to \cite[Table 2]{H-classif} since
\begin{itemize}[leftmargin=*]
    \item $q_{jj} \neq q$ (otherwise we have a vertex with label 1 connected with another vertex);
    \item the extreme vertices have labels $q_{ii}, q_{uu} \neq -1$ whose product is $q_{ii}q_{uu}=q_{jj}$ (the label of the middle vertex) and the product of the two edges is $\qti_{ij}\qti_{ju}=q_{jj}^2$;
    \item if $q_{ii}=\qti_{ij}^{-1}$ and $q_{jj}=-1$, then $q_{ii}^3=1$ and $q_{uu}\qti_{ju}=-1$ ($q_{uu}\neq \qti_{ju}^{-1}$).
\end{itemize}
\item $q^{-m(m+1)}q_{jj}^2=1$; that is, $q_{jj}=\pm q^{\frac{m(m+1)}{2}}$. By hypothesis, $q^{m+2} \neq 1$.
Hence the diagram of $\bq'$ has the form
\begin{align*}
& \xymatrix@C=30pt{\overset{vq^{\frac{(m+1)(m+2)}{2}}}{\underset{u}{\circ}} \ar@{-}[rr]^{q^{m+2}}&&\overset{q}{\underset{i}{\circ}}   \ar@{-}[rr]^{q^{-m}}& &\overset{vq^{\frac{m(m+1)}{2}}}{\underset{j}{\circ}}, } & &v\in \{-1,1\}.
\end{align*}
We check that this diagram does not belong to \cite[Table 2]{H-classif} since
\begin{itemize}[leftmargin=*]
\item the vertex in the middle has label $q_{ii}\neq -1$,
\item the product of the two edges is $\qti_{ij}\qti_{iu}=q_{ii}^2$ (the square of the vertex in the middle),
\item the vertex on the left has label $q_{uu}=\qti_{ij}\qti_{iu}q_{ii}$.
\end{itemize}
\end{enumerate}

\medbreak
All the possibilities lead to a contradiction, so $x_u=0$.
\epf

Next we analyze what happens when $\bq$ does not fulfill the hypothesis of Lemma \ref{lem:qs}.

\begin{lemma}\label{lem:dqs}
Let $i,j\in\I_\theta$ be such that
\begin{itemize}[leftmargin=*]
\begin{multicols}{3}
\item $c_{ij}<0$,
\item $\ord q_{ii}=2-c_{ij}$, 
\item $q_{ii}^{c_{ij}(1-c_{ij})}q_{jj}^2=1$.
\end{multicols}
\end{itemize}
Then $q_{ii}=q_{jj}=\qti_{ij}^{-1}\in\G_3'$ and $c_{ij}=-1$.
\end{lemma}

\pf
By hypothesis, $1=q_{ii}^{c_{ij}(1-c_{ij})}q_{jj}^2=q_{ii}^{-c_{ij}}q_{jj}^2=q_{ii}^{-2}q_{jj}^2$.
Hence $q_{ii}= \pm q_{jj}$. We discard the case $q_{jj}= -q_{ii}$ using Remark \ref{lem:q-q2-(-q)} since the  subdiagram with vertices $i,j$ is
$\xymatrix@C=30pt{\overset{q_{ii}}{\circ} \ar@{-}[r]^{q_{ii}^{c_{ij}}}\ &\overset{-q_{ii}}{\circ}}$. Therefore $q_{ii}=q_{jj}$ and the Dynkin diagram is $\xymatrix@C=30pt{\overset{q_{ii}}{\circ} \ar@{-}[r]^{q_{ii}^{c_{ij}}}\ &\overset{q_{ii}}{\circ}}$. Thus the braiding matrix is of Cartan type with Cartan matrix $\left( \begin{smallmatrix} 2 & c_{ij} \\ c_{ij} & 2 \end{smallmatrix} \right)$. As $\GK \toba_{\bq} <\infty$ we have that $c_{ij}=-1$ and a fortiori $q_{ii}=q_{jj}=\qti_{ij}^{-1}\in\G_3'$.
\epf

Hence the quantum Serre relations hold when $c_{ij}<0$, $q_{ii}^{1-c_{ij}}\ne1$ and the Dynkin diagram spanned by $i,j$ is not of the form $\xymatrix@C=30pt{\overset{\zeta}{\circ} \ar@{-}[r]^{\ztu}\ &\overset{\zeta}{\circ}}$, $\zeta\in\G_3'$. We get now a sufficient condition for the validity of the quantum Serre relations for this exceptional diagram, assuming the existence of another vertex connected either with $i$ or $j$.

\begin{lemma}\label{lem:dqs2}
Let $i,j\in\I_\theta$ be such that $q_{ii}=q_{jj}=\qti_{ij}^{-1}\in\G_3'$, and there exists $k\in\I_\theta-\{i,j\}$ such that $\qti_{ik}^2\qti_{jk}\ne 1$.
Then $(\adc x_i)^2x_j=0$.
\end{lemma}
\pf
Suppose that $x_u:=(\adc x_i)^2x_j \neq 0$. As $\qti_{uk}=\qti_{ik}^2\qti_{jk}\ne 1$ and $q_{uu}=1$, the subdiagram of $\bq'$  with vertices $u$ and $k$ is $\xymatrix@C=30pt{\overset{q_{kk}}{\circ} \ar@{-}[r]^{\qti_{ik}^2\qti_{jk}} &\overset{1}{\circ}}$. So $\GK\toba_{\bq'}=\infty$ by Lemma \ref{lem:1connected}, a contradiction.
\epf

\begin{remark}\label{rem:qs}
The hypothesis on $k\in\I_\theta-\{i,j\}$ of Lemma \ref{lem:dqs2} is fulfilled for example when:
\begin{enumerate}[leftmargin=*,label=\rm{(\alph*)}]
\begin{multicols}{2}
\item either $\qti_{ik}\neq \pm 1$, $\qti_{jk}=1$, or
\item $\qti_{jk}\neq 1$, $\qti_{ik}=1$.
\end{multicols}
\end{enumerate}
\end{remark}

\subsection{Quantum Serre relations, $c_{ij}<0$, $q_{ii}^{1-c_{ij}}=1$}

When $q_{ii}^{1-c_{ij}}=1$ and $i$ is not a Cartan vertex, the relations $(\adc x_i)^{1-c_{ij}}x_j$ hold in the Nichols algebra but are not minimal relations since they follow from $x_i^{1-c_{ij}}=0$. 
Anyway $(\adc x_i)^{1-c_{ij}}x_j$ is primitive in the tensor algebra and belongs to the defining ideal of the distinguished pre-Nichols algebra.

First we work with conditions depending only on $i$, $j$, and later we involve a third vertex.

\begin{lemma}\label{lem:QSR-orden-bajo}
Let $i,j\in\I_\theta$ such that $q_{ii}^{1-c_{ij}}=1$ and $\qti_{ij}=q_{ii}$. If either $q_{jj}\ne-1$ or  $c_{ij} \le - 2$, then $(\adc x_i)^{1-c_{ij}}x_j=0$.
\end{lemma}
\pf
Suppose that $x_u:=(\adc x_i)^{1-c_{ij}}x_j \neq 0$.
As $\qti_{ju}=q_{jj}^2$, $\qti_{ui}=\qti_{ij}$ and $q_{uu}=q_{ii}^{1-c_{ij}}q_{jj}=q_{jj}$, the Dynkin diagram of $W=\ku x_i \oplus \ku x_j \oplus \ku x_u$ is
\begin{align*}
\xymatrix@C40pt@R-25pt{& \overset{q_{jj}}{\underset{u}{\circ}} &\\
\overset{q_{jj}}{\underset{j}{\circ}} \ar@{-}[rr]^{q_{ii}} \ar@{-}[ur]^{q_{jj}^2}&  & \overset{q_{ii}}{\underset{i}{\circ}} \ar@{-}[ul]_{q_{ii}}. }
\end{align*}
Assume first that $q_{jj}\neq -1$, so $\qti_{ij}, \qti_{i u}, \qti_{uj} \neq 1$. 
If $q_{ii}\ne -1$, then the three vertices have labels $\ne-1$, so $\GK \toba (W)=\infty$ by Lemma \ref{lem:rank3}. If $q_{ii}=-1$, then the diagram does not belong to \cite[Table 2]{H-classif}.

It remains to consider the case $q_{jj}=-1$. By hypothesis $\ord q_{ii}\geq 3$ and the Dynkin diagram of $W$ is 
\begin{align*}
\xymatrix@C=40pt{\overset{-1}{\underset{j}{\circ}} \ar@{-}[r]^{q_{ii}} & \overset{q_{ii}}{\underset{i}{\circ}} \ar@{-}[r]^{q_{ii}} & \overset{-1}{\underset{u}{\circ}},}
\end{align*}
which does not belong to \cite[Table 2]{H-classif}. In any case we get a contradiction, so $x_u=0$.
\epf

\begin{lemma}\label{lem:diagonal:qs-m=1-q^m+1=1-extra-vertex}
Let $i,j\in \I_{\theta}$ such that $q_{ii}=\qti_{ij}=-1$, and assume there exists $k\ne i, j$ such that $\qti_{jk}, \qti_{ik}^2\qti_{jk}\neq 1$.
Then $(\adc x_i)^2x_j=0$.
\end{lemma}
\pf
Notice that $c_{ij}=-1$. Suppose that $x_u:=(\adc x_i)^2x_j \neq 0$. By direct computation, $q_{uk}=\qti_{ik}^2\qti_{jk}\neq 1$, $\qti_{iu}=q_{ii}^{5}=-1$, so the Dynkin diagram of $\bq'$ contains a $4$-cycle, a contradiction with Lemma \ref{lem:cycles}. Thus $x_u= 0$.
\epf

\subsection{The relation $[x_{ijk},x_j]_c$}
Next subsections include results about the validity of other relations. First we consider the relation $[x_{ijk},x_j]_c$, which holds in the Nichols algebra when $\qti_{ij}=\qti_{jk}^{-1}$, $\qti_{ik}=1$, $q_{jj}=-1$. 
We begin the subsection stating a formula of the coproduct of $[x_{ijk},x_j]_c$ which becomes important to determine if this element is primitive.

\begin{lemma}\label{lem:-supera-coproduct-[x123,x2]}
If $i,j,k\in\I_\theta$ satisfy $q_{jj}=-1$, $\qti_{ij}=\qti_{jk}^{-1}$, $\qti_{ik}=1$, then the following formula holds in $T(V)$:
\begin{align*}
\Delta& \big([x_{ijk}, x_j]_c \big)
=[x_{ijk}, x_j]_c\ot 1 + 1 \ot [x_{ijk}, x_j]_c + (1-\qti_{ij})\qti_{ij} q_{kj} \ x_i \ot x_{jjk} 
\\
&+ (1-\qti_{jk}) q_{kj} \  [x_{ij}, x_j]_c \ot x_k - (1-\qti_{jk}) q_{ij}^2 q_{kj} \  x_j \ot x_{jik} + 2(1-\qti_{jk}) q_{ij}^2 q_{kj} \ x_j^2 \otimes x_{ik}.
\end{align*}
\end{lemma}

\pf
This is an straightforward computation using the conditions on $\bq$.
\epf

\begin{lemma}\label{lem:diagonal-[x_ijk,x_j]-primitive}
Let $i,j,k\in \I_{\theta}$ be such that $q_{jj}=-1$, $\qti_{ik}=1$ and $\qti_{ij}=\qti_{jk}^{-1}\neq \pm 1 $. Then
\begin{enumerate}[leftmargin=*,label=\rm{(\alph*)}]
\begin{multicols}{2}
\item\label{item:diagonal-x_jik=0} $x_{jik}=0$.
\item\label{item:diagonal-[x_ijk,x_j]-primitive} $[x_{ijk},x_j]_c\in \Pc(\toba)$.
\end{multicols}
\end{enumerate}
\end{lemma}
\pf
\ref{item:diagonal-x_jik=0}
Notice that $x_{jik} \in \Pc(T(V))$ by \cite[Lemma 2.7]{ASa}, because $x_{ik}$ is primitive and $\wbq(\alpha_j, \alpha_i + \alpha_k) = 1$. Also,
\begin{align*}
\bq (\alpha_i+\alpha_k+\alpha_j, \alpha_i+\alpha_k+\alpha_j) &= -q_{ii}q_{kk}, &
\wbq(\alpha_i+\alpha_k+\alpha_j, \alpha_t) &= q_{jt}q_{tj} &&t \in \{i,j,k\}.
\end{align*}
Let $q=\qti_{jk}$. If $x_{jik}\neq 0$ in $\toba$, then $\Pc(\toba)$ contains the $4$-dimensional subspace $\ku x_i \oplus \ku x_j \oplus \ku x_k \oplus \ku x_{jik} $, which has Dynkin diagram
\begin{align*}
\xymatrix@C50pt@R-20pt{ 
& \overset{-q_{ii}q_{kk}}{\underset{jik}{\circ}} \ar  @{-}[dl]_{q^{-1}}^{} \ar  @{-}[dr]^{q}_{} 
\\
\overset{q_{ii}}{\underset{i}{\circ}} \ar  @{-}[r]^{\quad q^{-1}}  & \overset{-1}{\underset{j}{\circ}} \ar  @{-}[r]^{q \quad} &
\overset{q_{kk}}{\underset{k}{\circ}}.}
\end{align*}
Now Lemma \ref{lem:cycles} implies that this Dynkin diagram is not arithmetic. Assuming the validity of Conjecture \ref{conj:AAH}, Lemma \ref{lem:subspace-primitives} says that $\GK \toba = \infty$, a contradiction.

\noindent \ref{item:diagonal-[x_ijk,x_j]-primitive} 
By Lemma \ref{lem:-supera-coproduct-[x123,x2]} and \ref{item:diagonal-x_jik=0}, it is enough to prove that $x_{j}^2=x_{jji}=x_{jjk}=[x_{ij},x_j]_c=0$. For, we observe that $j\in\I_{\theta}$ is not a Cartan vertex since $\qti_{ij}\ne \pm 1$, so $x_j^2=0$ by Lemma \ref{lem:xi-no-Cartan}: this relation implies that $x_{jji}=x_{jjk}=[x_{ij},x_j]_c=0$.
\epf

\begin{lemma}\label{lem:diagonal-[x_ijk,x_j]}
Let $i,j,k\in\I_\theta$ be such that $q_{jj}=-1$, $\qti_{ik}=1$, and $\qti_{ij}=\qti_{jk}^{\,-1} \neq \pm 1$.
If either $q_{ii}=-1$ or $q_{kk}=-1$, then $[x_{ijk},x_j]_c=0$.
\end{lemma}
\pf
By Lemma \ref{lem:diagonal-[x_ijk,x_j]-primitive}, $x_u:=[x_{ijk},x_j]_c\in \Pc(\toba)$. Suppose that $x_u\neq 0$. Set $q=\qti_{jk}=\qti_{ij}^{\,-1}$. It is enough to consider the case $q_{ii}=-1$ (the other is analogous). By direct computation,
\begin{align*}
q_{uu}&=-q_{kk}, & \qti_{ui}&=q^{-2}, & \qti_{uj}&=1, & \qti_{uk}&=q^2q_{kk}^2.
\end{align*}
Hence the Dynkin diagram of the subspace $W=\ku x_i \oplus \ku x_j\oplus \ku x_k \oplus \ku x_u\subset \Pc(\toba)$ is
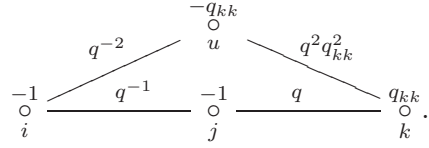
\begin{align*}
\xymatrix@C50pt@R-20pt{&\overset{-q_{kk}}{\underset{u}{\circ}}&\\
\overset{-1}{\underset{i}{\circ}} \ar@{-}[r]^{\quad q^{-1}} \ar@{-}[ur]^{q^{-2}}&  \overset{-1}{\underset{j}{\circ}} \ar@{-}[r]^{q\quad} & \overset{q_{kk}}{\underset{k}{\circ}} \ar@{-}[ul]_{q^2q_{kk}^2}.
}
\end{align*}

If $q^2q_{kk}^2 \neq 1$, then the diagram of $W$ contains a 4-cycle, so from Lemma \ref{lem:cycles} (and assuming Conjecture \ref{conj:AAH}), we get $\GK \toba(W)=\infty$, a contradiction. 

If $q^2q_{kk}^2=1$, that is $q_{kk}=\pm q^{-1}$, then the previous Dynkin diagram becomes
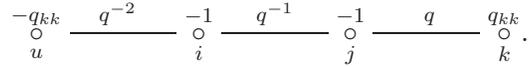
\begin{align*}
\xymatrix@C=40pt{\overset{-q_{kk}}{\underset{u}{\circ}} \ar@{-}[r]^{q^{-2}}&\overset{-1}{\underset{i}{\circ}}  \ar@{-}[r]^{q^{-1}}&\overset{-1}{\underset{j}{\circ}} \ar@{-}[r]^{q} & \overset{q_{kk}}{\underset{k}{\circ}}.
}
\end{align*}
We see that this diagram does not belong to \cite[Table 3]{H-classif} since
\begin{itemize}[leftmargin=*]
\item the vertices in the middle are $-1$,
\item the labels of the extreme vertices are opposite and both different of $-1$.
\end{itemize}
In any case we get a contradiction, so $x_u=[x_{ijk},x_j]_c=0$.
\epf

\begin{lemma}\label{lem:cor2}
Let $i,j,k\in\I_\theta$ be such that $q_{jj}=-1$, $\qti_{ik}=1$ and $\qti_{ij}=\qti_{jk}^{\,-1} \neq \pm 1$.
If $q_{ii}q_{kk}=1$ and there exists $\ell\in\I_\theta$ such that one of the following conditions holds
\begin{enumerate}[leftmargin=*,label=\rm{(\alph*)}]
\begin{multicols}{3}
\item\label{item:cor2-a} $\qti_{i\ell}\neq 1=\qti_{j\ell}=\qti_{k\ell}$,
\item\label{item:cor2-b} $\qti_{j\ell}^2\neq 1=\qti_{i\ell}=\qti_{k\ell}$,
\item\label{item:cor2-c} $\qti_{k\ell}\neq 1=\qti_{j\ell}=\qti_{i\ell}$,
\end{multicols}
\end{enumerate}
then $[x_{ijk},x_j]_c=0$.
\end{lemma}

\pf
By Lemma \ref{lem:diagonal-[x_ijk,x_j]-primitive}, $x_u:=[x_{ijk},x_j]_c\in \Pc(\toba)$.
Suppose that $x_u \neq 0$. Set $q=\qti_{jk}=\qti_{ij}^{-1}$. As $q_{uu}=q_{ii}q_{kk}=1$ and $\qti_{u\ell}=\qti_{i\ell}\qti_{j\ell}^2\qti_{k\ell}\ne 1$, the subdiagram of $\bq'$ with vertices $u$ and $\ell$ is $\xymatrix@C=30pt{\overset{1}{\circ} \ar@{-}[r]^{\qti_{u\ell}}& \overset{q_{\ell\ell}}{\circ}}$. Thus $\GK\toba_{\bq'}=\infty$ by Lemma \ref{lem:1connected} and we get a contradiction.
\epf

\subsection{Other relations} Next we consider other relations listed in \cite[Theorem 3.1]{An-crelle} that appear in the defining ideal of Nichols algebras of super type.

\begin{lemma}\label{lem:diagonal-[x_iij,x_ij]}
Let $i,j\in \I_{\theta}$ be such that $q_{ii}\qti_{ij}\in\G_3'\cup\G_6'$, $q_{jj}=-1,$ and either $q_{ii}\in\G_3'$ or $c_{ij} \le - 3$. If $[x_{iij},x_{ij}]_c\in\Pc(\toba)$, then $[x_{iij},x_{ij}]_c=0$.
\end{lemma}
\pf
We follow the proof of \cite[Lemma 4.3 (i)]{An-crelle}.
Suppose that $x_u:=[x_{iij},x_{ij}]_c\neq 0$. Since $q_{uu}=q_{ii}^9\qti_{ij}^{\,6}=q_{ii}^3$, $\qti_{ui}=\qti_{ij}^{\,3}q_{jj}^4=\qti_{ij}^{\,3}$, $\qti_{uj}=q_{ii}^{\,6}\qti_{ij}^{\,2}=\qti_{ij}^{\,-4}$, the Dynkin diagram of  $\ku x_j \oplus \ku x_i \oplus \ku x_u$ is
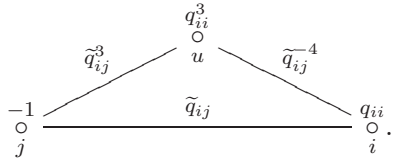
\begin{align*}
\xymatrix@C50pt@R-20pt{ & {\overset{q_{ii}^3}{\underset{u}{\circ}}}&\\
\overset{-1}{\underset{j}{\circ}}\ar@{-}[rr]^{\qti_{ij}} \ar@{-}[ur]^{\qti_{ij}^3} & & \overset{q_{ii}}{\underset{i}{\circ}}\ar@{-}[ul]_{\qti_{ij}^{-4}}. }
\end{align*}

If $q_{ii}\in\G_3'$, then $\GK \toba_{\bq'}=\infty$ by Lemma \ref{lem:1connected}, since either $\qti_{ij}^{\,3}\ne 1$ or  $\qti_{ij}^{\,4}\ne 1$.

If $c_{ij} \le - 3$ then this diagram is connected, and as we are assuming Conjecture \ref{conj:AAH}, it belongs to \cite[Table 2]{H-classif}. We get a contradiction since the unique diagram in \emph{loc. cit.} with some  $c_{rs}\le -3$ is the first one in row 7, but this forces $\qti_{ij}^{3}=1$,  $\qti_{ij}=q_{ii}^{-3}=-1$. In any case we get a contradiction so $x_u=0$.
\epf

\begin{lemma}\label{lem:diagonal-[x_iijk,x_ij]}
Let $i,j,k\in \I_{\theta}$ be such that $q_{ii}=\pm\qti_{ij}\in\G_3'$, $\qti_{ik}=1$ and one of the following conditions on $q_{jj}$, $\qti_{ij}$, $\qti_{jk}$ hold:
\begin{enumerate}[leftmargin=*,label=\rm{(\alph*)}]
\begin{multicols}{2}
\item\label{item:diagonal-[x_iijk,x_ij]-a} $q_{jj}=-1$, $\qti_{ij}=\qti_{jk}^{\, -1}$, or
\item\label{item:diagonal-[x_iijk,x_ij]-b} $q_{jj}^{-1}=\qti_{ij}=\qti_{jk}$.
\end{multicols}
\end{enumerate}
If $[x_{iijk},x_{ij}]_c\in\Pc(\toba)$, then $[x_{iijk},x_{ij}]_c=0$.
\end{lemma}
\pf
Suppose that $x_u:=[x_{iijk},x_{ij}]_c\neq 0$. Set $\xi=q_{ii}$. By direct computation,
\begin{align*}
\qti_{uk}&=\qti_{jk}^{\,2}q_{kk}^2,& \qti_{uj}&=q_{jj}^4\qti_{ij}^{\,3}\qti_{jk}, &
\qti_{ui}&=\xi^2\neq 1, & q_{uu}&=q_{jj}^4q_{kk}\qti_{jk}^{\,2}.
\end{align*}

Assume first that \ref{item:diagonal-[x_iijk,x_ij]-a} holds. The Dynkin diagram of $\ku x_i \oplus \ku x_j\oplus \ku x_k\oplus \ku x_u$ is
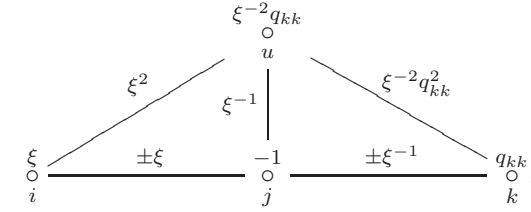
\begin{align*}
\xymatrix@C=30pt{
&&\overset{\xi^{-2}q_{kk}}{\underset{u}{\circ}}\ar@{-}[dll]_{\xi^2}\ar@{-}[drr]^{\xi^{-2}q_{kk}^2}&&\\
\overset{\xi}{\underset{i}{\circ}}\ar@{-}[rr]^{\pm\xi}&&\overset{-1}{\underset{j}{\circ}}\ar@{-}[rr]^{\pm\xi^{-1}}\ar@{-}[u]^{\xi^{-1}}&&\overset{q_{kk}}{\underset{k}{\circ}},}
\end{align*}
so $x_i,x_j,x_u$ determine a triangle: by Lemma \ref{lem:rank3}, $1=q_{ii}\qti_{ij}=\pm\xi^2$, a contradiction.

Assume now that \ref{item:diagonal-[x_iijk,x_ij]-b} holds. The subdiagram of $\bq'$ spanned by the vertices $i,j,k,u$ is
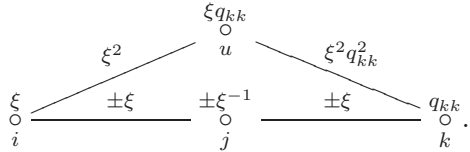
\begin{align*}
\xymatrix@C60pt@R-20pt{
&\overset{\xi q_{kk}}{\underset{u}{\circ}}\ar@{-}[dl]_{\xi^2}\ar@{-}[dr]^{\xi^{2}q_{kk}^2}&\\
\overset{\xi}{\underset{i}{\circ}}\ar@{-}[r]^{\pm\xi}&\overset{\pm \xi^{-1} }{\underset{j}{\circ}}\ar@{-}[r]^{\pm\xi}&\overset{q_{kk}}{\underset{k}{\circ}}.}
\end{align*}
If $\qti_{uk}=\xi^2q_{kk}^2\neq 1$, then the diagram above is a $4$-cycle, a contradiction to Lemma \ref{lem:cycles}.
If $q_{kk}=\xi^{-1}$, then $q_{uu}=1$ and $\qti_{iu}\neq 1$, a contradiction with Lemma \ref{lem:1connected}. Finally we assume $q_{kk}=-\xi^{-1}$. The diagram above becomes
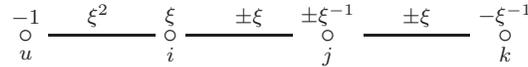
\begin{align*}
\xymatrix@C=40pt{
\overset{-1}{\underset{u}{\circ}}\ar@{-}[r]^{\xi^2}&\overset{\xi}{\underset{i}{\circ}}\ar@{-}[r]^{\pm\xi}&\overset{\pm \xi^{-1} }{\underset{j}{\circ}}\ar@{-}[r]^{\pm\xi}&\overset{-\xi^{-1}}{\underset{k}{\circ}}}
\end{align*}
which does not belong to \cite[Table 3]{H-classif} since
\begin{itemize}[leftmargin=*]
\item there is a unique vertex with label $=-1$, and this vertex is an extreme;
\item $q_{uu}=-1$, $q_{ii}=\qti_{ui}^{\,-1}\in \G_3'$ and $q_{kk}\in\G_6'$.
\end{itemize}
In any case we get a contradiction, so $x_u=0$.
\epf

\begin{lemma}\label{lem:diagonal-xij^2}
Let $i,j\in \I_{\theta}$ be such that $q_{ii}=\qti_{ij}=q_{jj}=-1$.
\begin{enumerate}[leftmargin=*,label=\rm{(\alph*)}]
\item\label{item:diagonal-xij^2-a} The following formula holds in $T(V)$:
\begin{align*}
\Delta (x_{ij}^2)&=x_{ij}^2\ot 1 + 1\ot x_{ij}^2 +2q_{ji} x_i^2 \ot x_j^2
-2q_{ji} [x_i, x_{ij}]_c \ot x_j - 2q_{ji} x_i\ot [x_{ij},x_j]_c.
\end{align*}

\item\label{item:diagonal-xij^2-b} 
If there exists $k\in\I_{\theta}-\{i,j\}$ such that $\qti_{ik}^2\qti_{jk}^2\neq 1$, then $x_{ij}^2=0$.
\end{enumerate}

\end{lemma}
\pf
For \ref{item:diagonal-xij^2-a} we note that $\Delta (x_{ij})=x_{ij}\ot 1+ 1\ot x_{ij}+2x_i\ot x_j$. Then
\begin{align*}
\Delta (x_{ij}^2)
&= x_{ij}^2\ot 1 + 1\ot x_{ij}^2 +2q_{ji} x_i^2 \ot x_j^2 + (1+q_{ii}\qti_{ij}q_{jj}) x_{ij}\ot x_{ij}
\\ & \quad -2q_{ji}x_ix_{ij}\ot x_j+ 2 x_{ij}x_i\ot x_j
+ 2x_i\ot \ot x_jx_{ij}- 2q_{ji} x_i\ot x_{ij}x_j
\\ &= x_{ij}^2\ot 1 + 1\ot x_{ij}^2 +2q_{ji} x_i^2 \ot x_j^2
-2q_{ji} [x_i, x_{ij}]_c \ot x_j - 2q_{ji} x_i\ot [x_{ij},x_j]_c.
\end{align*}

For \ref{item:diagonal-xij^2-b}, we suppose that $x_u:=x_{ij}^2\neq 0$. We check first that $x_{ij}^2\in\Pc(\toba)$.
By hypothesis, either $\qti_{ik}^2\neq 1$ or $\qti_{jk}^2\neq 1$: we may assume that $\qti_{jk}^2\neq 1$ so $j$ is a not a Cartan vertex. By Lemma \ref{lem:xi-no-Cartan}, $x_j^2=0$, and this relation implies that $[x_{ij},x_j]_c=0$.
If $\qti_{ik}^2\neq 1$, then the same result implies that $x_i^2=0$, thus $[x_i, x_{ij}]_c=0$. Otherwise, $\qti_{ik}^2\qti_{jk}=\qti_{jk}\neq 1$ and then Lemma \ref{lem:diagonal:qs-m=1-q^m+1=1-extra-vertex} says that
$[x_i, x_{ij}]_c=0$. Hence $x_{ij}^2\in\Pc(\toba)$.

By direct computation, $q_{uu}=1$, $\qti_{uk}=\qti_{ik}^2\qti_{jk}^2\neq 1$. Then the Dynkin diagram of $\Bbbk x_{k}+\Bbbk x_{u}$ is $\xymatrix@C40pt{\overset{q_{kk}}{\circ}\ar@{-}[r]^{\qti_{ik}^2\qti_{jk}^2}&\overset{1}{\circ}} $, and by Lemma \ref{lem:1connected} we get a contradiction. Thus $x_{ij}^2= 0$.
\epf

\begin{lemma}\label{lem:diagonal:[[xijk,xj]xj]}
Let $i,j,k\in\I_{\theta}$ be such that $q_{jj}=\qti_{ij}^{-1}=\qti_{jk}\in \G_{3}'$, $\qti_{ik}=1$, and either $q_{ii}\ne -1$ or else $q_{kk}\ne -1$.
Then $[[x_{ijk},x_{j}]_c,x_j]_c=0$.
\end{lemma}
\pf
First we prove that $x_u=[[x_{ijk},x_{j}]_c,x_j]_c\in\Pc(\toba)$. We have that:
\begin{itemize}[leftmargin=*]
\item $x_{ik}=0$ in $\toba$ applying Lemma \ref{lem:ij};
\item $x_{jji}=0$ in $\toba$ applying Lemma \ref{lem:qs} if $q_{ii}\ne q_{jj}$, and Lemma \ref{lem:dqs2} otherwise; and
\item $x_{jjjk}=0$ in $\toba$ by Lemma \ref{lem:QSR-orden-bajo} (here $-c_{jk}=2$).
\end{itemize}
By direct computation,
\begin{align*}
x_u&=[[x_{ijk},x_j]_c,x_j]_c= x_ix_jx_kx_j^2 -q_{jk}x_ix_kx_j^3 -q_{ij}q_{ik}x_jx_kx_ix_j^2 +q_{jk}q_{ij}q_{ik}x_kx_jx_ix_j^2 
\\
&-q_{ij}q_{jj}q_{kj}(1+q_{jj})x_jx_ix_jx_kx_j+q_{ij}q_{jj}\qti_{jk}(1+q_{jj})x_jx_ix_kx_j^2
\\
&+q_{ij}^2q_{jj}q_{kj}q_{ik}(1+q_{jj})x_j^2x_kx_ix_j-q_{ij}^2q_{ik}q_{jj}\qti_{jk}(1+q_{jj})x_jx_kx_jx_ix_j
\\
&+q_{ij}^2q_{kj}^2x_j^2x_ix_jx_k-q_{jk}q_{ij}^2q_{kj}^2x_j^2x_ix_kx_j
-q_{ik}q_{ij}^3q_{kj}^2x_j^3x_kx_i+q_{jk}q_{ik}q_{ij}^3q_{kj}^2x_j^2x_kx_jx_i.
\end{align*}
Using explicit formulas for the comultiplication of each summand we get
\begin{align*}
\Delta & (x_u)=x_u \ot 1+1\ot x_u +3q_{ij}^3q_{kj}^2(1-\qti_{jk})x_j^3\ot x_{ik}
\\
& +q_{kj}^2(1-\qti_{jk}) (x_ix_j^3-q_{ij}^3 x_j^3x_i)\ot x_k+(\qti_{ij}^2-1)q_{kj}^2x_i\ot (x_j^3x_k-q_{jk}^3 x_kx_j^3).
\end{align*}
Notice that $x_{jjji}=(\adc x_j)x_{jji}=0$, so $x_ix_j^3=q_{ij}^3 x_j^3x_i$. Similarly, $x_j^3x_k=q_{jk}^3 x_kx_j^3$.
Using these relations and $x_{ik}=0$, we get that $x_u\in\Pc(\toba)$.
\smallbreak

Suppose that $x_u\neq 0$. We have that
\begin{align*}
\qti_{uk}&=q_{kk}^2, & \qti_{uj}&=1, & \qti_{ui}&=q_{ii}^2, & q_{uu}&=q_{ii}q_{kk}.
\end{align*}
Hence the subdiagram of $\bq'$ spanned by the vertices $i,j,k,u$ is
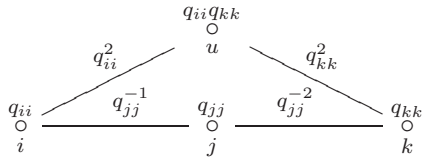
\begin{align*}
\xymatrix@C50pt@R-15pt{&\overset{q_{ii}q_{kk}}{\underset{u}{\circ}}\ar@{-}[dl]_{q_{ii}^2}\ar@{-}[dr]^{q_{kk}^2}&\\
\overset{q_{ii}}{\underset{i}{\circ}} \ar@{-}[r]^{\quad q_{jj}^{-1}} & \overset{q_{jj}}{\underset{j}{\circ}}\ar@{-}[r]^{q_{jj}^{-2}\quad} &\overset{q_{kk}}{\underset{k}{\circ}}}
\end{align*}
\begin{itemize}[leftmargin=*]
\item If $q_{kk}=-1$, then $q_{ii} \ne -1$ by hypothesis and the subdiagram with vertices $i$, $u$ is as in Remark \ref{lem:q-q2-(-q)}, then $\GK \toba_{\bq'}=\infty$. 

\item If $q_{ii}=-1$, then $q_{kk} \ne -1$ and $\GK \toba_{\bq'}=\infty$ by Remark \ref{lem:q-q2-(-q)}. 

\item If $q_{ii},q_{kk} \ne -1$, the subdiagram with vertices $i,j,k,u$ is a $4$-cycle, so $\GK \toba_{\bq'}=\infty$ by Lemma \ref{lem:cycles}.
\end{itemize}
In any case we get a contradiction, thus $x_u=0$.
\epf

\begin{lemma}\label{lem:diagonal-[[x_ij,x_ijk],x_j]}
Let $i,j,k\in \I_{\theta}$ be such that $\qti_{ik}=1$, $q_{jj}=q_{ii}=-1$, $\qti_{ij}^{\,2}=\qti_{jk}^{\,-1}\ne 1$, and either
$q_{kk}\ne-1$ or $\qti_{ij}^{\, 3}\ne 1$. If $[[x_{ij},x_{ijk}]_c,x_j]_c\in \Pc(\toba)$, then $[[x_{ij},x_{ijk}]_c,x_j]_c=0$.
\end{lemma}
\pf
Suppose that $x_u:=[[x_{ij},x_{ijk}]_c,x_j]_c\neq 0$. We have that 
\begin{align*}
    \qti_{uk}&=\qti_{jk}^{\,3}q_{kk}^2,& \qti_{uj}&=1, & \qti_{ui}&=\qti_{ij}^{\, 3}, & q_{uu}&=-q_{kk}.
\end{align*}
Then the Dynkin diagram of $\bq'$ contains that of $\ku x_i \oplus \ku x_j \oplus\ku x_k \oplus\ku x_u$, which is 
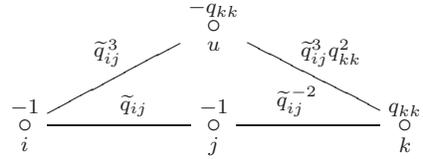
\begin{align*}
\xymatrix@C50pt@R-15pt{&\overset{-q_{kk}}{\underset{u}{\circ}}\ar@{-}[dl]_{\qti_{ij}^{\,3}}\ar@{-}[dr]^{\qti_{ij}^3q_{kk}^2}&\\
\overset{-1}{\underset{i}{\circ}} \ar@{-}[r]^{\quad\qti_{ij}} & \overset{-1}{\underset{j}{\circ}}\ar@{-}[r]^{\qti_{ij}^{\, -2}\quad} &\overset{q_{kk}}{\underset{k}{\circ}}.}
\end{align*}
\begin{itemize}[leftmargin=*]
\item If $q_{kk}=-1$, then $\qti_{ij}^3 \ne 1$ by hypothesis and the subdiagram with vertices $i,u$ is as in Lemma \ref{lem:1connected}, thus $\GK \toba_{\bq'}=\infty$. 

\item If $\qti_{ij}^3=1$, then $q_{kk}^2 \ne 1$ and the subdiagram with vertices $i,u$ is as in Remark \ref{lem:q-q2-(-q)}, hence $\GK \toba_{\bq'}=\infty$. 

\item If $q_{kk}^2, \qti_{ij}^3 \ne 1$, then the Dynkin diagram above is a $4$-cycle, so $\GK \toba_{\bq'}=\infty$ by Lemma \ref{lem:cycles}.
\end{itemize}
In any case we get a contradiction, thus $x_u=0$.
\epf

\begin{lemma}\label{lem:diagonal-[x_ij,x_ijk]}
Let $i,j,k\in \I_{\theta}$ be such that $q_{ii}=\qti_{ij}=-1$, $q_{jj}=\qti_{jk}^{\,-1}\neq -1$, $\qti_{ik}=1$ and either $q_{kk}\ne q_{jj}$ or else $q_{jj}\notin\G_{3}'$. Then $[x_{ij},x_{ijk}]_c=0$.
\end{lemma}
\pf
Suppose that $x_u:=[x_{ij},x_{ijk}]_c\neq 0$. First we prove that $x_u\in \Pc(\toba)$. It is enough to prove that all defining relations of $\toba_{\bq}$ of degree lower that $x_u$ hold in $\toba$:
\begin{itemize}[leftmargin=*]
\item If $i$ is not a Cartan vertex, then $x_i^2=0$ by Lemma \ref{lem:xi-no-Cartan}.

\item If $i$ is a Cartan vertex, then $x_{iij} =0$ in $\toba$ by Lemma \ref{lem:diagonal:qs-m=1-q^m+1=1-extra-vertex}: indeed, $\qti_{jk}=\qti_{ik}^2\qti_{jk}\ne 1$.

\item By hypothesis, $q_{kk}\ne q_{jj}$ or $q_{jj}\notin\G_{3}'$, therefore $x_{jjk}=0$ in $\toba$ by Lemma \ref{lem:qs}.

\item If $q_{kk}\ne -1$, then $x_{ik}=0$ by Lemma \ref{lem:ij}. If $q_{kk}=-1$, then $x_{ik}=0$ by Lemma \ref{lem:ij3} since $\qti_{ij}\qti_{jk}=-\qti_{jk}\ne 1$.
\end{itemize}

By direct computation,
\begin{align*}
\qti_{uk}&=q_{jj}^{-2}q_{kk}^2,& \qti_{uj}&=q_{jj}^3\neq 1, & \qti_{ui}&=1, &\qquad q_{uu}&=q_{jj}^2q_{kk}.
\end{align*}
Hence the subdiagram of $\bq'$ spanned by the vertices $i,j,k,u$ is
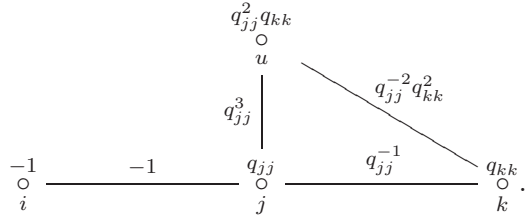
\begin{align*}
\xymatrix@C=30pt{&&{\overset{q_{jj}^2q_{kk}}{\underset{u}{\circ}}}\ar@{-}[d]_{q_{jj}^3}\ar@{-}[drr]^{q_{jj}^{-2}q_{kk}^2}&&\\
\overset{-1}{\underset{i}{\circ}} \ar@{-}[rr]^{-1} && \overset{q_{jj}}{\underset{j}{\circ}}\ar@{-}[rr]^{q_{jj}^{-1}} &&\overset{q_{kk}}{\underset{k}{\circ}}.
}
\end{align*}
If $q_{kk}= -q_{jj}$, then the subdiagram with vertices $j,k$ does not belong in \cite[Table 1]{H-classif}, so $q_{kk}\ne -q_{jj}$. We split the study in the following three cases:

\begin{itemize}[leftmargin=*]
\item $q_{jj}^3 \neq 1$, $q_{kk}^2 \ne q_{jj}^2$: the subdiagram with vertices $j,k,u$ is a 3-cycle. By Lemma \ref{lem:rank3}, $1=\qti_{jk}\qti_{uk}\qti_{uj}=q_{kk}^2$, so $q_{kk}=-1$, and $1=\qti_{uk}q_{uu}=-1$, a contradiction.

\item $q_{jj}^3=1$, $q_{kk}^2 \ne q_{jj}^2$. The subdiagram spanned by the vertices $i,j,k$ is
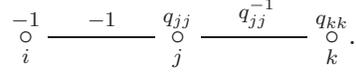
\begin{align*}
\xymatrix@C=40pt{\overset{-1}{\underset{i}{\circ}} \ar@{-}[r]^{-1} & \overset{q_{jj}}{\underset{j}{\circ}}\ar@{-}[r]^{q_{jj}^{-1}} &\overset{q_{kk}}{\underset{k}{\circ}}.}
\end{align*}
Looking at \cite[Table 2]{H-classif}, the unique possibility is $q_{kk}=-1$. This says that the Dynkin diagram of $\Bbbk x_k+\Bbbk x_u$ is $\xymatrix@C=30pt{\overset{-1}{{\circ}}\ar@{-}[r]^{q_{jj}^{-2}} & {\overset{-q_{jj}^2}{{\circ}}}}$, which does not belong at \cite[Table 1]{H-classif}.

\item $q_{jj}^3 \neq 1$, $q_{kk}^2=q_{jj}^2$: that is, $q_{kk}=q_{jj}$. The subdiagram spanned by $i,j,k,u$ is:
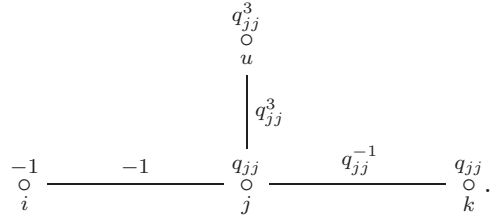
\begin{align*}
\xymatrix@C=30pt{&&{\overset{q_{jj}^3}{\underset{u}{\circ}}}\ar@{-}[d]^{q_{jj}^3}&&\\
\overset{-1}{\underset{i}{\circ}} \ar@{-}[rr]^{-1}&&\overset{q_{jj}}{\underset{j}{\circ}}\ar@{-}[rr]^{q_{jj}^{-1}} &&\overset{q_{jj}}{\underset{k}{\circ}}.}
\end{align*}
As $q_{jj}\notin\G_2\cup\G_3$, this diagram does not belong to \cite[Table 3]{H-classif}.
\end{itemize}
In any case we get a contradiction, thus $x_u=0$.
\epf

\begin{lemma}\label{lem:diagonal-[[[x_ijkl,x_k],x_i],x_k]}
Let $i,j,k,\ell\in \I_{\theta}$ be such that $q_{kk}= -1$, $q_{jj}\qti_{ij}=q_{jj}\qti_{jk}=1$, $\qti_{ik}=\qti_{il}=\qti_{j\ell}=1$ and $\qti_{jk}^2=\qti_{k\ell}^{-1}=q_{\ell\ell}$. If $[[[x_{ijk\ell},x_{k}]_c,x_{j}]_c,x_{k}]_c \in \Pc(\toba)$, then $[[[x_{ijk\ell},x_{k}]_c,x_{j}]_c,x_{k}]_c=0$.
\end{lemma}
\pf
Suppose that $x_u:=[[[x_{ijk\ell},x_{k}]_c,x_{j}]_c,x_{k}]_c\neq 0$. We set $q=\qti_{ij}$ and $p=q_{ii}$. We have
\begin{align*}
q_{uu}&=-p, & \qti_{u\ell}&=q^{-2}\neq 1, & \qti_{uk}&=1, & \qti_{uj}&=q^4, &\qti_{ui}&=p^{2}q^2,
\end{align*}
Hence the Dynkin diagram of $\bq'$ contains the following one:
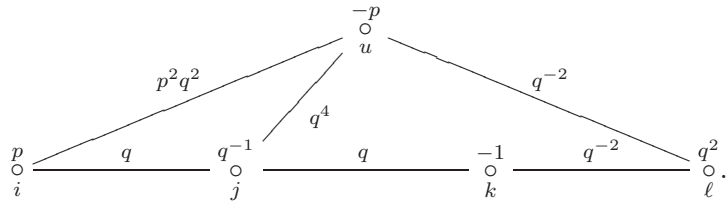
\begin{align*}
\xymatrix@C=30pt{&&& \overset{-p}{\underset{u}{\circ}}\ar@{-}[dl]^{q^4}\ar@{-}[dlll]_{{p}^{2}q^2}\ar@{-}[drrr]^{q^{-2}}&&&\\
\overset{p}{\underset{i}{\circ}} \ar@{-}[rr]^{q} &&  \overset{q^{-1}}{\underset{j}{\circ}} \ar@{-}[rr]^{q} && \overset{-1}{\underset{k}{\circ}}\ar@{-}[rr]^{q^{-2}} &&\overset{q^2}{\underset{\ell}{\circ}}.}
\end{align*}
If either $q\notin\G_4'$ or $p \ne\pm q^{-1}$, then $\GK\toba=\infty$ by Lemma \ref{lem:cycles}. Next we assume that $q\in\G_4'$, $p=\pm q^{-1}$. Then the diagram with vertices $j, k, \ell, u$ is
\begin{align*}
\xymatrix@C=40pt{\overset{q^{-1}}{\underset{j}{\circ}} \ar@{-}[r]^{q} & \overset{-1}{\underset{k}{\circ}}\ar@{-}[r]^{-1} &\overset{-1}{\underset{l}{\circ}}\ar@{-}[r]^{-1}
& \overset{\pm q}{\underset{u}{\circ}}. }
\end{align*}
This diagram does not belong to \cite[Table 3]{H-classif} since:
\begin{itemize}[leftmargin=*]
\item the labels of the vertices of degree one are roots of order $4$;
\item each vertex of degree two and the edge between them are labelled with $-1$;
\item the label of one of extreme vertices is not the inverse of the adjacent edge.
\end{itemize}
As we are assuming Conjecture \ref{conj:AAH} we get $x_u=0$.
\epf

\begin{lemma}\label{lem:diagonal-x_ijk}
Let $i,j,k\in \I_{\theta}$ be such that $\qti_{ij},\qti_{ik},\qti_{jk}\neq 1$, $\qti_{ij}\qti_{ik}\qti_{jk}=1$. Then
$$ x_{ijk} = q_{ij}(1-\qti_{jk})x_{j}x_{ik}-\frac{1-\qti_{jk}}{q_{kj}(1-\qti_{ik})}[x_{ik},x_{j}]_c.$$
\end{lemma}
\pf
Suppose that $x_u:=x_{ijk}- q_{ij}(1-q^2)x_{j}x_{ik} +q_{jk}(1+q^{-1})[x_{ik},x_{j}]_c\neq 0$. By direct computation, $x_u\in \Pc(\toba)$. The subdiagram of $\bq'$ with vertices $i,j,k,u$ is
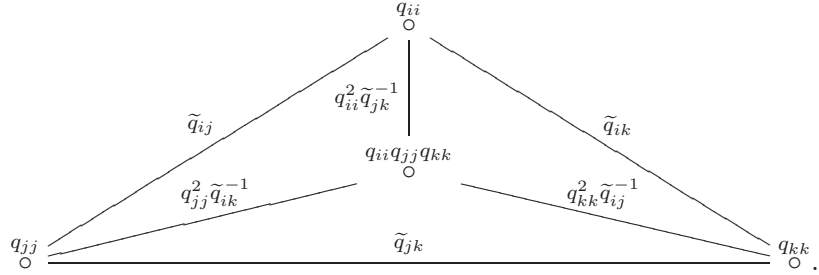
\begin{align*}
\xymatrix@C=55pt@R=15pt{&&\overset{q_{ii}}{\circ}\ar@{-}[dddll]_{\qti_{ij}} \ar@{-}[dd]_{q_{ii}^2\qti_{jk}^{\, -1}}\ar@{-}[dddrr]^{\qti_{ik}}&&\\
&&&&\\
&&\overset{q_{ii}q_{jj}q_{kk}}{\circ}\ar@{-}[dll]_{\quad q_{jj}^2\qti_{ik}^{\, -1}}\ar@{-}[drr]^{q_{kk}^2\qti_{ij}^{\, -1} \quad }&&\\\
\overset{q_{jj}}{\circ}\ar@{-}[rrrr]^{\qti_{jk}}&&&&\overset{q_{kk}}{\circ}.}
\end{align*}
By Lemma \ref{lem:rank3} there exists $\ell\in \{i,j,k\}$ such that $q_{\ell\ell}=-1$. If there exist $\ell_1 \ne \ell_2\in \{i,j,k\}$ such that $q_{\ell_1\ell_1}=q_{\ell_2\ell_2}=-1$, then $\qti_{u\ell_1}, \qti_{u\ell_2}\neq 1$. Thus the previous diagram contains a $4$-cycle, so $\GK \toba=\infty$ by Lemma \ref{lem:cycles}. 

Next we assume that there exists a unique $\ell\in \{i,j,k\}$ such that $q_{\ell\ell}=-1$; relabeling the vertices, we can assume that $q_{ii}=-1\neq q_{jj},q_{kk}$, and also $q_{jj}^2\qti_{ik}^{-1}=1=q_{kk}^2\qti_{ij}^{\,-1}$ (otherwise the diagram still contains a 4-cycle since $q_{ui}=\qti_{jk}^{\, -1}\neq 1$).
By Lemma \ref{lem:rank3}, $q_{jj}\qti_{ij}=1=q_{kk}\qti_{ik}$. 
Set $q=\qti_{jk}$. Notice that
\begin{align*}
&1=q_{jj}^2\qti_{ik}^{\,-1}=\qti_{ij}^{\,-2}\qti_{ik}^{\,-1}=\qti_{ij}^{\,-1}q, &
&1=q_{kk}^2\qti_{ij}^{\,-1}=\qti_{ik}^{\,-2}\qti_{ij}^{\,-1}=\qti_{ik}^{\,-1}q,
\end{align*}
so $\qti_{ij}=\qti_{ik}=q$, and therefore $q\in\G_3'$. 
Now the previous diagram becomes
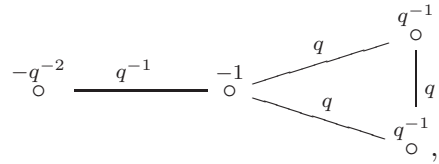
\begin{align*}
\xymatrix@C50pt@R-25pt{ & &\overset{q^{-1}}{\circ} \ar@{-}[dd]^{q}\\
\overset{-q^{-2}}{\circ} \ar@{-}[r]^{q^{-1}}& \overset{-1}{\circ} \ar@{-}[ur]^{q} \ar@{-}[dr]^{q}&\\
 &&\overset{q^{-1}}{\circ},}
\end{align*}
which does not belong to \cite[Table 3]{H-classif}. This contradiction shows that $x_u=0$.
\epf

\begin{lemma}\label{lem:diagonal-[[[xijk,xj],[xijkl,xj]],xjk]} 
Let $i,j,k,\ell \in \I_{\theta}$ and $q\in\Bbbk-(\G_2\cup\G_3)$ be such that 
\begin{align*}
q_{\ell \ell}&=\qti_{k\ell}^{\,-1}=q_{kk}=\qti_{jk}^{\,-1}=q^2, &
q_{jj}&=-1, & q_{ii}&=\qti_{ij}^{\,-1}=q^{-3}, & \qti_{ik}&=\qti_{i\ell}=\qti_{j\ell}=1.
\end{align*}
If $[[[x_{ijk},x_j]_c,[x_{ijk\ell},x_j]_c]_c,x_{jk}]_c\in\Pc(\toba)$, then 
$[[[x_{ijk},x_j]_c,[x_{ijk\ell},x_j]_c]_c,x_{jk}]_c=0$.
\end{lemma}
\pf
Suppose that $x_u:=[[[x_{ijk},x_j]_c,[x_{ijk\ell},x_j]_c]_c,x_{jk}]_c\neq 0$. By direct computation,
\begin{align*}
\qti_{ui}&=q_{ii}^4\qti_{ij}^{\, 5}=q^{-3}\neq 1, & \qti_{u\ell}&=q_{\ell\ell}^2\qti_{k\ell}^{\, 3}=q^{-2}\neq 1.
\end{align*}
Then the subdiagram with vertices $i,j,k,\ell,u$ contains a 5-cycle. Therefore $\GK \toba =\infty$ by Lemma \ref{lem:cycles}, a contradiction.
\epf

\begin{lemma}\label{lem:diagonal-[[[xijkl,xj]xk]} 
Let $i,j,k,\ell \in \I_{\theta}$ and $q\in\Bbbk-(\G_2\cup\G_3)$ be such that 
\begin{align*}
q_{ii}&=\qti_{ij}^{\,-1}=q^2, &
q_{kk}&=-1, & q_{\ell\ell}^{-1}&=\qti_{k\ell}=q^{3}, & \qti_{jk}^{-1}=q_{jj}=q,&&\qti_{ik}=\qti_{i\ell}=\qti_{j\ell}=1.\qquad
\end{align*}
If $[[x_{ijk\ell},x_j]_c,x_k]_c-q_{jk}(q^2-q)[[x_{ijk\ell},x_k]_c,x_j]_c\in\Pc(\toba)$, then 
$$[[x_{ijk\ell},x_j]_c,x_k]_c-q_{jk}(q^2-q)[[x_{ijk\ell},x_k]_c,x_j]_c=0.$$
\end{lemma}
\pf
Suppose that $x_u:=[[x_{ijk\ell},x_j]_c,x_k]_c-q_{jk}(q^2-q)[[x_{ijk\ell},x_k]_c,x_j]_c\neq 0$. We have that
$\qti_{ui}=\qti_{uj}=\qti_{u\ell}=1$ and $\qti_{uk}=q_{uu}=q$. Then the subdiagram with vertices $i,j,k,\ell,u$ is
\begin{align*}
\xymatrix@C=30pt@R-5pt{&&\overset{q}{\underset{u}{\circ}} &\\
\overset{q^2}{\underset{i}{\circ}} \ar@{-}[r]^{q^{-2}} & \overset{q}{\underset{j}{\circ}} \ar@{-}[r]^{q^{-1}} & 
\overset{-1}{\underset{k}{\circ}} \ar@{-}[r]^{q^3} \ar@{-}[u]_{q} & \overset{q^{-3}}{\underset{\ell}{\circ}}.}
\end{align*}

If $q\notin \G_6'$, then this diagram has a unique vertex labelled with $-1$ and has degree three. If $q\in\G_6'$, then $-c_{uk}=5$. In any case case this diagram does not belong to \cite[Table 3]{H-classif} and we get a contradiction. Thus $x_u=0$.
\epf

\begin{lemma}\label{lem:diagonal-[[[xijk,xj],xj],xj]}
Let $i,j,k\in \I_{\theta}$ such that $q_{jj}=\qti_{ij}^3=\qti_{jk}\in \G_4'$ and $\qti_{ik}=1$. If either $q_{ii}\neq -1$ or else $q_{kk}\neq -1$, then $[[[x_{ijk},x_j]_c,x_j]_c,x_j]_c=0$.
\end{lemma}
\pf
First we check that $x_u:=[[[x_{ijk},x_j]_c,x_j]_c,x_j]_c \in\Pc(\toba)$. Using previous results,
\begin{itemize}[leftmargin=*]
\item $x_{ik}=0$ in $\toba$ by Lemma \ref{lem:ij}.
\item $x_{jji}=0$ in $\toba$: this follows either by Lemma \ref{lem:qs} if $q_{ii}\ne q_{jj}$ or by Lemma \ref{lem:dqs2} if $q_{ii}=q_{jj}$.
\item $x_{jjjjk}=0$ in $\toba$ applying Lemma \ref{lem:QSR-orden-bajo}: here $-c_{jk}=3$.
\end{itemize}
Let $q:=q_{jj}\in\G_4'$, $z=[x_{ijk},x_j]_c$, $y=[[x_{ijk},x_j]_c,x_j]_c$. 
Using the relations $x_{ik}=0$ and $x_{ij}x_j=qq_{ij}x_jx_{ij}$, we compute recursively the following coproducts:
\begin{align*}
\Delta (x_{jk})=&x_{jk}\ot 1+(1-\qti_{jk})x_j\ot x_k +1\ot x_{jk},
\\
\Delta  (x_{ijk})=& x_{ijk}\ot 1+ (1-q)x_{ij}\ot x_k+ (1-q^{-1})x_i\ot x_{jk}+1\ot x_{ijk},
\\
\Delta (z) =& z \ot 1 + 1\ot z +2x_{ijk}\ot x_j
+ (q-1)q_{kj}x_{ij}\ot x_{jk}
\\ & 
+ (1-q)x_{ij}\ot (2x_kx_j-q_{kj}q x_{jk})
+ (1+q)x_i\ot (2 x_{jk}x_j-q_{kj}x_{jjk}),
\\
\Delta  (y) =& y \ot 1 +1\ot y +2qz\ot x_j +2(1+q)x_{ijk}\ot x_j^2
\\ 
& +q_{kj}x_{ij} \ot \big((q-1)q_{kj}x_{jjk}-4q x_{jk}x_j -4qq_{jk}x_kx_j^2\big)
\\ 
& + (1+q)x_i\ot ( 2(1+q) x_{jk}x_j^2-2qq_{kj} x_{jjk}x_j +qq_{kj}^2x_{jjjk} ),
\\
\Delta  (x_u) = &x_u \ot 1 +1\ot x_u+(q-1)q_{kj}^3 x_i\ot x_{jjjjk}.
\end{align*}
Since $x_{jjjjk}=0$ in $\toba$, it follows that $x_u\in\Pc(\toba)$.
Suppose that $x_u\neq 0$. As $\qti_{ui}=q_{ii}^2$, $\qti_{uk}=q_{kk}^2$, $\qti_{uj}=1$, $q_{uu}=q_{ii}q_{kk}$, the Dynkin subdiagram of $\bq'$ with vertices $i,j,k,u$ is
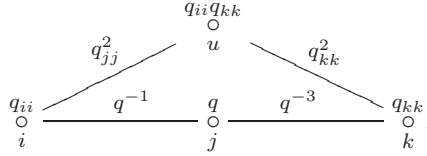
\begin{align*}
\xymatrix@C50pt@R-15pt{&\overset{q_{ii}q_{kk}}{\underset{u}{\circ}}\ar@{-}[dl]_{q_{jj}^2}\ar@{-}[dr]^{q_{kk}^2}&\\
\overset{q_{ii}}{\underset{i}{\circ}} \ar@{-}[r]^{\quad q^{-1}}&\overset{q}{\underset{j}{\circ}} \ar@{-}[r]^{q^{-3}\quad}& \overset{q_{kk}}{\underset{k}{\circ}}.}
\end{align*}
If $q_{ii}\neq -1 \neq q_{kk}$, then  the diagram above is a $4$-cycle so $\GK\toba=\infty$ by Lemma \ref{lem:cycles}. Otherwise, either $q_{ii}=-1$ or $q_{kk}=-1$, and the other vertex $\ell$ has label $q_{\ell \ell}\neq -1$. The subdiagram corresponding to $\ell,u$ is
$\xymatrix@C=30pt{\overset{q_{\ell\ell}}{{\circ}} \ar@{-}[r]^{q_{\ell\ell}^{2}}&\overset{-q_{\ell\ell}}{{\circ}}}$, so $\GK\toba=\infty$ by Remark \ref{lem:q-q2-(-q)}. In both cases we get a contradiction, so $x_u=0$.
\epf

\begin{lemma}\label{lem:diagonal-[[x_ij,[x_ij,x_ijk]],x_j]} 
Let $i,j,k \in \I_{\theta}$ be such that $q_{ii}=q_{jj}=-1$, $\qti_{ik}=1$, $\qti_{ij}=q$, $q_{kk}^{-1}=q^{-3}=\qti_{jk}$ for some $q\in\Bbbk-(\G_2\cup\G_3)$. If $x_u:=[[x_{ij},[x_{ij},x_{ijk}]_c]_c,x_j]_c\in\Pc(\toba)$, then $x_u=0$.
\end{lemma}
\pf
Suppose on the contrary that $x_u\neq 0$. As $\qti_{ui}=q^4$, $\qti_{uj}=1$, $\qti_{uk}=q^{-6}$, $q_{uu}=-q^3$,
the Dynkin subdiagram with vertices $i,j,k,u$ is
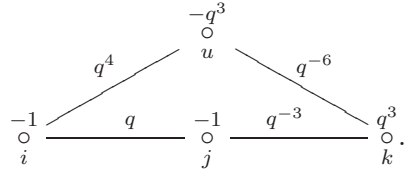
\begin{align*}
\xymatrix@C50pt@R-15pt{&\overset{-q^3}{\underset{u}{\circ}}\ar@{-}[ld]_{q^4}\ar@{-}[rd]^{q^{-6}}&\\
\overset{-1}{\underset{i}{\circ}} \ar@{-}[r]^{\quad q} & \overset{-1}{\underset{j}{\circ}} \ar@{-}[r]^{q^{-3}\quad} & 
\overset{q^3}{\underset{k}{\circ}}.}
\end{align*}
If $q\in\G_6'$, then the subdiagram with vertices $i,u$ is $\xymatrix@C=30pt{\overset{-1}{{\circ}} \ar@{-}[r]^{q^{4}} & \overset{1}{{\circ}}}$, so $\GK \toba = \infty$ by Lemma \ref{lem:1connected} .
Otherwise the subdiagram with vertices $k,u$ is $\xymatrix@C=30pt{\overset{q^3}{{\circ}} \ar@{-}[r]^{q^{-6}} & \overset{-q^3}{{\circ}}}$ where $q^{3}\ne \pm 1$; now Remark \ref{lem:q-q2-(-q)} gives $\GK \toba = \infty$.
In any case we get a contradiction, so $x_u=0$.
\epf

\begin{lemma}\label{lem:diagonal-[x_i,[x_ijk,x_j]]} 
Let $i,j,k \in \I_{\theta}$ be such that $q_{ii}=-q^{-1}$, $\qti_{ik}=1$, $q_{jj}=-1$, $\qti_{ij}=q^2$, $\qti_{jk}=q^{-3}=q_{kk}^{-1}$, for some $q\in\Bbbk-(\G_2\cup\G_3)$.

If $x_u :=[x_i,[x_{ijk},x_j]_c]_c -\frac{q_{ij}q_{kj}}{1+q}[x_{ij},x_{ijk}]_c +(q^{-1}-q^{-2})q_{ij}q_{ik}x_{ijk}x_{ij}\in\Pc(\toba)$, then $x_u=0$.
\end{lemma}
\pf
Suppose on the contrary that $x_u\ne 0$. As $\qti_{ui}=\qti_{uk}=1$, $\qti_{uj}=q_{uu}=q$, the subdiagram with vertices $i,j,k,u$ is
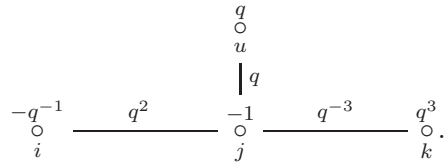
\begin{align*}
\xymatrix@R=10pt{&&\overset{q}{\underset{u}{\circ}} \ar@{-}[d]^{q}&&\\
\overset{-q^{-1}}{\underset{i}{\circ}} \ar@{-}[rr]^{q^2} && 
\overset{-1}{\underset{j}{\circ}} \ar@{-}[rr]^{q^{-3}} && 
\overset{q^3}{\underset{k}{\circ}}.  \\
}
\end{align*}
As either $-c_{uj}=5$ if $q\in\G_6'$, or the unique vertex with label $-1$ has degree 3, this diagram does not belong to \cite[Table 3]{H-classif} and we get a contradiction with Conjecture \ref{conj:AAH}.
\epf

\section{Eminent pre-Nichols algebras of Cartan type $G_2$} \label{section:Cartan-typeG2}

Consider a braiding  $\bq$ of Cartan type $G_2$, so the Dynkin diagram is $\xymatrix @C=15pt{ \underset{ 1 }{\overset{q}{\circ}} \ar  @{-}[rr]^{q^{-3}}  & & \underset{ 2 }{\overset{ q^3}{\circ} } }$, where $q$ is a root of unity of order $N>3$.  If $N\ne 4, 6$, then the distinguished pre-Nichols algebra of  $\bq$ is eminent  by \cite[Lemma 4.13]{ASa}. Here we extend the result for $N=4, 6$. In order to do that, we give first a minimal presentation of these Nichols algebras.

\subsection{Minimal presentation of $\toba_{\bq}$ when $N=4$ or $6$}\label{subsection:CartanG2-N=4,6-minimal}
Let $x_{11212} :=[x_{112}, x_{12}]_c$. By \cite[Theorem 3.1]{An-crelle}, the algebra $\toba_{\bq}$ is presented by generators $x_1, x_2$ and relations 
\begin{align}\label{eq:CartanG2-N=4,6-presentation}
\begin{aligned}
&x_{11112}, && x_{221}, &&  x_{\alpha}^{N_\alpha},  \ \alpha \in \varDelta^+,\\
&[x_{1112}, x_{112}]_c, && [x_{112}, x_{11212}]_c, && [x_{11212}, x_{12}]_c, &&[x_1,x_{11212}]_c-q_{12}\frac{q^4-q}{q+1} x_{112}^2.
\end{aligned}
\end{align}
Compare this presentation with \cite[Theorem 5.25]{standard} and \cite{AA17}.

\begin{lemma}\label{lem:CartanG2-N=4-minimal}
Assume that $\bq$ is of Cartan type $G_2$ and $N=4$. Then $\toba_{\bq}$ is minimally presented by generators $x_1, x_2$ and relations 
\begin{align*}
x_{221}, && [x_{1112}, x_{112}]_c, && x_{\alpha}^{N_\alpha}, \ \alpha \in \varDelta^+.
\end{align*}
\end{lemma}
\pf
Let $\toba$ be the algebra generated by $x_1, x_2$ subject to these relations. Since $q\in \G_4'$, we have $x_1^4=0$, hence $x_{11112}=0$ in $\toba$. We verify with \texttt{GAP} that relations $[x_{112}, x_{11212}]_c$, $[x_1,x_{11212}]_c-q_{12}(q^4-q)(2)^{-1}_{q} x_{112}^2$ and $[x_{11212}, x_{12}]_c$ also vanish in $\toba$. Hence $\toba=\toba_{\bq}$.

By a degree argument, to verify the minimality of this presentation it is enough to prove that $[x_{1112}, x_{112}]_c$ does not vanish in the algebra presented by the relations $x_{221}$, $x_{1}^{4}$,  $x_{2}^{4}$, which is checked using \texttt{GAP}.
\epf

\begin{lemma}\label{lem:CartanG2-N=6-minimal}
Assume that $\bq$ is of Cartan type $G_2$ and $N=6$. Then $\toba_{\bq}$ is minimally presented by generators $x_1, x_2$ and relations 
\begin{align*}
x_{11112}, && x_{221}, && [x_{11212}, x_{12}]_c, && x_{\alpha}^{N_\alpha}, \ \alpha \in \varDelta^+.
\end{align*}
\end{lemma}
\pf
We argue as in Lemma \ref{lem:CartanG2-N=4-minimal}.
\epf

\subsection{Eminent pre-Nichols algebras of Cartan type $G_2$}\label{subsubsection:CartanG2-N=4-eminent}

The goal now is to prove that the distinguished pre-Nichols algebra is eminent. By \cite[Lemma 4.13]{ASa} we already know that relations $x_{11112}=0$ and $x_{221}=0$ hold in any finite $\GK$ pre-Nichols algebra.

\begin{lemma}
Let $\toba$ be a pre-Nichols algebra of $\bq$. The following hold:
\begin{enumerate}[leftmargin=*,label=\rm{(\alph*)}]
\item\label{item:CartanG2-N=4-qsr-computations} If  $ x_{11112}= 0 = x_{221}$ hold in $\toba$, then the following relations also hold:
\begin{align}
\label{eq:CartanG2-N=4-qsr-[x12,x2]}
[x_{12},x_2]_c&=0
\\
\label{eq:CartanG2-N=4-qsr-[x112,x2]}
[x_{112},x_{2}]_c&= q_{12} q (q^2-1)x_{12}^2 
\\
\label{eq:CartanG2-N=4-qsr-[x1112,x2]}
[x_{1112},x_{2}]_c&= q_{12} q (q^2-q-1) x_{11212}+ q_{12}^2 q^2 (q^3-1)x_{12} x_{112}
\\
\label{eq:CartanG2-N=4-qsr-[x1112,x12]}
[x_{1112},x_{12}]_c&=[x_{1},x_{11212}]_c=q_{12} (q^3-1) (2)_{q^{-1}}^{-1} x_{12}^2
\\
\label{eq:CartanG2-N=4-qsr-[x1^3,x112]}
[x_{1}^3,x_{112}]_c&= q_{12}^2 q ^4 (3)_q x_{1112} x_1^2 
\end{align}
\end{enumerate}
\end{lemma}

\pf
Follow from rutinary computations.
\epf

Assume now $N=4$. By  Lemma \ref{lem:CartanG2-N=4-minimal}, the distinguished pre-Nichols algebra of $\bq$ is 
$$\wtoba_{\bq}=T(V)/\langle x_{11112}, x_{221}, [x_{1112}, x_{112}]_c \rangle;$$

\begin{lemma}\label{lem:CartanG2-N=4}
Assume $N=4$ and let $\toba$ be a pre-Nichols algebra of $\bq$. The following hold:
\begin{enumerate}[leftmargin=*,label=\rm{(\alph*)}]
\item\label{item:CartanG2-N=4-[x1112,x112]-primitive} If  $ x_{11112}= 0 = x_{221}$ hold in $\toba$, then  $[x_{1112}, x_{112}]_c$ is primitive in $\toba$.
\item\label{item:CartanG2-N=4-[x1112,x112]=0} If $\GK (\toba) < \infty$, then $[x_{1112}, x_{112}]_c= 0$ hold in $\toba$. Hence $\wtoba_{\bq}$ is eminent.
\end{enumerate}
\end{lemma}

\pf
\ref{item:CartanG2-N=4-[x1112,x112]-primitive} Using \cite[Lemma 4.23]{H-Lusztig-isos} and that the coproduct is braided multiplicative,
\begin{align*}
&\Delta([x_{1112}, x_{112}]_c) =x_{1112}x_{112} \ot 1 
+ \sum_{u=0}^{2} \binom{2}{u}_q \prod_{t=1}^{u} (1-q^{-1-t}) x_{1112} x_1^u \ot (\ad_c x_1)^{2-u}x_2
\\
&+\sum_{r=0}^{3} \binom{3}{r}_q \prod_{s=1}^{r} (1-q^{-s}) q_{11}^{2(3-r)} q_{12}^{3-r}q_{21}^2q_{22} x_1^r x_{112}\ot (\ad_cx_1)^{3-r}x_2 
\\
&+\sum_{\substack{0\leq r \leq 3\\ 0\leq u \leq 2}} \binom{3}{r}_q  \binom{2}{u}_q \prod_{\substack{1\leq s \leq r\\ 1\leq t \leq u}} (q^{-u}-q^{-s})  (q^3-q^{-1-t})  q_{21}^u x_1^{r+u} \ot (\ad_cx_1)^{3-r}x_2  (\ad_cx_1)^{2-u}x_2 
\\
&-q^3q_{12}  x_{112}x_{1112} \ot 1 
- q^3q_{12}\sum_{r=0}^{3} \binom{3}{r}_q \prod_{s=1}^{r} (1-q^{-s}) x_{112}x_1^r \ot (\ad_cx_1)^{3-r}x_2 
\\
&-q^3 q_{12} \sum_{u=0}^{2} \binom{2}{u}_q \prod_{t=1}^{u} (1-q^{-1-t}) q_{11}^{3(2-u)}q_{12}^{2-u} q_{21}^3 q_{22} x_1^u x_{1112} \ot (\ad_c x_1)^{2-u}x_2
\\
&-\sum_{\substack{0\leq r \leq 3\\ 0\leq u \leq 2}} 
\binom{3}{r}_q  \binom{2}{u}_q
\prod_{\substack{1\leq s \leq r\\ 1\leq t \leq u}} (q^2-q^{-s})  (q^{-r}-q^{-1-t})  q_{21}^{r-1}
x_1^{r+u} \ot (\ad_cx_1)^{2-u}x_2   (\ad_cx_1)^{3-r}x_2.
\end{align*}
Notice that the terms on  the homogeneous components $\toba_{(5,2)} \ot \ku $ and $\ku\ot\toba_{(5,2)} $ are $[x_{1112},x_{112}]_c \ot 1$ and $1\ot [x_{1112},x_{112}]_c$, respectively. Next we show that the other homogeneous components vanish.

\noindent Component  in $\toba_{(3,1)} \ot \toba_{(2,1)}$:
\begin{align*}
(x_{1112} + (1-q^{-1})q^2 (x_1x_{112} - q^2 q_{12} x_{112}x_1) - q^{12}\qti_{12}^3 x_{1112} ) \ot x_{112} = 0.
\end{align*}

\noindent Component  in $\toba_{(4,1)} \ot \toba_{(1,1)}$:
\begin{align*}
& ( (2)_q (1-q^{-2})  x_{1112}x_1 + (3)_q(1-q^{-1})(1-q^{-2}) q^2q_{21} x_1^2x_{112}
\\
& - q^3q_{12} (3)_q (1-q^{-1}) (1-q^{-2}) x_{112} x_1^2 - q^ 3q_{12}(2)_q(1-q^{-2}) q^3 q_{21}^2 x_1x_{1112}) \ot  x_{12}\\
=&(2)_q (1-q^{-2}) ( x_{1} x_{112} x_1 - q_{11}^2q_{12} x_{112} x_1^2 + q^3 q_{21} x_1^2x_{112} \\
&- q_{12} x_{112} x_1^2 -q_{21}q^3 x_1^2 x_{112} + q^2 x_{1} x_{112} x_1) \ot x_{12} = 0. 
\end{align*}

\noindent Component  in $\toba_{(5,1)} \ot \toba_{(0,1)}$: 
\begin{align*}
&(1-q^{-2})(1-q^{-3}) (x_{1112}x_1^2 - q^3 q_{12} x_1^2 x_{1112} + (1-q^{-1})q_{21}^2 q_{22} [x_{1}^3, x_{112}]_c) \ot x_2\\
=&(1-q^{-2})(1-q^{-3}) (1-q^3 +(1-q^{-1})q^2) x_{1112}x_1^2 \ot x_2 = 0,
\end{align*}
where we are using that $[x_1^2, x_{1112}]_c=0$ and \eqref{eq:CartanG2-N=4-qsr-[x1^3,x112]}.

\noindent Component  in $\toba_{(2,1)} \ot \toba_{(3,1)}$: 
$(q_{11}^6 q_{12}^3 q_{21}^2 q_{22} - q^3 q_{12}) x_{112} \ot x_{1112} = 0$.

\noindent Component  in $\toba_{(1,0)} \ot \toba_{(4,2)}$: using \eqref{eq:CartanG2-N=4-qsr-[x1112,x12]},
\begin{align*}
&(1-q^{-2}) x_1\ot ( (1-q^{-1}) (3)_q x_{112}^2 + (2)_q q_{11}^3 q_{21} [x_{1112},x_{12}]_c  ) \\
=&(1-q^{-2}) x_1\ot  ((3)_q (1-q^{-1}) + (2)_q q^3 ) x_{112}^2 = 0.
\end{align*}

\noindent Component  in $\toba_{(2,0)} \ot \toba_{(3,2)}$: by definition of $x_{11212}$ and \eqref{eq:CartanG2-N=4-qsr-[x1112,x2]} we get,
\begin{align*}
&(1-q^{-2})qx_1^2 \ot (q^2(2)^2_q q_{21} x_{11212}  + (2)_q(x_{12 }x_{112} - (2)_q q_{21} x_{112}x_{12} + (1-q) q q_{21}^2 [x_{1112}, x_2]_c) \\
=&(1-q^{-2})qx_1^2 \ot (( q^2q_{21}(2)_q^2 +(1-q)q^3q_{21}(q^2-q-1))  x_{11212}  \\
& \qquad \qquad \qquad  \qquad + ( (2)_q + q(1-q) (q^3-1)) x_{12} x_{112} - (2)_q q_{21} x_{112}x_{12}) \\
=&(1-q^{-2})q(2)_q  x_1^2 \ot (q_{21} x_{11212} - q x_{12}x_{112} - q_{21}x_{112}x_{12}) = 0. 
\end{align*}

\noindent Component  in $\toba_{(3,0)} \ot \toba_{(2,2)}$: using 
\eqref{eq:CartanG2-N=4-qsr-[x112,x2]},
\begin{align*}
&(1-q^{-1})(1-q^{-2})(1-q) x_1^3 \ot (-q(1+q)x_2 x_{112} + q_{21}^2(1+q) x_{112}x_2) \\
&+(1-q^{-1})(1-q^{-2})^2(2)_q q x_1^3 \ot ( q_{11}q_{21} x_{12}^2 -q_{11}^2q_{21}x_{12}^2) =0.
\end{align*}

\noindent Component  in $\toba_{(4,0)} \ot \toba_{(1,2)}$: by 
\eqref{eq:CartanG2-N=4-qsr-[x12,x2]},
\begin{align*}
&(1-q^{-1})(1-q^{-2})^2(1-q^{-3}) x_1^4 \ot \\
&\qquad \qquad \qquad \qquad \qquad ((2)_q q_{21} x_2 x_12 - (2)_q q^3 q_{21}^2 x_{12} x_2 + (3)_q q^2 q_{21}^2x_{12}x_2 - (3)_q q_{21}x_{2}x_{12}) \\
=&-(1-q^{-1})(1-q^{-2})^2(1-q^{-3}) q_{21}^2 x_1^4 \ot [x_{12}, x_2]_c = 0.
\end{align*}
\noindent Component  in $\toba_{(5,0)} \ot \toba_{(0,2)}$:
$(1-q^{-1})(1-q^{-2})^2(1-q^{-3})^2x_1^5 \ot (q_{21}^2  -q^3 q_{12}q_{21}^3) x_2^2=0$.

\ref{item:CartanG2-N=4-[x1112,x112]=0} Assume $[x_{1112}, x_{112}]_c \neq 0$. By \ref{item:CartanG2-N=4-[x1112,x112]-primitive} we get $\ku x_1\oplus \ku x_{2} \oplus \ku [x_{1112}, x_{112}]_c \subset \Pc(\toba)$, where the braiding has Dynkin diagram
$$ \xymatrix @C=15pt{ \underset{ 1 }{\overset{q}{\circ}} \ar  @{-}[rr]^{q^{-3}}  & & \underset{ 2 }{\overset{ q^3}{\circ} } 
\ar  @{-}[rr]^{q^{-3}}  & & \underset{ 1112112 }{\overset{q^3}{\circ} }},$$
which is of affine Cartan type $D_4^{(3)}$. By \cite[Theorem 1.2 (a)]{AAH-diag} the Nichols algebra with this braiding has infinite $\GK$, and by \cite[Lemma 2.7]{ASa} it follows $\GK \toba = \infty$.
\epf

\medbreak
Assume now $N=6$. By  Lemma \ref{lem:CartanG2-N=4-minimal}, the distinguished pre-Nichols algebra of $\bq$ is 
$$\wtoba_{\bq}=T(V)/\langle x_{11112}, x_{221}, [x_{11212}, x_{12}]_c \rangle;$$

\begin{lemma}\label{lem:CartanG2-N=6}
Assume $N=6$ and let $\toba$ be a pre-Nichols algebra of $\bq$. The following hold:
\begin{enumerate}[leftmargin=*,label=\rm{(\alph*)}]
\item\label{item:CartanG2-N=6-[x11212,x12]-primitive} If  $ x_{11112}= 0 = x_{221}$ hold in $\toba$, then  $[x_{11212}, x_{12}]_c$ is primitive in $\toba$.
\item\label{item:CartanG2-N=6-[x11212,x12]=0} If $\GK (\toba) < \infty$, then $[x_{11212}, x_{12}]_c= 0$ hold in $\toba$. Hence $\wtoba_{\bq}$ is eminent.
\end{enumerate}
\end{lemma}

\pf
\ref{item:CartanG2-N=6-[x11212,x12]-primitive} We compute $\Delta(x_{11212})$ in the tensor algebra and get
\begin{align*}
\Delta(x_{11212})=&x_{11212}\ot 1 + 1 \ot x_{11212} + 2q^2 x_{112} \ot x_{12} + 2q^2 q_{21} x_1 \ot [x_{112},x_2]_c 
\\
&+2 [x_{112},x_1]_c \ot x_2 + 2(1-q^{-2})q_{21}q^3 [x_1, x_{12}]_c \ot x_2 + 
\\
&+ 4(1-q^{-1})(1-q^{-2})q_{21} x_1^3 \ot x_2^2 + 2(q^2-1) x_1^2 \ot [x_2,x_{12}]_c.
\end{align*}
Using this and $\Delta(x_{12})=x_{12}\ot 1 + 1\ot x_{12} + 2x_1\ot x_2$, we compute $\Delta([x_{11212}, x_{12}]_c)$ in $\toba$. The components of degrees $\toba_{(4,3)} \ot \ku $ and $\ku\ot\toba_{(4,3)} $ are $[x_{11212},x_{12}]_c \ot 1$ and $1\ot [x_{11212},x_{12}]_c$, respectively. The other homogeneous components vanish, most of them due to commutations between powers of $x_2$ and powers of $x_{12}$, which follow from \eqref{eq:CartanG2-N=4-qsr-[x12,x2]}.

\ref{item:CartanG2-N=6-[x11212,x12]=0} If $[x_{11212}, x_{12}]_c \neq 0$ in $\toba$, then by \ref{item:CartanG2-N=6-[x11212,x12]-primitive} we get $\ku x_1\oplus \ku x_{2} \oplus \ku [x_{11212}, x_{12}]_c \subset \Pc(\toba)$, where the braiding has Dynkin diagram
$$ \xymatrix @C=15pt{ \underset{ 1121212 }{\overset{q}{\circ} } \ar  @{-}[rr]^{q^{-1}}  & & \underset{ 1 }{\overset{q}{\circ}} \ar  @{-}[rr]^{q^{-3}}  & & \underset{ 2 }{\overset{ q^3}{\circ} } },$$
which is of indefinite Cartan type. By Lemma \ref{lem:subspace-primitives} it follows $\GK \toba = \infty$.
\epf

\begin{theorem}\label{th:CartanG2-eminent}
If $\bq$ is of Cartan type $G_2$, then the distinguished pre-Nichols algebra $\wtoba_{\bq}$ is eminent.
\end{theorem}
\pf
If $N\ne 4,6$ this statement is included in \cite[Theorem 1.3 (a)]{ASa}. The cases $N=4,6$ were treated in Lemmas \ref{lem:CartanG2-N=4} and \ref{lem:CartanG2-N=6}.
\epf

\section{Eminent pre-Nichols algebras of super type} \label{section:super-type}

In this Section we describe case-by-case all eminent pre-Nichols of braiding of super type. We prove that the associated distinguished pre-Nichols algebra is eminent except for two exceptional cases in type A; for these exceptions we construct eminent pre-Nichols algebras and show that they are central extension of the distinguished one.

\subsection{Type A} \label{subsection:type-supera}
Here $\bq$ is a braiding of type $\superqa{\theta}{q}{\J}$, where $q\in \ku^\times$ is a root of unity of order $N>2$ and $\emptyset \ne \J \subseteq \I_{\theta}$ \cite[\S 5.1]{AA17}.  This means that the Dynkin diagram of $\bq$ is
\begin{align*}
\xymatrix{ \overset{q_{11}}{\underset{\ }{\circ}}\ar  @{-}[rr]^{\widetilde{q}_{12}}  &&
\overset{q_{22}}{\underset{\ }{\circ}}\ar@{.}[r] &  \overset{\quad q_{\theta-1 \theta-1}}{\underset{\ }{\circ}} \ar  @{-}[rr]^{\widetilde{q}_{\theta-1 \theta}}  &&
\overset{\quad q_{\theta\theta}}{\underset{\ }{\circ}}}
\end{align*}
where $q_{ii}$, $\qti_{ij}$ satisfy the following conditions:
\begin{enumerate}[leftmargin=*,label=\rm{(\Alph*)}]
\item \label{item:supera-diagram-condition-theta} $q_{\theta\theta}^2\widetilde{q}_{\theta-1 \theta}=q$;

\smallbreak
\item \label{item:supera-diagram-condition-i-in-J}
 if $i\in\J$, then $q_{ii}=-1$, $\widetilde{q}_{i-1 i}=\widetilde{q}_{i i+1}^{-1}$;

\smallbreak
\item \label{item:supera-diagram-condition-i-notin-J}
if $i\notin\J$, then $\widetilde{q}_{i-1\, i}= q_{ii}^{-1} = \widetilde{q}_{i\, i+1}$ (as long as $i\pm 1\in\I_{\theta}$).
\end{enumerate}
Given $q$ and $\J$, these conditions determine the diagram. One can deduce that:
\begin{itemize} [leftmargin=*]
\item $q_{ii}=q^{\pm1}\ne \pm1$ for $i\in \I-\J$, and $\qti_{ii+1}=q^{\pm1}\ne \pm1$ for all $i<\theta$;
\item if $i\in\J$ but $i-1,i+1\notin \J$, then $q_{i-1i-1}q_{i+1i+1}=1$.
\end{itemize} 
 These two remarks will allow us to apply Lemmas \ref{lem:diagonal-[x_ijk,x_j]} and \ref{lem:cor2}.

\medbreak
The distinguished pre-Nichols  algebra is  presented by generators $(x_i)_{i\in \I}$ and relations
\begin{align}\label{eq:Asuper-preNichols-rels}
\begin{aligned}
& x_{ij}=0, \quad i < j - 1; & x_{ii(i\pm1)} &= 0, \quad
q_{ii}\neq-1;  \\
&x_i^2=0, \quad q_{ii}=-1; & [x_{(i-1i+1)},x_i]_c&=0, \quad q_{ii}=-1.
\end{aligned}
\end{align}

\medbreak
Let $\toba$ be a pre-Nichols algebra of $\bq$ such that $\GK\toba<\infty$.
Next we prove that each relation in \eqref{eq:Asuper-preNichols-rels} must hold in $\toba$, except for two pairs $(\theta,\J)$.

\begin{lemma}\label{lem:ijsupera}
Assume that $\bq$ is not of type $\superqa{3}{q}{\{1,2,3\}}$. If $i,j\in \I_{\theta}$, $i<j-1$, then $x_{ij}=0$ in $\toba$. 
\end{lemma}
\pf
Suppose on the contrary that $x_{ij}\neq 0$ in $\toba$. By Lemma \ref{lem:ij}, $q_{ii}=q_{jj}=-1$, that is $i,j\in \J$. Also, $i=1$ by Lemma \ref{lem:ij3}: otherwise $i-1\in\I$ and $\qti_{ii-1}\qti_{ji-1}=\qti_{ii-1}\neq 1$. Analogously $j=\theta$, and $i=j-2$: if $i<j-2$, then $\qti_{ij-1}\qti_{jj-1}=\qti_{jj-1}\neq 1$. 

Hence $i=1$, $j=3=\theta$ and $i,j\in \J$. Thus $2\notin \J$ since $\bq$ is not of type $\superqa{3}{q}{\{1,2,3\}}$. But then $\qti_{12}\qti_{32}=q^{2}\neq 1$, a contradiction with Lemma \ref{lem:ij3}. Therefore $x_{ij}=0$. 
\epf

\begin{lemma}\label{lem:qssupera}
If $i\in \I_{\theta}-\J$, then $x_{ii(i\pm 1)}=0$ in $\toba$.
\end{lemma}
\pf
Let $i\in \I_{\theta-1}$ be such that $i\not \in \J$, and set $j=i+1$. Suppose that $x_{iij}\neq 0$. By Lemmas \ref{lem:qs} and \ref{lem:dqs},
$q_{ii}=q_{jj}=\qti_{ij}^{-1}\in \G_3'$. 
If $i>1$, then $\qti_{(i-1)i}\neq \pm 1$, $\qti_{(i-1)j}=1$, and we get a contradiction with Lemma \ref{lem:dqs2}. If $j<\theta$, then $\qti_{j(j+1)}\neq \pm 1$, $\qti_{i(j+1)}=1$, and again we get a contradiction with Lemma \ref{lem:dqs2}. Otherwise $i=1$, $\theta=2$ and $\J=\emptyset$, a contradiction. Therefore $x_{iij}=0$. 
The proof of $x_{ii(i-1)}=0$ for $i\in \I_{2,\theta}-\J$ is analogous.
\epf

\begin{lemma}\label{lem:corsupera}
If $\bq$ is not of type $\superqa{3}{q}{\{2\}}$ and $i\in \I_{2,\theta-1}\cap \J$, then $[x_{(i-1,i+1)},x_i]_c=0$.
\end{lemma}
\pf
Let $i\in\I_{2,\theta-1}\cap \J$. 
If either $i-1 \notin \J$ or $i+1 \notin \J$, then $[x_{(i-1,i+1)},x_i]_c= 0$ by Lemma \ref{lem:diagonal-[x_ijk,x_j]}.
If $2<i$ (respectively $i+1<\theta$), then $\ell=i-2$ (respectively $\ell=i+2$) satisfies \ref{item:cor2-a} (respectively \ref{item:cor2-b}) of Lemma \ref{lem:cor2}, hence $[x_{(i-1,i+1)},x_i]_c=0$. Otherwise $i=2$, $\theta=3$, $\J=\{2\}$, so $\bq$ is of type $\superqa{3}{q}{\{2\}}$.
\epf

\begin{theorem}\label{thm:supera-distinguished-eminent}
Let $\bq$ be of type $\superqa{\theta}{q}{\J}$, where the pair $(\theta,\J)$ is not one of the following:
\begin{enumerate}[leftmargin=*,label=\rm{(\Roman*)}]
\item\label{item:supera-exception-J=2} $\theta=3$,  $\J=\{2\}$.
\item\label{item:supera-exception-J=I3} $\theta=3$,  $\J=\{1,2,3\}$.
\end{enumerate}
Then the distinguished pre-Nichols algebra $\wtoba_{\bq}$ is eminent.
\end{theorem}
\pf
By Lemmas \ref{lem:xi-no-Cartan}, \ref{lem:ijsupera}, \ref{lem:qssupera} and \ref{lem:corsupera} the defining relations of $\wtoba_{\bq}$ hold in any finite $\GK$ pre-Nichols algebra $\toba$ of $\bq$, hence $\wtoba_{\bq}$ projects onto $\toba$. As $\GK \wtoba_{\bq}<\infty$, it follows that $\wtoba_{\bq}$ is eminent.
\epf

%

Next we find eminent pre-Nichols algebras for braidings of type \ref{item:supera-exception-J=2} and \ref{item:supera-exception-J=I3}. 

\subsubsection{Type $\superqa{3}{q}{\{2\}}$} \label{subsubsec:supera-exception-J=2}
In this case the distinguished pre-Nichols algebra is
$$\wtoba_{\bq}=T(V)/\langle x_2^2,x_{13},x_{112},x_{223},[x_{123},x_2]_c\rangle.$$
We consider the following algebra:
\begin{align}\label{eq:eminent-supera-caseII}
\htoba_{\bq}=T(V)/\langle x_2^2,x_{13},x_{112},x_{332} \rangle,
\end{align}
which is a braided Hopf algebra since $x_2^2$, $x_{13}$, $x_{112}$, and $x_{332}$ are primitive in $T(V)$.
Also the projections from $T(V)$ induce a surjective map of braided Hopf algebras $\pi: \htoba_{\bq} \twoheadrightarrow \wtoba_{\bq}$.
Next we prove that the pre-Nichols algebra $\htoba_{\bq}$ is eminent.

\begin{prop}\label{prop:eminent-supera-caseII}
If $\bq$ of type $\superqa{3}{q}{\{2\}}$, then $\htoba_{\bq}$ as in \eqref{eq:eminent-supera-caseII} an eminent pre-Nichols algebra of $\bq$. 

Let $x_u=[x_{123},x_2]_c$. The set
\begin{align}\label{eq:eminent-supera-caseII-basePBW}
B=\big\{x_3^ax_{23}^bx_2^cx_u^dx_{123}^ex_{12}^fx_1^g: \, b,c,e,f \in \{0,1\}, \, a,d,g\in \N_{0} \big\}
\end{align}
is a basis of $\htoba_{\bq}$, so $\GK \htoba_{\bq}=3$.
\end{prop}
\pf
We split the proof in steps. Let $\toba$ be a finite $\GK$ pre-Nichols algebra of $\bq$.

\begin{step}
The projection $T(V) \twoheadrightarrow \toba$ induces a projection
$\htoba_{\bq} \twoheadrightarrow \toba$ of braided Hopf algebras.
\end{step}

Indeed the defining relations of $\htoba_{\bq}$ annihilate in $\toba$ by Lemmas
\ref{lem:xi-no-Cartan}, \ref{lem:ij} and \ref{lem:qs}.

\smallbreak

To see that $\htoba_{\bq}$ is eminent, it remains to prove that $\GK \htoba_{\bq}<\infty$.
For, we will check that $B$ is a basis of $\htoba_{\bq}$ and that the later is a braided central extension of $\wtoba_{\bq}$.

\begin{step}
The following relations hold in $\htoba_{\bq}$:
\begin{align}\label{eq:eminent-supera-caseII-PRVextra}
x_{12}^2 &=0, & x_{123}^2 &=0, & x_{23}^2 &=0.
\end{align}
\end{step}

Using the relations $x_{112}=0$ and $x_2^2=0$ we compute
\begin{align*}
(1+q_{11})x_{12}^2 &=  (1+q_{11}) (x_1x_2x_1x_2-q_{12}x_2x_1^2x_2+q_{12}^2x_2x_1x_2x_1)
\\
&=q_{11}q_{12}x_2x_1^2x_2-q_{12}(1+q_{11})x_2x_1^2x_2+q_{12}x_2x_1^2x_2=0
\end{align*}
As $q_{11}\neq -1$, we get $x_{12}^2=0$. Analogously $x_{23}^2=0$. 

Using \eqref{eq:braided-commutator-iteration} and the defining relations we check that the following also hold:
\begin{align}\label{eq:eminent-supera-caseII-PBW1}
x_{123}&=x_{12}x_3-q_{13}q_{23}x_3x_{12},
\\ \label{eq:eminent-supera-caseII-PBW2}
x_{1}x_{12}&= q_{11}q_{12}x_{12}x_{1}, & x_1x_3&=q_{13}x_3x_1,
\\
\label{eq:eminent-supera-caseII-PBW3}
x_1x_{123} &=x_1(x_{12}x_3-q_{13}q_{23}x_3x_{12})=q_{11}q_{12}q_{13}x_{123}x_1
\end{align}
Using now that $\adc x_1$ is a skew-derivation and that $(\adc x_1)x_{123}=0$ by \eqref{eq:eminent-supera-caseII-PBW3},
\begin{align*}
0 & = (\adc x_1)^2 (x_{23}^2) = (\adc x_1)(x_{123}x_{23}+q_{12}q_{13}x_{23}x_{123})
=(1+q_{11})q_{12}q_{13}x_{123}^2.
\end{align*}
As $q_{11}\neq -1$, we have that $x_{123}^2=0$.

\begin{step}
$\htoba_{\bq}$ is spanned by $B$.
\end{step}

Let $I$ be the subspace spanned by $B$. It suffices to prove that $I$ is a left ideal of $\htoba_{\bq}$, which in turn follows from the following statement: $x_iI\subset I$ for all $i\in\I_3$. 

We note that $x_3I \subset I$ by definition. Next we check that $x_{23}x_3=q_{23}q_{33}x_3x_{23}$, so $x_{23}I \subset I$. Also, $x_2x_3=x_{23}+q_{23}x_3x_2$, $x_2x_{23}=-q_{23}x_{23}x_2$, so $x_2I \subset I$.

It remains to check that $x_1I \subset I$. Using the defining relations we check that
\begin{align*}
x_{12}x_2 &= -q_{12}x_2x_{12},
&
x_{12}x_{23} &= q_{12}(\qti_{23}-1)x_2x_{123}-q_{12}q_{13}q_{23} x_{23}x_{12}-q_{23}x_u,
\\
x_{123}x_{3} &= -q_{13}q_{23}x_3x_{123},
&
x_{12}x_{123} &=  -q_{13}q_{23}x_{123}x_{12},
\qquad 
x_{123}x_{23} = -q_{12}q_{13}x_{23}x_{123},
\\
x_1x_{u} &= q_{11}q_{12}^2q_{13} x_ux_1,
&
x_ux_3 &= q_{13}q_{23}^2q_{33}x_3x_u,
\qquad
x_ux_2 = q_{12}q_{32}x_2x_u.
\\
x_{12}x_u &= q_{12}q_{13}q_{23}x_{u}x_{12},
&
x_{123}x_u &= -q_{12}q_{22}q_{32}x_ux_{123},
\qquad
x_ux_{23} = q_{12}q_{13}q_{23} x_{23}x_u,
\end{align*}
Using these relations, \eqref{eq:eminent-supera-caseII-PRVextra}, \eqref{eq:eminent-supera-caseII-PBW1}, \eqref{eq:eminent-supera-caseII-PBW2}, \eqref{eq:eminent-supera-caseII-PBW3}, and the definition of the PBW generators we see that $I$ is stable by left multiplication by $x_{12}$, $x_{123}$ and $x_u$, and so by $x_1$.

Next we check that $\htoba_{\bq}$ fits into an exact sequence of braided Hopf algebras, in order to prove that $B$ is linearly independent. Let $\Zc$ be the subalgebra generated by $x_u= [x_{123},x_2]_c$.

\begin{step}
$\Zc$ is a skew-central Hopf subalgebra isomorphic to a polynomial algebra in one variable. 
\end{step}

As $x_u$ is primitive and $c(x_u \ot x_u) = x_u\ot x_u$, $\Zc$ is a Hopf subalgebra of $\htoba_{\bq}$ isomorphic to a polynomial algebra in one variable. It is central since $x_i x_u= q_{i1}q_{i2}^2q_{i3} x_u x_i$ for all $i\in\I_3$.

\begin{step}
There is a degree-preserving extension of braided Hopf algebras $\Zc \hookrightarrow \htoba_{\bq} \twoheadrightarrow \wtoba_{\bq}$.
\end{step}

Let $\Zc'=\htoba_{\bq}^{\co \pi}$. As $\Zc$ is normal (since it is central) and $x_u \in \Pc(\wtoba_{\bq})\cap \ker\pi$, we have that $\Zc \subseteq \Zc'$. By Lemma \ref{lem:extension-braided-graded} and the known PBW basis of $\wtoba_{\bq}$ we have that
\begin{align*}
\mathcal{H}_{\htoba_{\bq}} &= \mathcal{H}_{\wtoba_{\bq}}\mathcal{H}_{\Zc'} \ge \mathcal{H}_{\wtoba_{\bq}}\mathcal{H}_{\Zc} =
\frac{(1+t_2t_3)(1+t_2)(1+t_1t_2t_3)(1+t_1t_2)}{(1-t_3)(1-t_1)}\cdot \frac{1}{1-t_1t_2^2t_3}.
\end{align*}
On the other hand, as  $\htoba_{\bq}$ is spanned by $B$,
\begin{align*}
\mathcal{H}_{\htoba_{\bq}} \le
\frac{(1+t_2t_3)(1+t_2)(1+t_1t_2t_3)(1+t_1t_2)}{(1-t_3)(1-t_1)(1-t_1t_2^2t_3)}.
\end{align*}
Thence the two series above are equal, i.e.
\begin{align}\label{eq:eminent-supera-caseII-Hilbertseries}
\mathcal{H}_{\htoba_{\bq}} =
\frac{(1+t_2t_3)(1+t_2)(1+t_1t_2t_3)(1+t_1t_2)}{(1-t_3)(1-t_1)(1-t_1t_2^2t_3)}.
\end{align}
so $\Zc=\Zc'$, and the claim follows.

\begin{step}
$B$ is a basis of $\htoba_{\bq}$ and $\GK\htoba_{\bq}=3$.
\end{step}
By \eqref{eq:eminent-supera-caseII-Hilbertseries} $B$ is a basis of $\htoba_{\bq}$. Also, the decomposition above of the Hilbert series says that $\GK \htoba_{\bq} =\GK \Zc+\GK \wtoba_{\bq}=3$.
\epf

\begin{remark}
Let $\Zc_{\bq}$ denote the subalgebra of $\htoba_{\bq}$ generated by 
$x_{u}$, $x_{1}^N$ and $x_{3}^N$. One can verify that this is subalgebra is skew-central (more precisely, it is annihilated by the braided adjoint action of $\htoba_{\bq}$), and that it fits in a degree-preserving extension of braided Hopf algebras $ \Zc_{\bq} \hookrightarrow \htoba_{\bq} \twoheadrightarrow \toba_{\bq}$.
\end{remark}

\subsubsection{Type $\superqa{3}{q}{\{1,2,3\}}$} \label{subsubsec:supera-exception-J=I3} 

In this case the distinguished pre-Nichols algebra is given by
$$\wtoba_{\bq}=T(V)/\langle x_1^2,x_2^2, x_3^2, x_{13}, [x_{123},x_2]_c \rangle.$$

\begin{remark}\label{rem:supera-exception-J=I3-x13-substitute}
Let $\toba$ be a pre-Nichols algebra of $\bq$ such that $\GK\toba<\infty$. 
The relations 
\begin{align*}
x_1^2 &=x_2^2=x_3^2=0, & x_{213}&=0, & [x_{123}, x_2]_c&=0
\end{align*}
hold in $\toba$ by Lemmas \ref{lem:xi-no-Cartan}, \ref{lem:diagonal-[x_ijk,x_j]-primitive} \ref{item:diagonal-x_jik=0} and \ref{lem:diagonal-[x_ijk,x_j]}. 
\end{remark}

\begin{remark}\label{rmk:supera-exception-J=I3-definition-eminent}
Let $\htoba_{\bq}$ denote the following quotient of $T(V)$:
\begin{align*}
\htoba_{\bq}=T(V)/\langle x_1^2, x_2^2, x_3^2, x_{213}, [x_{123},x_2]_c \rangle.
\end{align*}
Note that the defining ideal is actually a Hopf ideal, see the proof of Lemma \ref{lem:diagonal-[x_ijk,x_j]-primitive}, so $\htoba_{\bq}$ is a pre-Nichols algebra of $\bq$; next we show that it is eminent.
\end{remark}

\begin{prop}\label{prop:supera-J=I3-eminent} 
The pre-Nichols algebra $\htoba_{\bq}$ is eminent, with $\GK \htoba_{\bq} = 3$ and basis
\begin{align}\label{item:supera-J=I3-eminent-basis}
B=\big\{x_3^ax_{23}^bx_2^cx_{13}^dx_{123}^ex_{12}^fx_1^g: \, a,c,e,g\in \{0,1\}, \, b,d,f \in \N_{0} \big\}.
\end{align}
\end{prop}

\pf  
By Remark \ref{rem:supera-exception-J=I3-x13-substitute}, every finite $\GK$ pre-Nichols algebra of $\bq$ is covered by $\htoba_{\bq}$. The rest of the proof is carried out in several parts.
\begin{stepo}
The following relations hold in $\htoba_{\bq}$:
\begin{align}
\label{eq:supera-J=I3-central-computations-xiij}
x_{1123}&=0; \qquad \qquad x_{iij}=0, \quad  i \ne j;
\\
\label{eq:supera-J=I3-central-computations-[x123,x23]}
[x_{23}, x_{13}]_c&= 0,  \qquad [ x_{123}, x_{23}]_c = 0, \qquad [ x_{123}, x_{13}]_c = 0, \qquad x_{123}^2=0.
\end{align}
\end{stepo}
For \eqref{eq:supera-J=I3-central-computations-xiij}, the relations $x_{iij} = 0$ for $i\ne j$ follow from the condition $q_{ii}=-1$ and the fact that $x_i^2=0$. Similarly,
\begin{align*}
x_{1123} = x_1^2 x_{23} - q_{12}q_{13} (1+q_{11}) x_1 x_{23} x_1 +q_{11}q_{12}^2q_{13}^3 x_{23} x_1^2 = 0.
\end{align*}
For \eqref{eq:supera-J=I3-central-computations-[x123,x23]}, notice that $x_{313}= 0 = x_{213} $ imply that $[x_{23}, x_{13}]_c = 0$; this last equality together with $[x_{23}, x_3]_c=0$ give
\begin{align*}
[x_{123}, x_3]_c= [x_1, [x_{23}, x_3]_c ]_c - q_{12} q_{13} x_{23} x_{13} + q_{23} q_{33} x_{13} x_{23} = q_{23} q_{33}(1-\qti_{12}) x_{13} x_{23}.
\end{align*}
Since $x_{213} = 0 = x_{223}$, it follows $ [ x_2, [x_{123}, x_3]_c]_c =0$. Then
\begin{align*}
[x_{123}, x_{23}]_c&= [[x_{123}, x_2]_c, x_3]_c+ q_{12}q_{22}q_{32}x_2[x_{123}, x_3]_c - q_{23} [x_{123}, x_3]_c x_2\\
&=  q_{12}q_{22}q_{32} [ x_2, [x_{123}, x_3]_c]_c =0.
\end{align*}
Next, use $[x_{23}, x_{13}]_c = 0 = x_{113}$ to get
\begin{align*}
[x_{123}, x_{13}]_c&= [x_1, [x_{23}, x_{13}]_c]_c - q_{12}q_{13} x_{23} x_{113} +q_{21} q_{23}q_{31} q_{33} x_{113} x_{23} = 0.
\end{align*}
Finally, use $x_{123} x_1= - (q_{12}q_{13})^{-1} x_1x_{123} $ and $x_{123} x_{23} = q_{12}q_{13}\qti_{23} x_{23} x_{123}$ to compute
\begin{align*}
x_{123}^2&=x_{123}(x_1x_{23} - q_{12}q_{13} x_{23}x_1) = -\qti_{23} (x_1x_{23} - q_{12}q_{13} x_{23}x_1) x_{123}  = -\qti_{23} x_{123}^2.
\end{align*}
Since $N=\ord \qti_{23} >2$, this implies $x_{123}^2 = 0$.

\begin{stepo}
The set $B$ linearly spans $\htoba_{\bq}$.
\end{stepo}
It is enough to show that the linear span $I$ of $B$ is a left ideal. The inclusion $x_3 L \subset L$ is clear, and $x_2 L \subset L$ follows from the commutation $x_{223} = 0$, cf. \eqref{eq:supera-J=I3-central-computations-xiij}. In order to verify that $x_1 L \subset L$, we argue inductively on $b \geq 1$ to get
\begin{align*}
x_1  x_{23}^b &= (q_{12}q_{13})^{b-1}(b)_{\qti_{23}} x_{23}^{b-1} x_{123} + (q_{12}q_{13})^{b}x_{23}^{b} x_1,
\\
x_1 x_3 x_{23}^b &= (q_{12}q_{13}q_{32}q_{33})^b x_{23}^b x_{13} + q_{12}^{b-1}q_{13}^b (b)_{\qti_{23}} x_3x_{23}^{b-1}x_{123} + q_{12}^bq_{13}^{b+1} x_3 x_{23}^b x_1.
\\
x_1  x_{23}^b x_2=& -q_{12}^bq_{13}^{b-1}q_{32}(b)_{\qti_{23}} x_{23}^{b-1} x_2 x_{123} + (q_{12}q_{13})^{b}x_{23}^{b} x_{12} + (q_{12}q_{13})^{b}q_{12 }x_{23}^{b} x_2x_1,
\\
x_1 x_3 x_{23}^b x_2=& q_{12}^{b+1}(q_{13}q_{32}q_{33})^b q_{32} x_{23}^b x_2 x_{13} - q_{12}^{b}q_{13}^b q_{32}(b)_{\qti_{23}} x_3x_{23}^{b-1}x_2x_{123} \\
&+ q_{12}^bq_{13}^{b+1} x_3 x_{23}^b x_{12} + (q_{12}q_{13})^{b+1} x_3 x_{23}^b x_2x_1.
\end{align*}
These equalities prove that $x_{1}B \in I$.

\medbreak 
The following is an auxiliary tool to compute the Hilbert series $\htoba_{\bq}$.
\begin{stepo}
The subalgebra $\Zc$ of $\htoba_{\bq}$ generated by $x_{13}$ is a skew-central Hopf subalgebra isomorphic to a polynomial algebra in one variable. 
\end{stepo}
The generator of $\Zc$ is annihilated by the braided adjoint action of $x_i$, $i\in\I_3$. In fact, $(\ad_c x_i) x_{13} = 0$ hold in $\htoba_{\bq}$ either by definition for $i=2$ or by \eqref{eq:supera-J=I3-central-computations-xiij} if $i=1,3$. As $x_{13}$ is primitive, $\Zc$ is a Hopf subalgebra. Since $x_{13}$ is a non-zero primitive element with $\bq (\alpha_1+\alpha_3,\alpha_1+\alpha_3)=1$, it generates a polynomial algebra.

\begin{stepo}
There is a degree-preserving extension of braided Hopf algebras $\Zc \hookrightarrow \htoba_{\bq} \twoheadrightarrow \wtoba_{\bq}$.
\end{stepo}
Let $\pi\colon \htoba_{\bq} \twoheadrightarrow \wtoba_{\bq}$ denote the canonical projection, and put $\Zc'=\htoba_{\bq}^{\co \pi}$. Notice that $\Zc \subset \Zc'$ because $x_{13}$ projects to zero, and Lemma \ref{lem:extension-braided-graded} gives
\begin{align*}
\mathcal{H}_{\htoba_{\bq}} &= \mathcal{H}_{\wtoba_{\bq}}\mathcal{H}_{\Zc'} \ge \mathcal{H}_{\wtoba_{\bq}}\mathcal{H}_{\Zc} =
\frac{(1+t_3)(1+t_2)(1+t_1t_2t_3)(1+t_1)}{(1-t_2t_3)(1-t_1t_2)} \cdot \frac{1}{1-t_1t_3}.
\end{align*}
On the other hand, since $B$ linearly spans $\htoba_{\bq}$ we get the opposite inequality, so
\begin{align}\label{item:supera-J=I3-eminent-Hilbert}
\mathcal{H}_{\htoba_{\bq}} =
\frac{(1+t_3)(1+t_2)(1+t_1t_2t_3)(1+t_1)}{(1-t_2t_3)(1-t_1t_3)(1-t_1t_2)},
\end{align}
and this warranties $\Zc=\Zc'$. The desired extension follows.

\begin{stepo}
The set $B$ is a basis of $\htoba_{\bq}$, which has Gelfand-Kirillov dimension $3$.
\end{stepo}
Now this follows as in the proof of Proposition \ref{prop:eminent-supera-caseII}.
\epf

\begin{remark}
Let $\Zc_{\bq}$ denote the subalgebra of $\htoba_{\bq}$ generated by 
$x_{13}$, $x_{12}^N$ and $x_{23}^N$. One can verify that this subalgebra is skew-central (more precisely, it is annihilated by the braided adjoint action of $\htoba_{\bq}$), and that here is a degree-preserving extension of braided Hopf algebras $\Zc_{\bq} \hookrightarrow \htoba_{\bq} \twoheadrightarrow \toba_{\bq} $.
\end{remark}

\subsection{Type B}
Here $\bq$ is of type $\superqb{\theta}{q}{\J}$, where $q\in \G_{\infty}$ has order $N \neq 2,4$ and $\emptyset \ne \J \subseteq \I_{\theta}$. Up to relabeling we may assume that $\theta\notin\J$ (see \cite[\S 5.2]{AA17}), so Dynkin diagram is:
\begin{align*}
\xymatrix@C=40pt{\superqa{\theta-1}{q^2}{\J} \ar@{-}[r]^{\qquad q^{-2}} & \overset{q}{\circ}.}
\end{align*}

If $N\neq 3$, then the pre-Nichols algebra $\wtoba_{\bq}$ has the following presentation. 
\begin{align}\label{eq:Bsuper-preNichols-relsNneq3}
\begin{aligned}
x_{\theta\theta\theta\theta-1}&=0; & x_{ij} &= 0, \quad i < j - 1; &  x_{iii\pm1}&= 0, \quad q_{ii}\neq-1,i\in \I_{\theta-1};  \\
&& [x_{(i-1i+1)},x_i]_c&=0, \quad q_{ii}=-1; & x_i^2&=0, \quad q_{ii}=-1.
\end{aligned}
\end{align}
And if $N=3$, then $\wtoba_{\bq}$ is defined by the following relations:
\begin{align}\label{eq:Bsuper-preNichols-relsN=3}
\begin{aligned}
x_{ij} &= 0, \quad i < j - 1;&  x_{iii\pm1}&= 0, \quad q_{ii}\neq-1,i\in \I_{\theta-1};\\
x_i^2&=0, \quad q_{ii}=-1;& [x_{(i-1i+1)},x_i]_c&=0, \quad q_{ii}=-1;\\
x_{\theta\theta\theta\theta-1}&=0; & [x_{\theta\theta(\theta-1)},x_{\theta(\theta-1)}]_c&=0, \quad q_{\theta-1\theta-1}=-1;\\
&& [x_{\theta\theta(\theta-1)(\theta-2)},x_{\theta(\theta-1)}]_c&=0.
\end{aligned}
\end{align}
Let $\toba$ denote a pre-Nichols algebra of $\bq$ with $\GK\toba<\infty$.
We will check that each relation in \eqref{eq:Bsuper-preNichols-relsNneq3}, respectively \eqref{eq:Bsuper-preNichols-relsN=3}, must hold in $\toba$.

\begin{lemma}\label{lem:superb-xij}
Let $i,j\in \I_{\theta}$, $i<j-1$. Then $x_{ij}=0$ in $\toba$.
\end{lemma}
\pf
We fix first $j=\theta$. Then $x_{i\theta}=0$ by Lemma \ref{lem:ij}. 
Assume now that $j<\theta$. If $q_{ii}\ne -1$ or $q_{jj}\ne -1$, then $x_{ij}=0$ by Lemma \ref{lem:ij}. If $q_{ii}=q_{jj}=-1$, then $x_{ij}= 0$ by Lemma \ref{lem:ij3}, since $\qti_{jj+1}\neq 1=\qti_{ij+1}$.
\epf

\begin{lemma}\label{lem:superb-qs-mij=1}
If $i\in \I_{\theta-1}-\J$, then $x_{iii\pm 1}=0$ in $\toba$. 
\end{lemma}
\pf
We check first that $x_{\theta-1\theta-1\theta}=0$ if $i=\theta-1 \notin \J$. If $q_{\theta-1\theta-1}^{3}=1$, i.e. $q\in\G_3'$, then
\begin{align*}
q_{\theta-1\theta-1}^{c_{\theta-1 \theta}(1-c_{\theta-1 \theta})}q_{\theta\theta}^2
=q^{-4}q^2=q \neq 1.
\end{align*}
Hence either $q_{\theta-1\theta-1}^{2-c_{\theta-1\theta}}\neq 1$ or $q_{\theta-1\theta-1}^{c_{\theta-1 \theta}(1-c_{\theta-1 \theta})}q_{\theta\theta}^2\ne 1$. By Lemma \ref{lem:qs}, $x_{\theta -1 \theta -1 \theta}=0$.

Now we take $i\in \I_{\theta-2}-\J$. If either $q\notin\G_3'$ or $\qti_{i\, i+1}=-1$, then $x_{iii+1}=0$ by Lemma \ref{lem:qs}. Otherwise $q_{ii}=q_{i+1i+1}=\qti_{ii+1}^{-1}\in \G_3'$. As $\qti_{i+1i+2}\neq 1$ and $\qti_{ii+2}=1$, Lemma \ref{lem:dqs2} applies to prove that $x_{iii+1}=0$. The proof for $i\in \I_{2,\theta-1}-\J$, $x_{iii-1}\neq 0$ follows analogously.
\epf

\begin{lemma}\label{lem:superb-qs-mij=2}
$x_{\theta\theta\theta\theta-1}=0$ in $\toba$.
\end{lemma}
\pf
Here $c_{\theta\theta-1}=-2$. 
If $N=3$, then $x_{\theta\theta\theta\theta-1}=0$ by Lemma \ref{lem:QSR-orden-bajo}.
Now we assume $N \ne 3$. Notice that
\begin{align*}
q_{\theta\theta}^{2-c_{\theta\theta-1}}&=q^4,
&
q_{\theta\theta}^{c_{\theta\theta-1}(1-c_{\theta\theta-1})}q_{\theta-1\theta-1}^2 &=
\begin{cases}
q^{-6}(-1)^2=q^{-6}, & \text{si }\theta-1\in\J,
\\
q^{-6}q^4=q^{-2}, & \text{si }\theta-1\notin\J.
\end{cases}
\end{align*}
Hence either $q_{\theta\theta}^{2-c_{\theta\theta-1}}$ or $q_{\theta\theta}^{c_{\theta\theta-1}(1-c_{\theta\theta-1})}q_{\theta-1\theta-1}^2$. By Lemma \ref{lem:qs}, $x_{\theta\theta\theta\theta -1}=0$.
\epf

\begin{lemma}\label{lem:superb-[xijk,xk]}
If $i\in \J$, then $[x_{(i-1,i+1)},x_i]_{c}=0$ in $\toba$.
\end{lemma}
\pf
If either $q_{i-1i-1}=-1$ or $q_{i+1i+1}=-1$, then $x_u=0$ by Lemma \ref{lem:diagonal-[x_ijk,x_j]}.
Hence we assume that $q_{i-1i-1},q_{i+1i+1}\ne -1$.
If $i<\theta-1$, then $i+2\in \I_{\theta}$, $q_{i-1i-1}q_{i+1i+1}=1$, $\qti_{i+1i+2}\neq 1$, $\qti_{i-1i+2}=\qti_{ii+2}=1$ and $\qti_{i-1i}=\qti_{ii+1}^{-1}\neq \pm 1$, so $x_u=0$ by Lemma \ref{lem:cor2}

Finally we consider $i=\theta-1$, $\theta-2\notin\J$. Suppose that $x_u=0$. As $q_{uu}=q^{-1}$, $\qti_{u(\theta-2)}=1$, $\qti_{u(\theta-1)}=1$, $\qti_{\theta u}=q^{-2}\ne 1$, the diagram of $\bq'$ contains the following subdiagram:
\begin{align*}
\xymatrix@C=40pt{\overset{q^{-2}}{\underset{\theta-2}{\circ}} \ar@{-}[r]^{q^{2}} &  \overset{-1}{{\underset{\theta-1}\circ}} \ar@{-}[r]^{q^{-2}} & \overset{q}{{\underset{\theta}\circ}} \ar@{-}[r]^{q^{-2}}
&\overset{q^{-1}}{{\underset{u}\circ}}.}
\end{align*}
This diagram does not belong to \cite[Table 3]{H-classif} since
the two extremal vertices and one in the middle have label $\ne -1$, and these three labels are pairwise different, a contradiction with Conjecture \ref{conj:AAH}.
Hence $x_u=0$ in any case.
\epf

\begin{lemma}\label{lem:superb-[x_ttt-1,x_tt-1]}
If $N=3$ and $\theta -1\in \J$, then $[x_{\theta\theta(\theta-1)},x_{\theta(\theta-1)}]_c=0$ in $\toba$.
\end{lemma}
\pf
By Lemmas \ref{lem:xi-no-Cartan}, \ref{lem:superb-qs-mij=1}, \ref{lem:superb-qs-mij=2} and \cite[Lemma 5.9 (a)]{standard}, $[x_{\theta\theta(\theta-1)},x_{\theta(\theta-1)}]_c\in \Pc(\toba)$. Now Lemma \ref{lem:diagonal-[x_iij,x_ij]} applies and we get $[x_{\theta\theta(\theta-1)},x_{\theta(\theta-1)}]_c=0$.
\epf

\begin{lemma}\label{lem:superb-[x_ttt-1t-2,x_tt-1]}
If $N=3$, then $[x_{\theta\theta(\theta-1)(\theta-2)},x_{\theta(\theta-1)}]_c=0$ in $\toba$.
\end{lemma}
\pf
By Lemmas \ref{lem:xi-no-Cartan}, \ref{lem:superb-qs-mij=1}, \ref{lem:superb-qs-mij=2} and \cite[Lemma 5.9 (b)]{standard}, $[x_{\theta\theta(\theta-1)(\theta-2)},x_{\theta(\theta-1)}]_c \in \Pc(\toba)$. Now Lemma \ref{lem:diagonal-[x_iijk,x_ij]} applies, so $[x_{\theta\theta(\theta-1)(\theta-2)},x_{\theta(\theta-1)}]_c=0$. 
\epf

\begin{theorem}\label{thm:superb-distinguished-eminent}
Let $\bq$ be of type $\superqb{\theta}{q}{\J}$. Then the distinguished pre-Nichols algebra $\wtoba_{\bq}$ is eminent.
\end{theorem}
\pf
If $N\neq 3$, then the statement follows by Lemmas \ref{lem:xi-no-Cartan}, \ref{lem:superb-qs-mij=1}, \ref{lem:superb-qs-mij=2} and \ref{lem:superb-[xijk,xk]}. If $N=3$, then we apply the same results together with Lemmas \ref{lem:superb-[x_ttt-1,x_tt-1]} and \ref{lem:superb-[x_ttt-1t-2,x_tt-1]}.
\epf

\subsection{Type D}
Here $\bq$ is of type $\superqd{\theta}{q}{\J}$, where $\theta \geq 3$, $q\in\ku^\times$ is a root of unity of order $N> 2$,  and $\emptyset \ne \J \subseteq \I_{\theta}$, see \cite[\S 5.3]{AA17}. The possible Dynkin diagrams are
\begin{align}
\label{eq:dynkin-Dtheta-super-c}
 &\xymatrix@R-6pt{
{\bf A}_{\theta-1}(q;\J) \ar  @{-}[r]^{\hspace*{1.1cm} q^{-2}}  & \overset{q^2}{\underset{\ }{\circ}}},
\\
\label{eq:dynkin-Dtheta-super-d1}
&
\begin{aligned}
\xymatrix@R-6pt{ & \overset{-1}{\circ} \ar  @{-}[d]_{q^{-1}}\ar  @{-}[dr]^{q^2} & \\
{\bf A}_{\theta-2}(q;\J\cap \I_{\theta - 2}) \hspace*{-1.2cm} & \ar  @{-}[r]^{ q^{-1}}  & \overset{-1}{\underset{\ }{\circ}},}
\end{aligned}
& \theta-1 \in \J, 
\\ 
\label{eq:dynkin-Dtheta-super-d2}
& \begin{aligned}
\xymatrix@R-6pt{ & \overset{q^{-1}}{\circ} \ar @{-}[d]_{q} & \\
{\bf A}_{\theta-2}(q^{-1};\J\cap \I_{\theta - 2}) \hspace*{-1.1cm} & \ar  @{-}[r]^{q}  & \overset{q^{-1}}{\underset{\ }{\circ}},}
\end{aligned}
& \theta-1 \notin \J.
\end{align}

\begin{theorem}\label{thm:superd-distinguished-eminent}
Let $\bq$ be of type $\superqd{\theta}{q}{\J}$. Then the distinguished pre-Nichols algebra $\wtoba_{\bq}$ is eminent.
\end{theorem}

We give a proof for each diagram above in Propositions \ref{prop:superc-distinguished-eminent}, \ref{prop:superd1-distinguished-eminent} and \ref{prop:superd2-distinguished-eminent}.

\subsubsection{The diagram \eqref{eq:dynkin-Dtheta-super-c}}\label{subsec:dynkin-Dtheta-super-c}
The presentation of the distinguished pre-Nichols algebra $\wtoba_{\bq}$ depends on $N$ and $q_{\theta-1 \theta-1}$ as we describe below:

\begin{itemize}[leftmargin=*]
\item $q_{\theta-1 \theta-1}\neq -1$, $N>4$: the set of defining relations is
\begin{align}\label{eq:rels-type-D-super-a-N>4}
\begin{aligned}
&x_{ij}= 0, & &i < j - 1; & 
&x_{iii\pm1}=0,& &i \in\I_{\theta - 1}-\J ;  
\\
&x_{\theta-1\theta-1\theta-2}=0; & &&
&x_{\theta-1\theta-1\theta-1\theta}=0; &&
\\
&x_{\theta\theta\theta-1}=0; & &&
&[x_{(i-1i+1)},x_i]_c=0,& &i\in\I_{2,\theta - 2}\cap \J ;\\
&x_i^2=0, \quad i\in \J ; & && 
&&  &
\end{aligned}
\end{align}

\item  $q_{\theta-1 \theta-1}\neq -1$, $N=4$: the set of defining relations is
\begin{align}\label{eq:rels-type-D-super-a-N=4}
\begin{aligned}
& [x_{(\theta-2\theta)},x_{\theta-1\theta}]_c=0; &
& x_{iii\pm1}=0,& &i\in\I_{\theta - 2}- \J ;  \\
& x_{\theta-1\theta-1\theta-2}=0; & 
& x_{ij} = 0,& &i < j - 1; \\
& x_{\theta-1\theta-1\theta-1\theta}=0; & 
& [x_{(i-1i+1)},x_i]_c=0,& &i\in\I_{2,\theta - 2}\cap \J ;\\
& x_i^2=0, \quad i\in \J ; &  &x_{\theta\theta\theta-1}=0,&
& .
\end{aligned}
\end{align}

\item  $q_{\theta-1 \theta-1}\neq -1$, $N=3$: the set of defining relations is
\begin{align}\label{eq:rels-type-D-super-a-N=3}
\begin{aligned}
&[[x_{(\theta-2\theta)},x_{\theta-1}]_c,x_{\theta-1}]_c=0;  & 
& x_{iii\pm1}=0, & &i\in\I_{\theta - 2}- \J;  \\
& x_{\theta-1\theta-1\theta-2}=0; &
& x_{ij} = 0, & &i < j - 1;\\
& x_{\theta\theta\theta-1}=0; & 
& [x_{(i-1i+1)},x_i]_c=0, && i\in\I_{2,\theta - 2} \cap \J ;\\
& x_i^2=0, \quad i\in \J ; & 
& x_{\theta-1\theta-1\theta-1\theta}=0, & &
\end{aligned}
\end{align}

\item  $q_{\theta-1 \theta-1}=-1$, $N \neq 4$: the set of defining relations is
\begin{align}\label{eq:rels-type-D-super-b-Nneq4}
\begin{aligned}
&[[x_{\theta-2\theta-1},x_{(\theta-2\theta)}]_c,x_{\theta-1}]_c=0, && \theta-2 \in \J ; \\
&[[[x_{(\theta-3\theta)},x_{\theta-1}]_c,x_{\theta-2}]_c,x_{\theta-1}]_c=0, && \theta-2 \notin \J ;\\
&[x_{(i-1i+1)},x_i]_c=0, \quad i\in\I_{2,\theta - 2}\cap \J; &&x_{ij} = 0, \quad i < j - 1; \\
&x_{iii\pm1}=0, \quad i\in\I_{\theta - 2}-\J ;  && x_{\theta\theta\theta-1}=0; \\
& x_i^2=0, \quad i\in \J ; &&
.
\end{aligned}
\end{align}

\item  $q_{\theta-1 \theta-1}=-1$, $N=4$: the set of defining relations is
\begin{align}\label{eq:rels-type-D-super-b-N=4}
\begin{aligned}
& [[[x_{(\theta-3\theta)},x_{\theta-1}]_c, x_{\theta-2}]_c,x_{\theta-1}]_c=0, && \theta-2 \notin \J ;
\\
& [[x_{\theta-2\theta-1},x_{(\theta-2\theta)}]_c,x_{\theta-1}]_c=0, && \theta-2 \in \J ;  \\
& x_{iii\pm1}=0, \quad  i\in\I_{\theta - 2}- \J ;  && x_{ij} = 0, \ \ i < j - 1;  \\
&[x_{(i-1i+1)},x_i]_c=0, \quad i\in\I_{2,\theta - 2}\cap \J; && x_i^2=0, \ \ i\in \J; \\
& x_{\theta\theta\theta-1}=0; && x_{\theta-1\theta}^2=0.
\end{aligned}
\end{align}
\end{itemize}

Fix a finite $\GK$ pre-Nichols algebra $\toba$ of $\bq$.
We check that each relation in \eqref{eq:rels-type-D-super-a-N>4}, \eqref{eq:rels-type-D-super-a-N=4}, \eqref{eq:rels-type-D-super-a-N=3}, \eqref{eq:rels-type-D-super-b-Nneq4} and \eqref{eq:rels-type-D-super-b-N=4} holds in $\toba$.

\begin{lemma}\label{lem:superc-xij}
If $i,j\in \I_{\theta}$, $i<j-1$, then $x_{ij}=0$ in $\toba$.
\end{lemma}
\pf
Let $j=\theta$. Notice that $q_{ii}\in\{q^{\pm 1},-1\}$. Hence, 
\begin{itemize}[leftmargin=*]
\item if either $N\ne 4$ or $q_{ii}\ne -1$, then $x_{i\theta}=0$ by Lemma \ref{lem:ij};
\item if $N=4$, $q_{ii}=-1$, then $q_{ii}q_{\theta\theta}=1$, $\qti_{ii+1}\qti_{i+1\theta}=\pm q^{\pm 1} \ne 1$. Thus $x_{i\theta}=0$ by Lemma \ref{lem:ij3}.
\end{itemize}
Assume now that $j<\theta$. If $q_{ii}\ne -1$ or $q_{jj}\ne 1$, then $x_{ij}=0$ by Lemma \ref{lem:ij}. If $q_{ii}=q_{jj}=-1$, then $x_{ij}= 0$ by Lemma \ref{lem:ij3}, since $\qti_{jj+1}\neq 1$ and $\qti_{ij+1}=1$.
\epf

\begin{lemma}\label{lem:superc-x_ttt-1}
Te relation $x_{\theta\theta\theta-1}=0$ holds in $\toba$.
\end{lemma}
\pf
Recall that $c_{\theta\theta-1}=-1$. If $N\ne 4$, the desired relation follows from Lemma \ref{lem:qs} since $q_{\theta\theta}^{1-c_{\theta\theta-1}}=q^4\neq 1$ and $q_{\theta\theta}^{c_{\theta\theta-1}(1-c_{\theta\theta-1})}q_{\theta-1\theta-1}^2=q^{-4}q_{\theta-1\theta-1}^2\in \{q^{-4},q^{-2}\}$. 

In the case $N=4$ we can directly apply Lemma \ref{lem:diagonal:qs-m=1-q^m+1=1-extra-vertex}, since  $q_{\theta\theta}^{1-c_{\theta\theta-1}}=1$,
$q_{\theta\theta}=-1$, $\qti_{\theta-1\theta}=-1$,  
$\qti_{\theta-2\theta-1}=q^{\pm 1}\ne 1$, and  
$\qti_{\theta\theta-2}^2\qti_{\theta-1\theta-2}=\qti_{\theta-1\theta-2}\ne 1$. 
\epf

\begin{lemma}\label{lem:superc-x_t-1t-1t-1t}
If $q_{\theta-1\theta-1}\neq -1$ then $x_{\theta-1\theta-1\theta-1\theta}=0$ in $\toba$.
\end{lemma}
\pf
Now $c_{\theta-1\theta}=-2$. When $N\ne 3$ the relation follows from Lemma \ref{lem:qs}, since $q_{\theta-1\theta-1}^{1-c_{\theta-1\theta}}=q^3\ne 1$ and $q_{\theta-1\theta-1}^{c_{\theta-1\theta}(1-c_{\theta-1\theta})}q_{\theta\theta}^2=q^{-6}q^4=q^{-2}\neq 1$. 

If $N=3$, we are under the hypothesis of Lemma \ref{lem:QSR-orden-bajo} (indeed, $q_{\theta-1\theta-1}^{1-c_{\theta-1\theta}}=q^3=1$, $q_{\theta-1\theta-1}=q=q^{-2}=\qti_{\theta-1\theta}$, and  $q_{\theta\theta}\neq -1$),  so $x_{\theta-1\theta-1\theta-1\theta}=0$.
\epf

\begin{lemma}\label{lem:superc-qs-mij=1}
\begin{enumerate}[leftmargin=*,label=\rm{(\alph*)}]
\item\label{item:superc-qs-mij=1-1} If $i\in \I_{\theta-2}-\J$, then $x_{iii+ 1}=0$ in $\toba$.
\item\label{item:superc-qs-mij=1-2} If $i\in \I_{\theta-1}-\J$, then $x_{iii- 1}=0$ en $\toba$.
\end{enumerate}
\end{lemma}
\pf
\ref{item:superc-qs-mij=1-1} Notice that $c_{i\,i+1}=-1$ and $q_{ii}^2\ne 1$. If $x_{iii+1}\ne 0$ in $\toba$ then Lemma \ref{lem:qs} implies that $q_{ii}^3=1$ and $q_{ii}^{-2}q_{i+1i+1}^2=1$. So Lemma \ref{lem:dqs} assures that $q_{ii}=q_{i+1i+1}=\qti_{ii+1}^{-1}\in \G_3'$. This last equality together with $\qti_{ii+2}^2\qti_{i+1i+2}=\qti_{i+1i+2}\ne 1$ are exactly the hypothesis of Lemma \ref{lem:dqs2}, which states that $x_{iii+1}\ne 0$, a contradiction.

\noindent
\ref{item:superc-qs-mij=1-2} Note that $c_{i\,i-1}=-1$ and $q_{ii}^2\ne 1$. 
If $x_{iii-1}\ne 0$ in $\toba$ then Lemma \ref{lem:qs} implies that $q_{ii}^3=1$ and $q_{ii}^{-2}q_{i-1i-1}^2=1$. 
So Lemma \ref{lem:dqs} assures that $q_{ii}=q_{i-1i-1}=\qti_{ii-1}^{-1}\in \G_3'$. This and $\qti_{ii+1}^2\qti_{i-1i+1}=\qti_{ii+1}^2\ne 1$ allow us to apply Lemma \ref{lem:dqs2} to get $x_{iii-1}\ne 0$, a contradiction.
\epf

\begin{lemma}\label{lem:superc-[x_ijk,x_j]}
Given $i\in \I_{2,\theta-2}\cap \J$ we have $[x_{(i-1,i+1)},x_i]_c=0$ in $\toba$.
\end{lemma}
\pf
If $[x_{(i-1,i+1)},x_i]_c\neq 0$, Lemma \ref{lem:diagonal-[x_ijk,x_j]} implies $i-1$, $i+1\notin \J$. Now condition  \ref{item:supera-diagram-condition-i-notin-J} gives
 $q_{i-1i-1}=\qti_{i-1i}^{-1}=\qti_{ii+1}=q_{i+1i+1}^{-1}$.  Note also that  $\qti_{i+1i+2}\neq 1 =\qti_{ii+2}=\qti_{i-1i+2}$ so Lemma \ref{lem:cor2}\ref{item:cor2-c} assures $[x_{(i-1,i+1)},x_i]_c=0$, a contradiction.
\epf

\begin{prop}\label{prop:superc-distinguished-eminent}
If $\bq$ has Dynkin diagram (\ref{eq:dynkin-Dtheta-super-c}), then $\wtoba_{\bq}$ is eminent.
\end{prop} 
\pf 
Fix a pre-Nichols algebra $\toba$ with finite $\GK$. From Lemmas \ref{lem:xi-no-Cartan},\ref{lem:superc-xij}, \ref{lem:superc-x_ttt-1}, \ref{lem:superc-x_t-1t-1t-1t}, \ref{lem:superc-qs-mij=1} and \ref{lem:superc-[x_ijk,x_j]} we know that the quantum Serre relations as well as $x_i^2=0$ (for $i \in \J$) and  
$[x_{(j-1,j+1)},x_j]_c=0$ (for $j\in \I_{2,\theta-2}\cap \J$) hold in $\toba$. 
Other previous results apply:
\begin{itemize}[leftmargin=*]
\item If $N=4$ and $q_{\theta-1\theta-1}=-1$, then $\qti_{\theta-1\theta-2}^2\qti_{\theta\theta-2}^2=q^2\ne 1$. So $x_{\theta-1\theta}^2=0$ by Lemma \ref{lem:diagonal-xij^2}. 

\item If $N=4$ and $q_{\theta-1\theta-1}\neq -1$, then $i=\theta$, $j=\theta-1$, $k=\theta-2$ fulfill the hypothesis of Lemma \ref{lem:diagonal-[x_ij,x_ijk]}, so $[x_{\theta\theta-1},x_{\theta\theta-1\theta-2}]_c=0$. Thus  $[x_{(\theta-2,\theta)},x_{\theta-1\theta}]_c=0$.

\item If $N=3$ and $q_{\theta-1\theta-1}\neq -1$, we can apply Lemma \ref{lem:diagonal:[[xijk,xj]xj]} to $i=\theta-2$, $j=\theta-1$, $k=\theta$ and get $[[x_{(\theta-2,\theta)},x_{\theta-1}]_c,x_{\theta-1}]_c=0$.
\end{itemize}

It only remain to verify that the elements
\begin{align*}
&[[[x_{(\theta-3,\theta)},x_{\theta-1}]_c,x_{\theta-3}]_c,x_{\theta-1}]_c, &
&[[x_{\theta-2\theta-1},x_{(\theta-2,\theta)}]_c,x_{\theta-1}]_c,
\end{align*}
vanish in $\toba$ when $\bq$ fulfills the corresponding constraints.
In each case the defining relations of $\wtoba_\bq$ with smaller degrees certainly hold in $\toba$, so these elements are primitive in $\toba$. Hence they annihilate in $\toba$ thanks to Lemmas \ref{lem:diagonal-[[[x_ijkl,x_k],x_i],x_k]} and \ref{lem:diagonal-[[x_ij,x_ijk],x_j]}, respectively.
\epf

\subsubsection{The diagram \eqref{eq:dynkin-Dtheta-super-d1}} \label{subsec:dynkin-Dtheta-super-d1}
The presentation of $\wtoba_{\bq}$ depends on $\bq$ as follows:
\begin{itemize}[leftmargin=*]
\item $q_{\theta-2 \theta-2}\neq -1$,  $N\neq 4$: the set of defining relations is
\begin{align}\label{eq:rels-type-D-super-c-Nneq4}
\begin{aligned}
\begin{aligned}
& x_{\theta-2\, \theta-2 \, \theta}=0;&
& x_{iii\pm1}=0, \quad i\in\I_{\theta - 2}- \J ; \\
& x_{ij} = 0, \quad  i < j - 1,\theta-2;  &
& [x_{(i-1i+1)},x_i]_c=0, \quad i\in\I_{\theta - 3}\cap \J ;\\
& x_i^2=0, \quad  i\in \J ; & 
\end{aligned}\\
\raggedright \begin{aligned}
x_{(\theta-2\theta)} = q_{\theta-2\theta-1}(1-q^2)x_{\theta-1}x_{\theta-2\theta}
-q_{\theta-1\theta}(1+q^{-1})[x_{\theta-2\theta},x_{\theta-1}]_c.
\end{aligned}
\end{aligned}
\end{align}

\item $q_{\theta-2 \theta-2}\neq -1$, $N=4$: the set of defining relations is
\begin{align}\label{eq:rels-type-D-super-c-N=4}
\begin{aligned}
\begin{aligned}
& x_{iii\pm1}=0, \quad i\in\I_{\theta - 2}-\J; && [x_{(i-1i+1)},x_i]_c=0, \quad  i\in\I_{\theta - 3} \cap \J ; \\
& x_{\theta-2\, \theta-2 \, \theta}=0; &&x_{ij}= 0, \quad  i < j - 1,\theta-2;  \\
& x_i^2=0, \quad  i\in \J; && x_{\theta\theta\theta-1}=0;\\
&x_{\theta-1\theta-1\theta}=0;&&
\end{aligned}
\\
\raggedright
\begin{aligned}
x_{(\theta-2\theta)}= q_{\theta-2\theta-1}(1-q^2)x_{\theta-1}x_{\theta-2\theta} -q_{\theta-1\theta}(1+q^{-1})[x_{\theta-2\theta},x_{\theta-1}]_c.
\end{aligned}
\end{aligned}
\end{align}

\item $q_{\theta-2 \theta-2}=-1$, $N\neq 4$: the set of defining relations is
\begin{align}\label{eq:rels-type-D-super-d-Nneq4}
\begin{aligned}
\begin{aligned}
& [x_{\theta-3\, \theta-2 \, \theta},x_{\theta-2}]_c=0; &
& x_{iii\pm1}=0, &&  i\in\I_{\theta - 2} \cap \J ; \\
& x_{ij} = 0, \quad i < j - 1,\theta-2;  &
& [x_{(i-1i+1)},x_i]_c=0, && i\in \I_{\theta - 3} \cap \J;\\
& x_i^2=0, \quad i\in \J; &
\end{aligned}
\\
\raggedright
\begin{aligned}
x_{(\theta-2\theta)} = q_{\theta-2\theta-1}(1-q^2)x_{\theta-1}x_{\theta-2\theta}
-q_{\theta-1\theta}(1+q^{-1})[x_{\theta-2\theta},x_{\theta-1}]_c.
\end{aligned}
\end{aligned}
\end{align}

\item $q_{\theta-2 \theta-2}=-1$, $N=4$: the set of defining relations is
\begin{align}\label{eq:rels-type-D-super-d-N=4}
\begin{aligned}
\begin{aligned}
&x_{iii\pm1}=0, \quad i\in\I_{\theta - 2}-\J;  && x_{ij}= 0, \quad i < j - 1,\theta-2;  
\\
& [x_{\theta-3\, \theta-2 \, \theta},x_{\theta-2}]_c=0; &&
[x_{(i-1i+1)},x_i]_c=0, \quad i\in\I_{\theta - 3}\cap \J; 
\\
&x_i^2=0, \quad i\in\J;&& x_{\theta\theta\theta-1}=0;\\
&x_{\theta-1\theta-1\theta}=0;&&
\end{aligned}
\\
\begin{aligned}
x_{(\theta-2\theta)}= 2 q_{\theta-2\theta-1} x_{\theta-1}x_{\theta-2\theta} -q_{\theta-1\theta}(1+q^{-1}) [x_{\theta-2\theta},x_{\theta-1}]_c.
\end{aligned}
\end{aligned}
\end{align}
\end{itemize}

Next we verify that all relations (\ref{eq:rels-type-D-super-c-Nneq4}), (\ref{eq:rels-type-D-super-c-N=4}), (\ref{eq:rels-type-D-super-d-Nneq4}) and (\ref{eq:rels-type-D-super-d-N=4}) hold in any finite $\GK$ pre-Nichols algebra $\toba$ of $\bq$.

\begin{lemma}\label{lem:superd1-xij}
For $i,j\in \I_{\theta}$ with $i<j-1,\theta-2$, we have $x_{ij}=0$  in $\toba$.
\end{lemma}
\pf The argument depends on the value of $q_{ii}q_{jj}$ as follows:
\begin{itemize}[leftmargin=*]
\item If $q_{ii}\ne -1$ or $q_{jj}\ne -1$, then $x_{ij}=0$ by Lemma \ref{lem:ij}.
\item If $q_{ii}=q_{jj}=-1$, by Lemma \ref{lem:ij3} it is enough to find $l\ne i,j$ such that  $\qti_{il} \qti_{jl} \ne 1$. Now
\begin{itemize}
\item for $j=\theta$, note that $\qti_{i\theta-1} \qti_{\theta\theta-1}= q^2 \ne 1$;
\item if $j=\theta-1$ then $\qti_{i\theta} \qti_{\theta-1\theta}= q^2 \ne 1$;
\item when $j<\theta-1$ we have $\qti_{ij+1}\qti_{jj+1}=\qti_{jj+1}\neq 1$.
\end{itemize}
\end{itemize}
So $x_{ij}=0$ in any case.
\epf

\begin{lemma}\label{lem:superd1-qs-m=1}
\begin{enumerate}[leftmargin=*,label=\rm{(\alph*)}]
\item\label{item:superd1-qs-mij=1-1} If $i\in \I_{\theta-2}-\J$, then $x_{iii\pm 1}=0$ in $\toba$.
\item\label{item:superd1-qs-mij=1-2} If $\theta-2\notin \J$, then $x_{\theta-2\theta-2\theta}=0$ in $\toba$.
\end{enumerate}
\end{lemma}
\pf
\ref{item:superd1-qs-mij=1-1} Note that $c_{i\,i\pm1}=-1$ and $q_{ii}^2\ne 1$. Assume $x_{iii+1}\ne0$, so Lemma \ref{lem:qs} gives $q_{ii}^{3}= 1$ and $q_{ii}^{-2}q_{i+1i+1}^2=1$. 
Hence $q_{ii}=q_{i+1i+1} = \qti_{ii+1}^{-1}\in \G_3'$ by Lemma \ref{lem:dqs}. Using this and $\qti_{ii+2}\qti_{i+1i+2}=\qti_{i+1i+2}\neq 1$,  we apply Lemma \ref{lem:dqs2} and get $x_{iii+1}=0$, a contradiction.
The same argument leads to $x_{iii-1}=0$, using in the last step that $\qti_{ii+1}^2\qti_{i-1i+1}=\qti_{ii+1}^2\neq 1$.

\smallbreak
\noindent\ref{item:superd1-qs-mij=1-2} Now $c_{\theta -2 \theta}=-1$, $q_{\theta -2 \theta -2 }^2\ne 1$. Since $q_{\theta-2\theta-2}^{-2} q_{\theta\theta}^2  = q_{\theta-2\theta-2}^{-2}\ne 1$, Lemma \ref{lem:qs} applies.
\epf

\begin{lemma}\label{lem:superd1-qs-m=1-N=4}
If $N=4$, then $x_{\theta\theta\theta-1}= 0$ and $x_{\theta-1\theta-1\theta}=0$ in $\toba$.
\end{lemma}
\pf
Note that if $\{i,j\}=\{\theta-1, \theta\}$, then Lemma \ref{lem:diagonal:qs-m=1-q^m+1=1-extra-vertex} applies. Indeed, $c_{ij}=-1$, $q_{ii}^2=1$ and $q_{ii}=\qti_{ij}=-1$; moreover  $k=\theta-2$ is such that $\qti_{jk}=q^{-1}\ne 1$ and $\qti_{ik}^2\qti_{jk}=q^{-3}\ne 1$. 
\epf
\begin{lemma}\label{lem:superd1-[x_ijk,x_j]}
\begin{enumerate}[leftmargin=*,label=\rm{(\alph*)}]
\item\label{item:superd1-[xijk,xj]-1} If $i\in \I_{2,\theta-3}\cap \J$, then  $[x_{(i-1,i+1)},x_i]_c=0$ in $\toba$.
\item \label{item:superd1-[xijk,xj]-2} If $\theta-2\in \J$, then $[x_{\theta-3\theta-2\theta},x_{\theta-2}]_c=[x_{(\theta-3\theta-1)},x_{\theta-2}]_c=0$ in $\toba$.
\end{enumerate}
\end{lemma}
\pf
\ref{item:superd1-[xijk,xj]-1} Let $i\in \I_{\theta-3}\cap \J$. If either $i-1\in\J$ or $i+1\in \J$, then $[x_{(i-1,i+1)},x_i]_c=0$ by Lemma \ref{lem:diagonal-[x_ijk,x_j]}. Assume now $i-1,i+1\notin \J$, so  $\qti_{i+1i+2}\neq 1 =\qti_{ii+2}=\qti_{i-1i+2}$, and Lemma \ref{lem:cor2} \ref{item:cor2-c} gives $[x_{(i-1,i+1)},x_i]_c=0$.
\medbreak

\noindent \ref{item:superd1-[xijk,xj]-2} If $k\in\{\theta-1,\theta\}$, then $i=\theta-3$, $j=\theta-2$, $k$ fulfill the hypothesis of Lemma \ref{lem:diagonal-[x_ijk,x_j]}, since $q_{\theta -2 \theta -2 }=-1$, $\qti_{\theta-3 k} = 1$, $\qti_{\theta-3 \theta-2 }= \qti_{\theta-2 k}^{-1} = q \ne \pm 1$ and $q_{kk}=-1$. 
\epf

\begin{prop}\label{prop:superd1-distinguished-eminent}
If $\bq$ has Dynkin diagram \eqref{eq:dynkin-Dtheta-super-d1}, then $\wtoba_{\bq}$ is eminent.
\end{prop} 
\pf 
All defining relations of $\wtoba_{\bq}$ hold in any finite $\GK$ pre-Nichols algebra of $\bq$ by Lemmas \ref{lem:xi-no-Cartan}, \ref{lem:diagonal-x_ijk}, \ref{lem:superd1-xij}, \ref{lem:superd1-qs-m=1}, \ref{lem:superd1-qs-m=1-N=4} and \ref{lem:superd1-[x_ijk,x_j]}.
\epf

\subsubsection{The diagram \eqref{eq:dynkin-Dtheta-super-d2}} \label{subsec:dynkin-Dtheta-super-d2}
The following is a presentation of $\wtoba_{\bq}$:

\begin{itemize}[leftmargin=*]
\item $q_{\theta-2 \theta-2}\neq -1$, $N\neq 4$:
\begin{align}\label{eq:rels-type-D-super-e-Nneq4}
\begin{aligned}
& x_{iii\pm1}=0, \   i\in\I_{\theta - 2}-\J; && x_{\theta-2\, \theta-2 \, \theta}=0; & 
& x_{\theta-1\theta-1\theta-2}=0; \\
& x_{ij}=0, \  i < j - 1, \theta-2; && x_{\theta-1\theta}=0; & 
& x_{\theta\theta\theta-2}=0; \\
& [x_{(i-1i+1)},x_i]_c=0, \  i\in \I_{\theta - 3}\cap \J; &
& x_i^2=0, \  i\in \J . &
\end{aligned}
\end{align}

\item $q_{\theta-2 \theta-2}\neq -1$, $N=4$:
\begin{align}\label{eq:rels-type-D-super-e-N=4}
\begin{aligned}
&x_{\theta\theta\theta-2}=0; &
& x_{iii\pm1}=0, &&   i\in\I_{\theta - 2}-\J ; \\
&x_{\theta-2\, \theta-2 \, \theta}=0; & 
& [x_{(i-1i+1)},x_i]_c=0, && i\in \I_{\theta - 3} \cap \J; \\
&x_{\theta-1\theta-1\theta-2}=0; & 
&x_{ij}= 0, &&  i < j - 1,\theta-2; \\
&x_{\theta-1\theta}=0; &
&x_i^2=0, \quad i\in \J.  & &
\end{aligned}
\end{align}

\item $q_{\theta-2 \theta-2}=-1$, $N\neq 4$:
\begin{align}\label{eq:rels-type-D-super-f-Nneq4}
\begin{aligned}
&x_{iii\pm1}=0, \   i \in \I_{\theta - 2}-\J ; && x_{\theta\theta\theta-2} =0;  && x_{\theta-1\theta-1\theta-2}=0; 
\\
& [x_{(i-1i+1)},x_i]_c=0, \  i\in\I_{\theta - 2}\cap \J; & &x_{\theta-1\theta}=0; & 
& [x_{\theta-3\, \theta-2 \, \theta},x_{\theta-2}]_c=0;
\\
& x_{ij} = 0, \  i < j - 1, \theta-2;
&&
x_i^2=0, \   i\in \J. & 
\end{aligned}
\end{align}

\item $q_{\theta-2 \theta-2}=-1$, $N=4$:
\begin{align}\label{eq:rels-type-D-super-f-N=4}
\begin{aligned}
&x_{\theta-1\theta}=0; & &x_{ij} = 0,&& i < j - 1, \theta-2; 
\\
& [x_{\theta-3\, \theta-2 \, \theta},x_{\theta-2}]_c=0; &
& [x_{(i-1i+1)},x_i]_c=0, && i\in \I_{\theta - 2} \cap \J;
\\
& x_{\theta\theta\theta-2}=0; &
& x_{iii\pm1}=0, &&   i \in\I_{\theta - 2}-\J ; 
\\
& x_{\theta-1\theta-1\theta-2}=0; && x_i^2=0, \quad i\in \J && 
\end{aligned}
\end{align}
\end{itemize}

Given a finite $\GK$ pre-Nichols algebra $\toba$ of $\bq$, we show that  the relations in (\ref{eq:rels-type-D-super-e-Nneq4}), (\ref{eq:rels-type-D-super-e-N=4}), (\ref{eq:rels-type-D-super-f-Nneq4}) and (\ref{eq:rels-type-D-super-f-N=4}) hold in $\toba$.

\begin{lemma}\label{lem:superd2-xij}
For $i,j\in \I_{\theta}$ with $i<j-1,\theta-2$ we have $x_{ij}=0$ in $\toba$; also $x_{\theta -1 \theta} = 0$.
\end{lemma}
\pf  The relation $x_{\theta -1 \theta} = 0$ follows from Lemma \ref{lem:ij}, since $q_{\theta \theta }\ne 1$. 

Consider now $i+1<j\leq\theta-1$. Then the vertices $i, j$ belong to the subdiagram spanned by $\I_{\theta-1}$, which is of type $\superqa{\theta-1}{q^{-1}}{\J}$. Since $\theta-1\notin \J$, this last diagram is not of type $\superqa{3}{p}{\{1,2,3\}}$, so the relation $x_{ij}=0$ follows from Lemma \ref{lem:ijsupera}.

For the case $i+1<j=\theta$, we may apply the same argument but taking care of the subdiagram spanned by $\I_{\theta-2} \cup \{\theta\}$ instead of $\I_\theta$, and using that $\theta \notin \J$.
\epf

\begin{lemma}\label{lem:superd2-qs-m=1}
\begin{enumerate}[leftmargin=*,label=\rm{(\alph*)}]
\item\label{item:superd2-qs-1} If $i\in \I_{\theta-2}-\J$, then $x_{iii+1}=0$ in $\toba$.
\item\label{item:superd2-qs-2} If $i\in \I_{2,\theta-1}-\J$, then $x_{iii-1}=0$ in $\toba$.
\item\label{item:superd2-qs-3} If $\theta-2\notin \J$, then $x_{\theta-2\theta-2\theta}=0$ in $\toba$.
\item \label{item:superd2-qs-4} The relation $x_{\theta\theta\theta-2}=0$ holds in $\toba$.
\end{enumerate}
\end{lemma}
\pf
Note that the relations considered are of the form $x_{iij}$ for some $i,j$ with $c_{ij}=-1$.

\noindent\ref{item:superd2-qs-1} 
In the case $q_{ii}^{3}\ne 1$ or $q_{ii}^{-2}q_{i+1i+1}^2\ne 1 $, we have  $x_{iii+1}=0$ by Lemma \ref{lem:qs}. In the other case, Lemma \ref{lem:dqs} gives $q_{ii}=q_{i+1i+1}=\qti_{ii+1}^{-1}\in \G_3'$. By Lemma \ref{lem:dqs2}, the relation $x_{iii+1}$ will follow if we find $k\ne i, i+1$ such that $\qti^2_{ik}\qti_{i+1k} \ne 1$. This is easily achieved:
\begin{itemize} [leftmargin=*]
\item if $i<\theta-2$, then $\qti^2_{i i+2}\qti_{i+1i+2} =\qti_{i+1i+2} =q^{\pm}\ne 1$;
\item if $i=\theta-2$, then $\qti^2_{\theta-2 \theta}\qti_{\theta-1 \theta}= \qti^2_{\theta-2 \theta}=q^{-2} \ne 1$.
\end{itemize}
\noindent \ref{item:superd2-qs-2} Follows similarly.

\smallbreak
\noindent\ref{item:superd2-qs-3}\ref{item:superd2-qs-4} 
The subdiagram spanned by $\{\theta-2, \theta-1, \theta \}$ is of type $\superqa{3}{q^{-1}}{\J'}$ with $\J'\subset\{\theta -2\}$. 
Hence both claims follow from Lemma \ref{lem:qssupera}.
\epf

\begin{lemma}\label{lem:superd2-[x_ijk,x_j]}
\begin{enumerate}[leftmargin=*,label=\rm{(\alph*)}]
\item\label{item:superd2-[xijk,xj]-1} If $i\in \I_{2,\theta-3}\cap \J$, then $[x_{(i-1,i+1)},x_i]_c=0$ in $\toba$.
\item \label{item:superd2-[xijk,xj]-2} If $\theta-2\in \J$, then $[x_{(\theta-3\theta-1)},x_{\theta-2}]_c=0$ and  $[x_{\theta-3\theta-2\theta},x_{\theta-2}]_c=0$ in $\toba$.
\end{enumerate}
\end{lemma}

\pf
\ref{item:superd2-[xijk,xj]-1} Lemma \ref{lem:diagonal-[x_ijk,x_j]} assures the desired relation if either $i-1\in\J$ or $i+1\in \J$. If that is not the case, it follows that $q_{i-1i-1}q_{i+1i+1}=1$; also, since $i<\theta-2$, we have $\qti_{i+1i+2}\neq 1 =\qti_{ii+2}=\qti_{i-1i+2}$. Thus Lemma \ref{lem:cor2} \ref{item:cor2-c} gives $[x_{(i-1,i+1)},x_i]_c=0$.

\smallbreak
\noindent\ref{item:superd2-[xijk,xj]-2} By Lemma \ref{lem:diagonal-[x_ijk,x_j]}, we only need to consider the case $\theta-3\notin\J$. Since also $\theta-1\notin\J$ we have $q_{i-3i-3}q_{i-1i-1}=1$. Now Lemma \ref{lem:cor2} \ref{item:cor2-b} gives $[x_{(\theta-3\theta-1)},x_{\theta-2}]_c=0$, as  $\qti_{\theta-2\theta}^2=q^{-2}\ne 1=\qti_{\theta-3\theta}=\qti_{\theta-1\theta}$. 
A similar argument shows that $[x_{\theta-3\theta-2\theta},x_{\theta-2}]_c=0$.
\epf

\begin{prop}\label{prop:superd2-distinguished-eminent}
If $\bq$ has Dynkin diagram (\ref{eq:dynkin-Dtheta-super-d2}), then $\wtoba_{\bq}$ is eminent.
\end{prop} 
\pf 
Follows directly from Lemmas \ref{lem:xi-no-Cartan}, \ref{lem:superd2-xij}, \ref{lem:superd2-qs-m=1}, \ref{lem:superd2-[x_ijk,x_j]}.
\epf

\subsection{Type $\superda{\alpha}$.} The possible Dynkin diagrams of $\bq$ are:
\begin{multicols}{2}
\begin{align}\label{eq:dynkin-D2-1-alpha1}
\begin{aligned}
\xymatrix@C30pt@R-25pt{& & \\ \overset{q}{\underset{1 }{\circ}}\ar  @{-}[r]^{q^{-1}}  &
\overset{-1}{\underset{2 }{\circ}} \ar  @{-}[r]^{r^{-1}}  & \overset{r}{\underset{3}{\circ}},}
\end{aligned}
\\
\label{eq:dynkin-D2-1-alpha2}
\begin{aligned}
\xymatrix@C30pt@R-25pt{
&\overset{-1} {\underset{3}{\circ}} \ar  @{-}[dl]_{s} \ar  @{-}[dr]^{r} & \\
\overset{-1}{\underset{1}{\circ}} \ar  @{-}[rr]^{q}& & \overset{-1}{\underset{2}{\circ}},}
\end{aligned}
\end{align}
\end{multicols}
\noindent
where $q,r,s\in \ku-\{1\}$ are such that $qrs=1$. It is required that, if $q=-1$, then $r,s \neq -1$; when $q\ne-1$, either $r=1$ or $s=1$. See \cite[\S 5.4]{AA17} for details. 

In this subsection we prove the following.
\begin{theorem}\label{thm:D2-1-alpha-dist-is-eminent}
Let $\bq$ of type $\superda{\alpha}$. Then the distinguished pre-Nichols algebra  $\wtoba_\bq$ is eminent.
\end{theorem}

\subsubsection{The diagram \eqref{eq:dynkin-D2-1-alpha1}} \label{subsubsec:type-D2-1-alpha-a}
If $q,r,s\neq -1$, then $\wtoba_{\bq}$ is presented by
\begin{align}\label{eq:rels-type-D2-1-alpha-a-neq-1}
\begin{aligned}
x_2^2&=0; & x_{112}&=0; & x_{332}&=0; & x_{13}&=0.
\end{aligned}
\end{align}
If $q=-1$ and $r,s\neq -1$, then $\wtoba_{\bq}$ is presented by the relations
\begin{align}\label{eq:rels-type-D2-1-alpha-b-q=-1}
\begin{aligned}
x_{112}&=0; & x_2^2&=0; & x_{12}^2&=0; & x_{332}&=0; & x_{13}&=0.
\end{aligned}
\end{align}

\begin{prop}\label{prop:type-D2-1-alpha-a-distinguished-eminent}
If $\bq$ has Dynkin diagram \eqref{eq:dynkin-D2-1-alpha1}, then  $\wtoba_{\bq}$ is eminent.
\end{prop} 
\pf
Let $\toba$ a be pre-Nichols algebra of $\bq$ such that $\GK\toba<\infty$.

\smallbreak
\noindent 
\ref{eq:rels-type-D2-1-alpha-a-neq-1} Assume $q,r,s\neq -1$. By Lemma \ref{lem:xi-no-Cartan}, $x_2^2=0$ . As either $q_{11}\ne -1$ or $q_{33}\ne -1$, the relation $x_{13}=0$ follows from Lemma \ref{lem:ij}. For $x_{112}=0$ we apply Lemma \ref{lem:qs} \ref{item:lem-qs-2} since $c_{12}=-1$, $q_{11}^2=q^2\ne1$ and $q_{11}^{-2}q_{22}^2=q^{-2}\ne 1$. Similarly, using $r^{\pm2}\ne1$ we have $x_{332}=0$.

\smallbreak
\noindent
\ref{eq:rels-type-D2-1-alpha-b-q=-1}
Consider now $q=-1$. 
Again $x_2^2=0$ by Lemma \ref{lem:xi-no-Cartan}. 
For $x_{13}$, note that either $q_{11}\ne -1$ or $q_{33}\ne -1$, so Lemma \ref{lem:ij} applies. 
The argument given in the previous paragraph for $x_{332}=0$ also works here. 
For $x_{112}=0$ we apply Lemma \ref{lem:diagonal:qs-m=1-q^m+1=1-extra-vertex} since $c_{12}=-1$, $q_{11}^2=q^2=1$, $\qti_{12}=-1$, $\qti_{23}=r^{-1}\ne 1$ and $\qti_{13}^2\qti_{23}=\qti_{23}\ne 1$. 
Finally $x_{12}^2=0$ follows from Lemma \ref{lem:diagonal-xij^2} \ref{item:diagonal-xij^2-b} since $q_{11}=\qti_{12}=q_{22}=-1$  and $\qti_{13}^2\qti_{23}^2=r^{-2}\ne 1$.
\epf

\subsubsection{The diagram \emph{(\ref{eq:dynkin-D2-1-alpha2})}}
\label{subsubsec:type-D2-1-alpha-c}
Up to relabeling, we can assume that $\qti_{13},\qti_{23}\ne -1$. In this case the defining relations of the distinguished pre-Nichols algebra are
\begin{align}\label{eq:rels-type-D2-1-alpha-c} 
\begin{aligned}
&x_1^2=0; \qquad x_2^2=0; \qquad x_3^2=0; \qquad x_{12}^2=0 \text{ si }\qti_{12}=-1;
\\
&x_{(13)}-\frac{1-s}{q_{23}(1-r)}[x_{13},x_2]_c-q_{12}(1-s)x_2x_{13}=0.
\end{aligned}\end{align}

\begin{prop}\label{prop:d1alph3-distinguished-eminent}
If $\bq$ has Dynkin diagram \eqref{eq:dynkin-D2-1-alpha2}, then $\wtoba_{\bq}$ is eminent.
\end{prop} 
\pf
It follows by Lemmas \ref{lem:xi-no-Cartan}, \ref{lem:diagonal-xij^2} and \ref{lem:diagonal-x_ijk}.
\epf

\subsection{Type super $\superf$.} Fix $\bq$ of type $\superf$ \cite[\S 5.5]{AA17}: the Dynkin diagrams are
\begin{align}\label{eq:dynkin-F4-super}
\begin{aligned}
&\xymatrix@C-4pt{ \overset{\,\, q^2}{\underset{1}{\circ}}\ar  @{-}[r]^{q^{-2}}  & \overset{\,\,
q^2}{\underset{2}{\circ}} \ar  @{-}[r]^{q^{-2}}  & \overset{q}{\underset{3}{\circ}} &
\overset{-1}{\underset{4}{\circ}} \ar  @{-}[l]_{q^{-1}}} &
&\xymatrix@C-4pt{ \overset{\,\, q^2}{\underset{1}{\circ}}\ar  @{-}[r]^{q^{-2}}  & \overset{\,\,
q^2}{\underset{2}{\circ}} \ar  @{-}[r]^{q^{-2}}  & \overset{-1}{\underset{3}{\circ}} &
\overset{-1}{\underset{4}{\circ}} \ar  @{-}[l]_{q}}
\\
&\xymatrix@R-6pt@C-4pt{ &  \overset{q}{\underset{4}{\circ}}\ar  @{-}[d]_{q^{-1}} \ar
@{-}[rd]^{q^{-1}} &
\\ \overset{q^2}{\underset{1}{\circ}} \ar  @{-}[r]^{q^{-2}}  &  \overset{-1}{\underset{2}{\circ}} \ar  @{-}[r]^{q^{2}}  
& \overset{-1}{\underset{3}{\circ}}} &
&\xymatrix@R-6pt@C-4pt{  & \overset{-1}{\underset{4}{\circ}} \ar  @{-}[d]^{q^{2}} \ar  @{-}[ld]_{q^{-3}} &
\\ 
\overset{-1}{\underset{3}{\circ}} \ar  @{-}[r]^{q}  & \overset{-1}{\underset{2}{\circ}} & \overset{q^2}{\underset{1}{\circ}} \ar  @{-}[l]_{q^{-2}}  }
\\
&\xymatrix@C-4pt{ \overset{q^{2}}{\underset{1}{\circ}} \ar @{-}[r]^{q^{-2}}  
& \overset{q^2}{\underset{2}{\circ}} \ar  @{-}[r]^{q^{-2}}  
& \overset{-1}{\underset{3}{\circ}} \ar  @{-}[r]^{q^{3}} & \overset{\,\, q^{-3}}{\underset{4}{\circ}},}&
& \xymatrix@C-4pt{\overset{q^{2}}{\underset{1}{\circ}} \ar @{-}[r]^{q^{-2}}  
& \overset{q}{\underset{2}{\circ}} \ar  @{-}[r]^{q^{-1}}  
& \overset{-1}{\underset{3}{\circ}}\ar  @{-}[r]^{q^{3}} & \overset{\,\, q^{-3}}{\underset{4}{\circ}},}
\end{aligned}
\end{align}
where $q\in \G'_{N}$, with $N\geq 4$. Treating each of these diagrams separately, we prove:

\begin{theorem}\label{thm:F-super-dist-is-eminent}
Let $\bq$ of type $\superf$. Then the distinguished pre-Nichols algebra is eminent.
\end{theorem}

\subsubsection{The diagram \emph{(\ref{eq:dynkin-F4-super}
a)}} \label{subsubsec:type-F4-super-a}
When $N>4$ the defining relations of $\wtoba_{\bq}$ are:
\begin{align}\label{eq:rels-type-F4-super-a-Nneq4}
\begin{aligned}
x_{13}&=0; & x_{14}&=0; & x_{24}&=0; & x_{3332}&=0; & x_{4}^2&=0;\\
x_{112}&=0; & x_{221}&=0; & x_{223}&=0; & x_{334}&=0.
\end{aligned}
\end{align}

For $N=4$, the distinguished pre-Nichols algebra has the following presentation:
\begin{align}\label{eq:rels-type-F4-super-a-N=4}
\begin{aligned}
x_{13}&=0; & x_{14}&=0; & x_{221}&=0; & x_{3332}&=0; & [x_{(13)},x_2]_c&=0; & x_{4}^2&=0;\\
x_{334}&=0; & x_{24}&=0; & x_{112}&=0; & x_{223}&=0; & [x_{23},x_{(24)}]_c&=0.
\end{aligned}
\end{align}

\begin{prop}\label{prop:F4-super-a-eminent}
If $\bq$ has diagram (\ref{eq:dynkin-F4-super} a), then $\wtoba_{\bq}$ is eminent.
\end{prop} 
\pf
Let $\toba$ be a finite $\GK$ pre-Nichols algebra of $\bq$.
\begin{itemize}[leftmargin=*]
\item The subdiagram spanned by $\{1,2,3\}$ is of Cartan type $B_3$, hence all the defining relations involving only $x_1$, $x_2$ and $x_3$ hold in $\toba$ by Theorem \ref{thm:AndSanmarco}. 
\item Similarly, the subdiagram spanned by $\{4,3,2\}$ is of type $\superqd{3}{q}{\{4\}}$  \eqref{eq:dynkin-Dtheta-super-c} hence all the defining relations supported in $x_2,x_3,x_4$ hold in $\toba$ by Proposition \ref{prop:superc-distinguished-eminent}.
\end{itemize}
Finally $x_{14}=0$ holds in $\toba$ by Lemma \ref{lem:ij} since $q_{11}= q^2 \ne -1$.
\epf

\subsubsection{The diagram \emph{(\ref{eq:dynkin-F4-super} b)}} 
When $N>4$, a set of defining relations of the $\wtoba_{\bq}$ is
\begin{align}\label{eq:rels-type-F4-super-b-Nneq4}
\begin{aligned}
x_{13}&=0; & x_{14}&=0; & x_{24}&=0; & x_{112}&=0;& &[[x_{43},x_{432}]_c,x_3]_c=0;\\
x_{3}^2&=0; & x_{4}^2&=0; & x_{221}&=0; & x_{223}&=0.
\end{aligned}
\end{align}
In the case $N=4$, a presentation of $\wtoba_{\bq}$ is
\begin{align}\label{eq:rels-type-F4-super-b-N=4}
\begin{aligned}
x_{13}&=0; & x_{14}&=0; & x_{3}^2&=0; & x_{112}&=0; & [x_{(13)},x_2]_c&=0; & x_{223}&=0; \\
x_{24}&=0; & x_{23}^2&=0; & x_{4}^2&=0; & x_{221}&=0; & [[x_{43},x_{432}]_c,x_3]_c&=0. \\
\end{aligned}
\end{align}

\begin{prop}\label{prop:F4-super-b-eminent}
If $\bq$ has diagram (\ref{eq:dynkin-F4-super} b), then $\wtoba_{\bq}$ is eminent.
\end{prop} 
\pf
Let $\toba$ be a pre-Nichols algebra of $\bq$ with $\GK \toba<\infty$.
\begin{itemize}[leftmargin=*]
\item On one hand $\{2, 3, 4\}$ span a subdiagram of the form $\xymatrix@R-6pt{\overset{q^2}{\underset{\ }{\circ}} \ar  @{-}[r]^{\hspace*{-0.7cm} q^{-2}} & {\bf A}_{2}(q;\{3,4\}) }$ hence all the relations in \ref{eq:rels-type-F4-super-b-Nneq4} and \ref {eq:rels-type-F4-super-b-N=4} supported in $x_2,x_3, x_4$ hold in $\toba$ by Proposition \ref{prop:superc-distinguished-eminent}.

\smallbreak
\item If $N>4$ the subdiagram spanned by $\{1, 2, 3\}$ is of type $\superqa{3}{q^2}{\{3\}}$, hence all the defining relations involving only $x_1$, $x_2$ and $x_3$ hold in $\toba$ by Theorem \ref{thm:supera-distinguished-eminent}. 

\smallbreak
\item If $N=4$ the vertices $\{1, 2, 3\}$ span a subdiagram of Cartan type $A_3$ at $-1$, so $x_{112}=x_{223}=0$ by \cite[Lemma 5.3]{ASa}. We get $x_{13}=0$ from Lemma \ref{lem:ij3}, since $q_{11}q_{33}=1$ and $\qti_{14}\qti_{34}=\qti_{34}\ne 1$.
The element $x_{221}$ must vanish in $\toba$ since it is primitive and the Dynkin diagram of  $V\oplus \Bbbk x_{221}$ is
$$\xymatrix@C40pt@R-15pt{ 
 \overset{-1}{\underset{221}{\circ}} & & & \\
\overset{-1}{\underset{1}{\circ}}\ar  @{-}[r]^{-1\quad}  & \overset{\, -1}{\underset{2}{\circ}} \ar  @{-}[r]^{-1} \ar@{-}[ul]_{-1} & \overset{-1}{\underset{3}{\circ}} & \overset{-1}{\underset{4}{\circ}} \ar  @{-}[l]_{q},}$$
\noindent
which does not belong to \cite[Table 4]{H-classif} because $q\ne-1$. Finally turn to $x_u=[x_{(13)},x_{2}]_c$. From Lemma  \ref{lem:-supera-coproduct-[x123,x2]} and the facts $x_{13}=x_{223}=[x_{12},x_2]_c=0$, it follows that $x_u\in \Pc(\toba)$. But the Dynking diagram of $\ku x_u + \ku x_4$ is $\xymatrix{ \overset{1}{\circ}\ar  @{-}[r]^{q}  & \overset{-1}{\circ}}$, so Lemma \ref{lem:1connected} gives $[x_{(13)},x_{2}]_c=0$.
\end{itemize}

\smallbreak
Finally $x_{14}=0$ holds in $\toba$ since either $q_{11}^2 \ne 1$ (for $N>4$) so Lemma \ref{lem:ij} applies, or else $q_{11}q_{44}=1$,  $\qti_{13}\qti_{43}\neq 1$ (for $N=4$) and Lemma \ref{lem:ij3} applies here.
\epf

\subsubsection{The diagram \emph{(\ref{eq:dynkin-F4-super} c)}} \label{subsubsec:type-F4-super-c}
A set of defining relations of $\wtoba_\bq$ when $N>4$ is
\begin{align}\label{eq:rels-type-F4-super-c-Nneq4}
\begin{aligned}
x_{13}&=0; & x_{14}&=0; & x_{112}&=0; & &x_{(24)}- q_{34}q[x_{24},x_3]_c-q_{23}(1-q^{-1})x_3x_{24}=0;
\\
x_{442}&=0; & x_{443}&=0; & x_{2}^2&=0; & &x_{3}^2=0; \qquad [x_{(13)},x_2]_c=0.
\end{aligned}
\end{align}
A presentation when $N=4$ is
\begin{align}\label{eq:rels-type-F4-super-c-N=4}
\begin{aligned}
x_{13}&=0; & x_{14}&=0;   & x_{12}^2&=0; && [x_{(13)},x_2]_c=0; \quad x_{2}^2=0; \quad x_{3}^2=0;\\
x_{442}&=0; & x_{443}&=0; & x_{112}&=0; && x_{(24)}-q_{34}q[x_{24},x_3]_c-q_{23}(1-q^{-1})x_3x_{24}=0.
\end{aligned}
\end{align}

\begin{prop}\label{prop:F4-super-c-eminent}
If $\bq$ has diagram (\ref{eq:dynkin-F4-super} c), then $\wtoba_{\bq}$ is eminent.
\end{prop} 
\pf
Let $\toba$ be a pre-Nichols algebra of $\bq$, $\GK \toba<\infty$.
\begin{itemize}[leftmargin=*]
\item The vertices $\{2,3,4\}$ determine a subdiagram of type \eqref{eq:dynkin-Dtheta-super-d1} hence all the defining relations in  \eqref{eq:rels-type-F4-super-c-Nneq4}  and \eqref{eq:rels-type-F4-super-c-N=4} involving only $x_2$, $x_3$ and $x_4$ hold in $\toba$ by Proposition \ref{prop:superd1-distinguished-eminent}.

\smallbreak
\item If $N>4$  the vertices $\{1,2,3\}$ determine a subdiagram of type $\superqa{3}{q^2}{\{2,3\}}$, hence the relations in \ref{eq:rels-type-F4-super-c-Nneq4} involving only $x_1$, $x_2$ and $x_3$ hold in $\toba$ by Theorem \ref{thm:supera-distinguished-eminent}. 

\smallbreak
\item If $N=4$ we check case-by-case that all the defining relations  \ref{eq:rels-type-F4-super-c-N=4} supported in $1,2,3$ hold in $\toba$. Since $q_{11}q_{33}=1$ and $\qti_{14}\qti_{34}\neq 1$, Lemma \ref{lem:ij3} gives $x_{13}=0$, while $x_2^2=x_3^2=0$ hold by Lemma \ref{lem:xi-no-Cartan}. Also $x_{112}=0$ in $\toba$ by Lemma \ref{lem:diagonal:qs-m=1-q^m+1=1-extra-vertex} with $k=3$, and $x_{12}^2=0$ by Lemma \ref{lem:diagonal-xij^2} with $k=4$. Now focus on $x_u=[x_{(13)},x_2]_c$ which is primitive in $\toba$ by Lemma \ref{lem:-supera-coproduct-[x123,x2]} since $[x_{12},x_2]_c=x_{223}=x_{13}=0$. The Dynking diagram of $\ku x_u + \ku x_4$ is $\xymatrix{ \overset{1}{\circ}\ar  @{-}[r]^{q^{-3}}  & \overset{q}{\circ}}$, so Lemma \ref{lem:1connected} gives $[x_{(13)},x_{2}]_c=0$.
\end{itemize}
Finally $x_{14}=0$ holds in $\toba$ by Lemma \ref{lem:ij} since $q_{44}\ne -1$.
\epf

\subsubsection{The diagram \emph{(\ref{eq:dynkin-F4-super} d)}} \label{subsubsec:type-F4-super-d}
A presentation of $\wtoba_\bq$ when $N>4$ is
\begin{align}\label{eq:rels-type-F4-super-d-Nneq4}
\begin{aligned}
\begin{aligned}
x_{13}&=0; & x_{14}&=0; & x_{2}^2&=0; & &[x_{124},x_2]_c=0; \\
x_{112}&=0; & x_{4}^2&=0; & x_{3}^2&=0; & &[[x_{32},x_{321}]_c,x_2]_c=0; \\
\end{aligned}
\\
x_{(24)}+ q_{34}\frac{1-q^3}{1-q^2}[x_{24},x_3]_c-q_{23}(1-q^{-3})x_3x_{24}=0.
\end{aligned}\end{align}
When $N=4$ a set of defining relations  of $\wtoba_{\bq}$ is
\begin{align}\label{eq:rels-type-F4-super-d-N=4}
\begin{aligned}
&\begin{aligned}
x_{13}&=0; & x_{14}&=0; & x_{2}^2&=0; & & [x_{124},x_2]_c=0; \qquad x_{112}=0;\\
x_{12}^2&=0; & x_{4}^2&=0; & x_{3}^2&=0; & & [[x_{32},x_{321}]_c,x_2]_c=0;
\end{aligned}
\\
&x_{(24)}+ q_{34}\frac{1-q^3}{1-q^2}[x_{24},x_3]_c-q_{23}(1-q^{-3})x_3x_{24}=0.
\end{aligned}
\end{align}

\begin{prop}\label{prop:F4-super-d-eminent}
If $\bq$ has   Dynkin  diagram (\ref{eq:dynkin-F4-super} d), then $\wtoba_{\bq}$ is eminent.
\end{prop} 
\pf
Let $\toba$ be a pre-Nichols algebra of $\bq$ with $\GK \toba<\infty$.
\begin{itemize}[leftmargin=*]
\item The vertices $\{1,2,3 \}$ span a subdiagram of type \eqref{eq:dynkin-Dtheta-super-c}, hence the relations in \eqref{eq:rels-type-F4-super-d-Nneq4} and \eqref{eq:rels-type-F4-super-d-N=4} supported in $x_1,x_2,x_3$ hold in $\toba$ by Proposition \ref{prop:superc-distinguished-eminent}.

\smallbreak
\item The subdiagram determined by $\{2,3,4 \}$ is of type \eqref{eq:dynkin-D2-1-alpha2}, hence  the relations in \eqref{eq:rels-type-F4-super-d-Nneq4} and \eqref{eq:rels-type-F4-super-d-N=4} supported in $x_2,x_3,x_4$ hold in $\toba$ by Proposition \ref{prop:d1alph3-distinguished-eminent}.

\smallbreak
\item If $N>4$, the relations in \eqref{eq:rels-type-F4-super-d-Nneq4} involving only  $x_1$,$x_2$ and $x_4$ hold in $\toba$ by Theorem \ref{thm:supera-distinguished-eminent}, since $\{1,2,4\}$  determine a subdiagram of type $\superqa{3}{q^2}{\{2,4\}}$.

\smallbreak
\item If $N=4$, we check that each relation in \eqref{eq:rels-type-F4-super-d-N=4} with support contained in $x_1,x_2,x_4$  hold in $\toba$. 
Note that $x_{14}=0$ by Lemma \ref{lem:ij3} since $q_{11}q_{44}=1$ and $\qti_{13}\qti_{43}=q^{-3}\ne 1$, while $x_2^2=x_4^2=0$ by Lemma \ref{lem:xi-no-Cartan}.
Also, $x_{112}=0$ by \cite[Lemma 5.3]{ASa}, and $x_{12}^2=0$ by Lemma \ref{lem:diagonal-xij^2} with $k=3$. Let $x_u=[x_{124},x_2]_c$ which is primitive in $\toba$ by Lemma \ref{lem:-supera-coproduct-[x123,x2]} since $[x_{12},x_2]_c=x_{224}=x_{14}=0$. The Dynking diagram of $\ku x_u + \ku x_3$ is $\xymatrix{ \overset{1}{\circ}\ar  @{-}[r]^{q^{-1}}  & \overset{-1}{\circ}}$, so Lemma \ref{lem:1connected} gives $[x_{124},x_{2}]_c=0$.
\end{itemize}
Hence the relations in \eqref{eq:rels-type-F4-super-d-Nneq4} and \eqref{eq:rels-type-F4-super-d-N=4} hold in $\toba$.
\epf

\subsubsection{The diagram \emph{(\ref{eq:dynkin-F4-super} e)}} \label{subsubsec:type-F4-super-e}
A presentation of $\wtoba_\bq$ when $N\ne4,6$ is
\begin{align}\label{eq:rels-type-F4-super-e-Nneq4,6}
\begin{aligned}
x_{13}&=0; & x_{14}&=0; & &x_{24}=0; \quad x_{112}=0; \\
x_{221}&=0; & x_{223}&=0; &
& [[[x_{432},x_3]_c,[x_{4321},x_3]_c]_c,x_{32}]_c=0;\\
x_{443}&=0; & x_{3}^2&=0; & &
\end{aligned}
\end{align}
A presentation of $\wtoba_\bq$ when $N=6$ is
\begin{align}\label{eq:rels-type-F4-super-e-N=6}
\begin{aligned}
x_{13}&=0; & x_{14}&=0; & & x_{24}=0; \quad x_{112}=0;  \\
x_{221}&=0; & x_{223}&=0; &
& [[[x_{432},x_3]_c,[x_{4321},x_3]_c]_c,x_{32}]_c=0;\\
x_{34}^2&=0; & x_{3}^2&=0; & &
\end{aligned}
\end{align}
A set of defining relations of $\wtoba_\bq$ when $N=4$ is
\begin{align}\label{eq:rels-type-F4-super-e-N=4}
\begin{aligned}
x_{13}&=0; & x_{14}&=0; & x_{221}&=0: & x_{3}^2&=0; & & [x_{(13)},x_2]_c=0; \quad x_{112}=0;\\
x_{24}&=0; & x_{443}&=0; & x_{223}&=0; & x_{23}^2&=0; & & [[[x_{432},x_3]_c,[x_{4321},x_3]_c]_c,x_{32}]_c=0.
\end{aligned}
\end{align}

\begin{prop}\label{prop:F4-super-e-eminent}
If $\bq$ has Dynkin  diagram (\ref{eq:dynkin-F4-super} e), then  $\wtoba_{\bq}$ is eminent.
\end{prop} 
\pf
Let $\toba$ be a finite $\GK$ pre-Nichols algebra of $\bq$.
\begin{itemize}[leftmargin=*]
\item The subdiagram spanned by $\{2,3,4\}$ is of type \eqref{eq:dynkin-D2-1-alpha1}, hence those relations in \eqref{eq:rels-type-F4-super-e-Nneq4,6}, \eqref{eq:rels-type-F4-super-e-N=6} and \eqref{eq:rels-type-F4-super-e-N=4}
that do not involve $x_1$ hold in $\toba$ by Proposition \ref {prop:type-D2-1-alpha-a-distinguished-eminent}.

\smallbreak
\item If $N>4$ the vertices $\{1,2,3\}$ determine a subdiagram of type $\superqa{3}{q^2}{\{3\}}$, hence the relations \eqref{eq:rels-type-F4-super-e-Nneq4,6}, \eqref{eq:rels-type-F4-super-e-N=6} that do not involve $x_4$ hold in $\toba$ by Theorem \ref{thm:supera-distinguished-eminent}. 

\smallbreak
\item If $N=4$ we check that each relation in \eqref{eq:rels-type-F4-super-e-N=6} involving only $x_1$, $x_2$ and $x_3$ hold in $\toba$. Note that $x_{13}=0$ by Lemma \ref{lem:ij3} since $q_{11}q_{33}=1$ and $\qti_{14}\qti_{34}=q^3\ne 1$, while $x_3^2=0$ hold by Lemma \ref{lem:xi-no-Cartan}.
Also $x_{112}=x_{223}=0$  by \cite[Lemma 5.3]{ASa} and $x_{23}^2=0$ by Lemma \ref{lem:diagonal-xij^2} with $k=4$. 
If $x_{221} \ne 0$ we get $V\oplus \ku x_{221}\subset \Pc(\toba)$ with Dynkin diagram 
$$\xymatrix@C40pt@R-15pt{ 
\overset{-1}{\underset{221}{\circ}} & & & \\
\overset{-1}{\underset{1}{\circ}}\ar  @{-}[r]^{-1\quad}  & \overset{\, -1}{\underset{2}{\circ}} \ar  @{-}[r]^{-1} \ar@{-}[ul]_{-1} & \overset{-1}{\underset{3}{\circ}} & \overset{q^3}{\underset{4}{\circ}} \ar  @{-}[l]_{q^{-3}},}$$
\noindent
which is not in \cite[Table 4]{H-classif}. We get a contradiction with Conjecture \ref{conj:AAH}, so $x_{221}=0$.
Now $[x_{(13)},x_2]_c$ is primitive in $\toba$ by Lemma \ref{lem:-supera-coproduct-[x123,x2]}. Moreover, the diagram of $\ku[x_{(13)},x_2]_c+\ku x_4$ is  $\xymatrix{ \overset{1}{\circ}\ar  @{-}[r]^{q^{3}}  & \overset{q^{-3}}{\circ}}$, so Lemma \ref{lem:1connected} gives $[x_{(13)},x_{2}]_c=0$.
\end{itemize}
Finally $[[[x_{432},x_3]_c,[x_{4321},x_3]_c]_c,x_{32}]_c=0$ in $\toba$ by Lemma  \ref{lem:diagonal-[[[xijk,xj],[xijkl,xj]],xjk]} and 
$x_{14}=0$ holds in $\toba$ by Lemma \ref{lem:ij} since $q_{11}q_{44}=q^{-1}\neq 1$.
\epf

\subsubsection{The diagram \emph{(\ref{eq:dynkin-F4-super} f)}} \label{subsubsec:type-F4-super-f}
A presentation of $\wtoba_\bq$ when $N\ne4,6$ is
\begin{align}\label{eq:rels-type-F4-super-f-Nneq46}
\begin{aligned}
x_{13}&=0; & x_{14}&=0; & x_{24}&=0; && x_{443}=0; \qquad x_{3}^2=0;
\\
x_{112}&=0; & x_{2221}&=0; & x_{223}&=0; && [[x_{(14)},x_2]_c,x_3]_c=q_{23}(q^2-q)[[x_{(14)},x_3]_c,x_2]_c.
\end{aligned}
\end{align}
If  $N=6$, a set of defining relations of the distinguished pre-Nichols algebra is the following:
\begin{align}\label{eq:rels-type-F4-super-f-N=6}
\begin{aligned}
x_{13}&=0; & x_{14}&=0; & x_{24}&=0; && x_{34}^2=0; \qquad x_{3}^2=0; \qquad x_{443}=0;
\\
x_{112}&=0; & x_{2221}&=0; & x_{223}&=0;  && [[x_{(14)},x_2]_c,x_3]_c=q_{23}(q^2-q)[[x_{(14)},x_3]_c,x_2]_c.
\end{aligned}
\end{align}
When $N=4$ the relations defining $\wtoba_{\bq}$ are 
\begin{align}\label{eq:rels-type-F4-super-f-N=4}
\begin{aligned}
x_{13}&=0; & x_{14}&=0; & x_{24}&=0;  && x_{3}^2=0; \qquad  x_{112}=0; \qquad   [x_{12},x_{(13)}]_c=0;\\
x_{223}&=0; & x_{2221}&=0; & x_{443}&=0;  
&& [[x_{(14)},x_2]_c,x_3]_c=q_{23}(q^2-q)[[x_{(14)},x_3]_c,x_2]_c.
\end{aligned}
\end{align}

\begin{prop}\label{prop:F4-super-f-eminent}
If $\bq$ has Dynkin  diagram (\ref{eq:dynkin-F4-super} f),  then $\wtoba_{\bq}$ is eminent.
\end{prop} 
\pf
Let $\toba$ be a pre-Nichols algebra of $\bq$ with $\GK \toba<\infty$.
\begin{itemize}[leftmargin=*]
\item The subdiagram spanned by the vertices $\{1, 2, 3\}$ is of type $\superqd{3}{q}{\{3\}}$, hence the defining relations with support contained in $1,2,3$ hold in $\toba$ by Proposition \ref{prop:superd1-distinguished-eminent}. 
\item The vertices $\{2, 3, 4\}$ determine a subdiagram of type \eqref{eq:dynkin-D2-1-alpha1}, hence all the defining relations with support contained in $x_2,x_3,x_4$ hold in $\toba$ by Proposition \ref{prop:type-D2-1-alpha-a-distinguished-eminent}.
\end{itemize}
Finally $[[x_{(14)},x_2]_c,x_3]_c=q_{23}(q^2-q)[[x_{(14)},x_3]_c,x_2]_c$ in $\toba$ by Lemma \ref{lem:diagonal-[[[xijkl,xj]xk]} and 
$x_{14}=0$ holds in $\toba$ by Lemma \ref{lem:ij} since $q_{11}q_{44}=q^{-1}\neq 1$.
\epf

\subsection{Type  $\superg$}\label{subsec:type-G-super} Here $N=\ord q > 3$ and the possible Dynkin diagrams are
\begin{align}\label{eq:dynkin-G-super}
\xymatrix@C-6pt{ \overset{-1}{\underset{1}{\circ}}\ar  @{-}[r]^{q^{-1}}  &
\overset{q}{\underset{2}{\circ}} \ar  @{-}[r]^{q^{-3}}  & \overset{\,\, q^3}{\underset{3}{\circ}},} 
\quad
\xymatrix@C-6pt{ \overset{-1}{\underset{1}{\circ}} \ar  @{-}[r]^{q}  &
\overset{-1}{\underset{2}{\circ}} \ar  @{-}[r]^{q^{-3}}  & \overset{\,\, q^3}{\underset{3}{\circ}},}
\quad
\xymatrix@C-6pt{ \overset{-q^{-1}}{\underset{1}{\circ}} \ar  @{-}[r]^{q^2}  &
\overset{-1}{\underset{2}{\circ}} \ar  @{-}[r]^{q^{-3}}  & \overset{\,\, q^3}{\underset{3}{\circ}},}
\quad
\begin{aligned}
 \xymatrix@C-12pt@R-8pt{ & \overset{-1}{\underset{3}{\circ}} & \\
\overset{q}{\underset{1}{\circ}} \ar  @{-}[ru]^{q^{-2}} \ar@{-}[rr]^{q^{-1}}
& &\overset{-1}{\underset{2}{\circ}} \ar  @{-}[ul]_{q^{3}}}
\end{aligned}
\end{align}
see \cite[\S 5.6]{AA17}.
We deal with each of these diagrams and conclude the following.
\begin{theorem}\label{thm:G-super-dist-is-eminent}
Let $\bq$ of type $\superg$. Then the distinguished pre-Nichols algebra is eminent.
\end{theorem}

\subsubsection{The diagram \emph{(\ref{eq:dynkin-G-super} a)}} \label{subsubsec:type-G-super-a} 
A presentation of $\wtoba_\bq$ when $N\ne4,6$ is
\begin{align}\label{eq:rels-type-G-super-a-Nne4,6}
\begin{aligned}
x_{13}&=0; & x_{221}&=0; &  x_{332}&=0; & x_{22223}&=0; & x_{1}^2&=0.
\end{aligned}
\end{align}
For $N=6$ the defining relations of the distinguished pre-Nichols algebra are
\begin{align}\label{eq:rels-type-G-super-a-N=6}
\begin{aligned}
x_{13}&=0; & x_{221}&=0; & [x_{32},x_{(321)}]_c&=0; &x_{332}=0;&\\
x_{1}^2&=0; & x_{22223}&=0; & [[x_{223},x_{23}]_c,x_{23}]_c&=0. & &
\end{aligned}
\end{align}
If $N=4$, then a presentation of $\wtoba_\bq$ is
\begin{align}\label{eq:rels-type-G-super-a-N=4}
\begin{aligned}
x_{13}&=0; & x_{221}&=0; & [[[x_{(13)},x_2]_c,x_2]_c,x_2]_c&=0; & x_{22223}&=0;\\
x_{1}^2&=0; & x_{332}&=0; &  [x_{2223},x_{223}]_c &=0.
\end{aligned}
\end{align}

\begin{prop}\label{prop:G3-super-a-eminent}
If $\bq$ has Dynkin diagram (\ref{eq:dynkin-G-super} a), then $\wtoba_{\bq}$ is eminent.
\end{prop} 
\pf
Let $\toba$ be a pre-Nichols algebra of $\bq$ such that $\GK \toba<\infty$.
\begin{itemize}[leftmargin=*]
\item The vertices $\{2, 3\}$ determine a subdiagram of Cartan type $G_2$, hence those relations in \eqref{eq:rels-type-G-super-a-Nne4,6}, \eqref{eq:rels-type-G-super-a-N=6} and \eqref{eq:rels-type-G-super-a-N=4} involving only $x_2$ and $x_3$ hold in $\toba$ by Theorem \ref{th:CartanG2-eminent}. 
\smallbreak
\item The subdiagram spanned by $\{1, 2\}$ is of type super A, so the defining relations supported in $x_1,x_2$ hold in $\toba$ by Theorem \ref{thm:supera-distinguished-eminent} .
\end{itemize}
Consider  the relations involving $x_1$, $x_2$ and $x_3$. If $N=6$ we have $[x_{32},x_{321}]=0$ by Lemma \ref{lem:diagonal-[x_ij,x_ijk]}, and if $N=4$ then $[[[x_{(13)},x_2]_c,x_2]_c,x_2]_c=0$ by Lemma \ref{lem:diagonal-[[[xijk,xj],xj],xj]}. Finally, $x_{13}=0$ either by Lemma \ref{lem:ij} if $N\ne 6$ or by Lemma \ref{lem:ij3} if $N=6$.
\epf

\subsubsection{The diagram \emph{(\ref{eq:dynkin-G-super}
b)}} \label{subsubsec:type-G-super-b}
A set of defining relations of $\wtoba_\bq$ is
\begin{align}\label{eq:rels-type-G-super-b-Nneq6}
x_{1}^2&=0; & x_{2}^2&=0; & x_{13}&=0; & x_{332}&=0; & [[x_{12},&[x_{12},x_{(13)}]_c]_c,x_2]_c=0.
\end{align}

\begin{prop}\label{prop:G3-super-b-eminent}
If $\bq$ has Dynkin diagram (\ref{eq:dynkin-G-super} b), then  $\wtoba_{\bq}$ is eminent.
\end{prop} 
\pf
Given a finite $\GK$ pre-Nichols algebra $\toba$, Lemma \ref{lem:xi-no-Cartan} gives  $x_1^2=x_2^2=0$. Also,
\begin{itemize}[leftmargin=*]
\item If $N\neq 6$, then $x_{332}=0$ by Lemma \ref{lem:qs} and $x_{13}=0$ by Lemma \ref{lem:ij}.
\item If $N=6$, then $x_{332}=0$ by Lemma \ref{lem:diagonal:qs-m=1-q^m+1=1-extra-vertex} and $x_{13}=0$ by Lemma \ref{lem:ij3}.
\end{itemize}
Finally, the relations of $\toba_{\bq}$ of degree lower than that of $x_u=[[x_{12},[x_{12},x_{(13)}]_c]_c,x_2]_c$ are satisfied in $\toba$, so $x_u\in \Pc(\toba)$. Then $x_u=0$ by Lemma \ref{lem:diagonal-[[x_ij,[x_ij,x_ijk]],x_j]}.
\epf

\subsubsection{The diagram \emph{(\ref{eq:dynkin-G-super} c)}} \label{subsubsec:type-G-super-c}
A set of defining relations of $\wtoba_{\bq}$ when $N\ne 6$ is
\begin{align}\label{eq:rels-type-G-super-c-Nneq6}
\begin{aligned}
x_{13}&=0; \quad x_{332}=0; \quad x_{1112}=0; \quad x_{2}^2=0; \\
[x_1, & [x_{123},x_2]_c]_c = \frac{q_{12}q_{32}}{1+q}[x_{12},x_{123}]_c-(q^{-1}-q^{-2})q_{12}q_{13} x_{123}x_{12}.
\end{aligned}
\end{align}
If $N=6$ the distinguished pre-Nichols algebra is presented by 
\begin{align}\label{eq:rels-type-G-super-c-N=6}
\begin{aligned}
x_{13}&=0; \quad [x_{112},x_{12}]_c=0; \quad x_{2}^2=0; \quad x_{332}=0; \quad x_{1112}=0; \\
[x_1, & [x_{123},x_2]_c]_c = \frac{q_{12}q_{32}}{1+q}[x_{12},x_{123}]_c-(q^{-1}-q^{-2})q_{12}q_{13} x_{123}x_{12}.
\end{aligned}
\end{align}

\begin{prop}\label{prop:G3-super-c-eminent}
If $\bq$ has Dynkin diagram (\ref{eq:dynkin-G-super} c), then $\wtoba_{\bq}$ is eminent.
\end{prop} 
\pf
Let $\toba$ be a finite $\GK$ pre-Nichols algebra of $\bq$. Then $x_2^2=0$ by Lemma \ref{lem:xi-no-Cartan}, 
and $x_{13}=0$ by Lemma \ref{lem:ij} if $N\ne 4$ or by Lemma \ref{lem:ij} if $N=4$. Now
\begin{itemize}[leftmargin=*]
\item If $N\neq 6$, then $x_{332}=x_{1112}=0$ by Lemma \ref{lem:qs}.
\item If $N=6$, then $x_{332}=0$ by Lemma \ref{lem:diagonal:qs-m=1-q^m+1=1-extra-vertex}, $x_{1112}=0$ by Lemma \ref{lem:QSR-orden-bajo} and $[x_{112},x_{12}]_c=0$ by Lemma \ref{lem:diagonal-[x_iij,x_ij]}.
\end{itemize}
Finally, the defining relations of $\toba_{\bq}$ of degree lower than that of
$$ x_u=[x_1,[x_{123},x_2]_c]_c = \frac{q_{12}q_{32}}{1+q}[x_{12},x_{123}]_c-(q^{-1}-q^{-2})q_{12}q_{13} x_{123}x_{12}$$ 
are satisfied in $\toba$, so $x_u\in \Pc(\toba)$. Then $x_u=0$ by Lemma \ref{lem:diagonal-[x_i,[x_ijk,x_j]]}. 
\epf

\subsubsection{The  diagram \emph{(\ref{eq:dynkin-G-super} d)}} \label{subsubsec:type-G-super-d}
A presentation of $\wtoba_\bq$ is
\begin{align}\label{eq:rels-type-G-super-d}
\begin{aligned}
& x_{1112}=0; \quad x_{2}^2=0; \quad x_3^2=0; \quad x_{113}=0; \\
& x_{(13)}+q^{-2}q_{23}\frac{1-q^3}{1-q}[x_{13},x_2]_c-q_{12}(1-q^3)x_2x_{13}=0.
\end{aligned}
\end{align}

\begin{prop}\label{prop:G3-super-d-eminent}
If $\bq$ has Dynkin diagram (\ref{eq:dynkin-G-super},d), then $\wtoba_{\bq}$ is eminent.
\end{prop} 
\pf
This follows from Lemmas \ref{lem:xi-no-Cartan}, \ref{lem:qs} and \ref{lem:diagonal-x_ijk}.
\epf

\section{Eminent pre-Nichols algebras of standard type} \label{section:standard-type}

Here we consider pre-Nichols algebras of standard (non-super) type \cite[\S 6]{AA17}. The associated Cartan matrix is either $B_{\theta}$ or $G_2$, we study each case in a separate subsection. The main results can be summarized as follows:

\begin{theorem}\label{thm:standard-dist-is-eminent}
Let $\bq$ be a braiding matrix of standard type. Then the distinguished pre-Nichols algebra $\wtoba_{\bq}$ is eminent.
\end{theorem}

\subsection{Type $B_{\theta}$.} Here $\bq$ is a braiding of standard stype $B_{\theta}$. Up to relabeling $\I_\theta$ \cite[\S 6.1]{AA17}, we can assume that the Dynkin diagram is
\begin{align}\label{eq:dynkin-B-stadard}
\xymatrix{ {\bf A}_{\theta-1}(-\ztu;\J) \ar @{-}[r]^(.7){-\zeta }  & \overset{\zeta}{\underset{\ }{\circ}}},
\end{align}
where $\zeta \in \G_3'$. A set of defining relations of the distinguished pre-Nichols algebra is
\begin{align}\label{eq:rels-type-B-standard}
\begin{aligned}
& x_{ij}= 0, && i < j - 1; &  
& [x_{(i-1i+1)},x_i]_c=0, && i\in\J;  
\\
& x_{ii(i\pm1)}= 0, && i\in\I_{\theta-1}-\J; & 
& [x_{\theta\theta(\theta-1)(\theta-2)}, x_{\theta(\theta-1)}]_c=0;  
\\
& x_i^2=0, && i\in\J; & & [x_{\theta\theta(\theta-1)}, x_{\theta(\theta-1)}]_c=0, && \theta-1\in\J; 
\\
&x_\theta^3=0. & &&
&&&
\end{aligned}
\end{align}
In the following lemmas $\toba$ denotes  a finite $\GK$ pre-Nichols algebra of $\bq$.

\begin{lemma}\label{lem:superbstandard-qs-mij=0}
If $i,j\in \I_{\theta}$ are such that $i<j-1$, then $x_{ij}=0$ in $\toba$.
\end{lemma}
\pf
Assume first $i+1<j<\theta$. If $q_{ii}q_{jj}\ne 1$ we have $x_{ij}=0$ by Lemma \ref{lem:ij}. If $q_{ii}q_{jj}=1$, since $\qti_{ij+1}\qti_{jj+1}=\qti_{jj+1}\neq 1$, we get $x_{ij}= 0$ from Lemma \ref{lem:ij3}.

If $i\in \I_{\theta-2}$, since $q_{\theta\theta}\neq-1$, Lemma \ref{lem:ij} gives $x_{i\theta}=0$.
\epf

\begin{lemma}\label{lem:superbstandard-qs-mij=1}
If $i\in \I_{\theta-1}-\J$, then $x_{iii\pm 1}=0$ in $\toba$. 
\end{lemma}
\pf
This follows from Lemma \ref{lem:qs}, as $c_{ij}=-1$ and $\ord q_{ii}=6\neq 2, 3$.
\epf

\begin{lemma}\label{lem:superbstandard-[xijk,xj]}
If $i\in \J$, then $[x_{(i-1,i+1)},x_i]_{c}=0$ in $\toba$.
\end{lemma}
\pf
Notice that Lemmas \ref{lem:diagonal-[x_ijk,x_j]-primitive}, \ref{lem:superbstandard-qs-mij=0} and \ref{lem:superbstandard-qs-mij=1} assure $x_u:=[x_{(i-1,i+1)},x_i]_{c}\in \Pc(\toba)$. 
Suppose that $x_u\ne 0$. By Lemma \ref{lem:diagonal-[x_ijk,x_j]}we have $i\pm 1\notin\J$. 
Moreover, if $i<\theta-1$ then $i+2\in \I_{\theta}$, $q_{i-1i-1}q_{i+1i+1} =1$, $\qti_{i+1i+2}\neq 1$, $\qti_{i-1i+2}=\qti_{ii+2}=1$ and $\qti_{i-1i}=\qti_{ii+1}^{-1}\neq \pm 1$, which contradicts Lemma \ref{lem:cor2}. 
Thus $i=\theta-1$.
Since $\theta-2\notin\J$ we have $q_{\theta-2\theta-2}=-\ztu$. Thus $q_{uu}=-1$, $\qti_{u(\theta-2)}=1$, $\qti_{u(\theta-1)}=1$, $\qti_{\theta u}=q^{-2}\ne 1$.
Summarizing, $\Pc(\toba)$ contains $\Bbbk x_{\theta -2}\oplus\Bbbk x_{\theta -1}\oplus\Bbbk x_{\theta} \oplus\Bbbk x_u$ which has Dynkin diagram
\begin{align*}
\xymatrix@C=40pt{\overset{-\zeta}{\underset{\theta-2}{\circ}} \ar@{-}[r]^{-\ztu} &  \overset{-1}{{\underset{\theta-1}\circ}} \ar@{-}[r]^{-\zeta} & \overset{\zeta}{{\underset{\theta}\circ}} \ar@{-}[r]^{\zeta}
&\overset{\zeta^{2}}{{\underset{u}\circ}}.}
\end{align*}
This diagram is not in \cite[Table 3]{H-classif} since
\begin{itemize}[leftmargin=*]
\item there are three vertices with labels $\ne -1$, and
\item these vertices have pairwise different labels.
\end{itemize}
So $x_u=0$ in $\toba$.
\epf

\begin{lemma}\label{lem:superbstandard-[x_ttt-1,x_tt-1]}\label{lem:superbstandard-[x_ttt-1t-2,x_tt-1]}
\begin{enumerate}[leftmargin=*,label=\rm{(\alph*)}]
\item Si $\theta -1\in \J$, entonces $[x_{\theta\theta(\theta-1)},x_{\theta(\theta-1)}]_c=0$ en $\toba$.

\item $[x_{\theta\theta(\theta-1)(\theta-2)},x_{\theta(\theta-1)}]_c=0\text{ en }\toba$.
\end{enumerate}
\end{lemma}
\pf
From \ref{lem:xi-no-Cartan}, \ref{lem:superbstandard-qs-mij=0}, \ref{lem:superbstandard-qs-mij=1} and \cite[Lemma 5.9]{standard}, it follows that 
\begin{align*}
[x_{\theta\theta(\theta-1)(\theta-2)},x_{\theta(\theta-1)}]_c &\in \Pc(\toba), & [x_{\theta\theta(\theta-1)},x_{\theta(\theta-1)}]_c &\in \Pc(\toba)\text{ si }\theta -1\in \J.
\end{align*}
Thus the desired relations follow from Lemmas \ref{lem:diagonal-[x_iij,x_ij]} and \ref{lem:diagonal-[x_iijk,x_ij]}.
\epf

\begin{theorem}
Let $\bq$ be a braiding of standard type $B_{\theta}$. Then the distinguished pre-Nichols algebra is eminent.
\end{theorem} 

\pf
Follows from Lemmas \ref{lem:xi-no-Cartan}, \ref{lem:superbstandard-qs-mij=0}, \ref{lem:superbstandard-qs-mij=1}, \ref{lem:superbstandard-[xijk,xj]} and \ref{lem:superbstandard-[x_ttt-1t-2,x_tt-1]}.
\epf

\subsection{Type $G_2$} \label{subsection:standardG2}
There are three different possibilities for the Dynkin diagram of a braiding $\bq$ with root system of standard type $G_2$, namely 
\begin{align} \label{eq:dynkin-G2-standard}
\xymatrix @C=15pt{ \underset{ 1 }{\overset{q^2}{\circ}} \ar  @{-}[rr]^{q}  & & \underset{ 2 }{\overset{ q^{-1}}{\circ} } }, &&
\xymatrix @C=15pt{ \underset{ 1 }{\overset{q^2}{\circ}} \ar  @{-}[rr]^{q^3}  & & \underset{ 2 }{\overset{ -1}{\circ} } }, &&
\xymatrix @C=15pt{ \underset{ 1 }{\overset{q}{\circ}} \ar  @{-}[rr]^{q^5}  & & \underset{ 2 }{\overset{ -1}{\circ} } },
\end{align}
where $q$ is a root of unity of order $8$ \cite[\S 6.2]{AA17}. Next we give minimal presentations of these Nichols algebras and prove that the corresponding distinguished pre-Nichols algebras are eminent.  

\subsubsection{The generalized Dynkin diagram \emph{(\ref{eq:dynkin-G2-standard} a)}}\label{subsubsec:standardG2-a}

\begin{lemma}\label{lem:standardG2-a}
Assume $\bq$ is of type $\xymatrix @C=15pt{ {\overset{q^2}{\circ}} \ar  @{-}[rr]^{q}  & & {\overset{ q^{-1}}{\circ}}}$ with $q\in\G'_8$.
\begin{enumerate}[leftmargin=*,label=\rm{(\alph*)}]
\item\label{item:standardG2-a-Nichols} The Nichols algebra is minimally presented by generators $x_1, x_2$ and relations
\begin{align}\label{eq:standardG2-a-presentation}
x_{1}^4, && x_{221}, &&  [x_{1112}, x_{112}]_c, &&  x_{2}^8, && x_{112}^{8}.
\end{align}
\item\label{item:standardG2-a-eminent} The distinguished pre-Nichols algebra is eminent.
\end{enumerate}
\end{lemma}

\pf
\ref{item:standardG2-a-Nichols} By \cite[Theorem 3.1]{An-crelle},  $\toba_\bq$ is presented by the relations \eqref{eq:standardG2-a-presentation} and 
\begin{align*}
[x_1,x_{11212}]_c + q_{12}q^2 x_{112}^2, && [x_{11212},x_{12}]_c, && [x_{112}, x_{11212}]_c,
\end{align*}
We verify using \texttt{GAP} that this last bunch of relations can be deduced from $x_{1}^4$, $x_{221}$ and $[x_{1112}, x_{112}]_c$. 
For the minimality, it is enough to show that $[x_{1112}, x_{112}]_c$ does not belong to the ideal of $T(V)$ generated by $x_{1}^4$ and $x_{221}$, which is again checked with \texttt{GAP}.

\ref{item:standardG2-a-eminent} Notice first that  the previous computation provides a minimal presentation:
\begin{align*}
\wtoba_{\bq}=T(V)/\langle x_{1}^4, \ x_{221}, \ [x_{1112}, x_{112}]_c \rangle.
\end{align*}
If $\toba$ is a finite $\GK$ pre-Nichols of $\bq$, then $x_1^4=0$ and $x_{221} = 0$ by Lemmas \ref{lem:xi-no-Cartan} and \ref{lem:qs}, respectively. Thus $[x_{1112}, x_{112}]_c\in \Pc(\toba)$.
It must be $[x_{1112}, x_{112}]_c =0$ by Lemmas \ref{lem:subspace-primitives} and \ref{lem:1connected}, because the Dynkin diagram of $\ku [x_{1112}, x_{112}]_c+\ku x_1 \subset \Pc(\toba)$ is  $\xymatrix @C=15pt{ {\overset{1}{\circ}} \ar  @{-}[rr]^{q^6}  & & {\overset{ q^2}{\circ}}}$.
\epf

\subsubsection{The generalized Dynkin diagram \emph{(\ref{eq:dynkin-G2-standard} b)}}\label{subsubsec:standardG2-b}
By \cite[\S 6.2.6]{AA17} we have 
\begin{align*}
\wtoba_{\bq}=T(V)/\langle x_{1}^4, \ x_{2}^2, \ [x_{1}, x_{11212}]_c + q_{12} (1-q)^{-1} x_{112}^2 \rangle.
\end{align*}

\begin{lemma}\label{lem:standardG2-b}
If $\bq$ is of type $\xymatrix @C=15pt{ {\overset{q^2}{\circ}} \ar  @{-}[rr]^{q^3}  & & {\overset{ -1}{\circ}}}$ with $q\in\G'_8$, then $\wtoba_{\bq}$ is eminent.
\end{lemma}

\pf
If $\toba$ is a finite $\GK$ pre-Nichols, then $x_1^4=x_2^2=0$ in $\toba$ by Lemma \ref{lem:xi-no-Cartan}. Then $x_u:=[x_{1}, x_{11212}]_c + q_{12} (1-q)^{-1} x_{112}^2$ is primitive in $\toba$, and it must vanish because the Dynkin diagram of $\ku x_u+\ku x_1 \subset \Pc(\toba)$ is  $\xymatrix @C=15pt{ {\overset{1}{\circ}} \ar  @{-}[rr]^{q^6}  & & {\overset{ q^2}{\circ}}}$.
\epf

\subsubsection{The generalized Dynkin diagram \emph{(\ref{eq:dynkin-G2-standard} c)}}\label{subsubsec:standardG2-c}
By \cite[\S 6.2.7]{AA17} we have 
\begin{align*}
\wtoba_{\bq}=T(V)/\langle x_{11112}, \ x_{2}^2, \ [ x_{11212}, x_{12}]_c  \rangle.
\end{align*}

\begin{lemma}\label{lem:standardG2-c}
If $\bq$ is of type $\xymatrix @C=15pt{ {\overset{q}{\circ}} \ar  @{-}[rr]^{q^5}  & & {\overset{ -1}{\circ}}}$ with $q\in\G'_8$, then $\wtoba_{\bq}$ is eminent.
\end{lemma}

\pf
If $\toba$ is a finite $\GK$ pre-Nichols, then $x_{11112}=x_2^2=0$ in $\toba$ by Lemmas \ref{lem:qs} and \ref{lem:xi-no-Cartan}. Then $ [x_{11212}, x_{12}]_c $ is primitive in $\toba$, and it must vanish because the Dynkin diagram of $\ku  [ x_{11212}, x_{12}]_c +\ku x_1 \subset \Pc(\toba)$ is  $\xymatrix @C=15pt{ {\overset{1}{\circ}} \ar  @{-}[rr]^{q^7}  & & {\overset{ q}{\circ}}}$.
\epf

\end{document}